\documentclass[amstex, amysymb, 12pt]{report}
\usepackage{amsmath, amsfonts, amssymb, amsthm}

\newtheorem{Pa}{Paper}[chapter]
\newtheorem{Tm}[Pa]{{\bf Theorem}}
\newtheorem{La}[Pa]{{\bf Lemma}}
\newtheorem{Cy}[Pa]{{\bf Corollary}}
\newtheorem{Rk}[Pa]{{\bf Remark}}
\newtheorem{Pn}[Pa]{{\bf Proposition}}
\newtheorem{Pb}[Pa]{{\bf Problem}}
\newtheorem{Dn}[Pa]{{\bf Definition}}

\newcommand{\tA}{{\widetilde A}}
\newcommand{\C}{{\mathbb C}}
\newcommand{\tC}{{\widetilde C}}
\newcommand{\bC}{{\mathbf C}}
\newcommand{\btC}{{\widetilde{\mathbf C}}}
\newcommand{\tE}{{\widetilde E}}

\newcommand{\N}{{\mathbf N}}
\newcommand{\tN}{{\widetilde{\mathbf N}}}
\newcommand{\Sch}{{\mathbf S}}

\newcommand{\nline}{{\vspace{\baselineskip}}}

\newcommand{\w}{{w}}
\newcommand{\ka}{{\kappa}}

\def\IC{\hbox{\rm C\kern-.43em
       \vrule depth 0ex height 1.4ex width .05em\kern.41em}}

\def\IQ{\hbox{\rm Q\kern-.43em
       \vrule depth 0ex height 1.4ex width .05em\kern.41em}}

\begin{document}
\pagestyle{empty}
\centerline{BOUNDARY NEVANLINNA-PICK INTERPOLATION FOR}
\centerline{GENERALIZED
NEVANLINNA FUNCTIONS}
\vskip 0.5in
\centerline{Paul Anthony Smith}
~
\vskip 1.0in
\pagenumbering{roman}
\begin{center}{\large\bf Abstract}
\end{center}
\smallskip
We formulate three boundary
Nevanlinna-Pick interpolation problems for
generalized Nevanlinna functions.  For each problem,
we parameterize the set of all solutions in terms
of a linear fractional transformation with an
extended Nevanlinna class parameter.

\vskip .5in 

\noindent
{\bf Keywords}:  Generalized Nevanlinna function, fundamental matrix
inequality, boundary interpolation, interpolation with inequalities,
missed interpolation values, lost negative squares.

\tableofcontents
\pagestyle{plain}

\chapter{Introduction}
\pagenumbering{arabic}

\section{Notations}
Throughout this thesis we shall use the symbols
$\C$, $\C^+$, $\mathbb{R}$, $\mathbb{D}$,
and $\mathbb{T}$ to denote the field of complex numbers,
the open upper half plane, the field
of real numbers, the open unit disk, and
the unit circle, respectively.
We shall sometimes write $w^*$ rather than $\bar{w}$ for
$w \in \C$.  

If a statement involving some domain holds
regardless of whether the domain
is taken to be $\C^+$
or to be $\mathbb{D}$, then we shall
write just $\Delta$ for the domain, with
the understanding that $\Delta$ may consistently
be replaced by either $\C^+$ or by $\mathbb{D}$.
We denote the boundary of $\Delta$ by
$\partial \Delta$.
Note that $\partial \C^+ =
\mathbb{R}$ and $\partial \mathbb{D} = \mathbb{T}$.

If $f$ is a function meromorphic on
$\Delta$, then
$\rho(f)$ stands for its domain of analyticity.  We regard
$\rho(f)$ as a subset of $\Delta$, even if $f$ admits
an analytic continuation into a larger domain containing
$\Delta$.

Nontangential boundary limits shall be of fundamental 
importance throughout this thesis.  
To avoid repetitive phrasing, we adopt several conventions.
Let
$f$ be a function meromorphic on a domain $\Delta$ and
let $x_0 \in \partial \Delta$.  We denote by
\begin{equation}
f(x_0) := {\displaystyle \lim_{z \to x_0} f(z)},
\label{r1.1}
\end{equation}
the nontangential boundary limit where $z$ tends to $x_0$
staying in a nontangential neighborhood of $x_0$ in
$\Delta$.

A special case of (\ref{r1.1}) is of particular interest in this
thesis.
Let $g(z) \equiv f(z) / (z - x_0)$.
We then define $g_{-1}(x_0) := f(x_0)$ whenever the limit
(\ref{r1.1}) exists.
As a convention, we take $g_{-1}(x) \equiv \infty$
for $x \in \partial \Delta$ when $g \equiv \infty$.

Occasionally,
we shall denote by $\mathbb{N}_n$ the set of the first $n$ natural numbers: 
$\mathbb{N}_n := \{1, \ldots, n\}$.  We use $\mathbb{Z}_+$ to denote
the set of nonnegative integers.  The symbol $\mathbf{e}_i$
stands for the $i$-th column of the 
$n \times n$ identity matrix $I_n$
(the appropriate choice of $n$ will be clear from the context).
We use $\C^{a \times b}$ to denote the set
of $a \times b$ matrices with entries in $\C$.  The symbol
``$\vee$'' indicates the logical disjunction.  We use $\Im\;z$
to denote the imaginary part of $z \in \C$.

\section{The Classical Nevanlinna-Pick Interpolation Problem}
A function $\varphi$ belongs to the \emph{Nevanlinna class}
$\N$ (sometimes called the \emph{Pick class})
provided that $\varphi$ is analytic on $\C^+$ and has
nonnegative imaginary part there.
Interpolation theory for Nevanlinna functions originates in 
\cite{Pick} and
\cite{Nev}, where the following problem was posed and solved:
\begin{Pb}
Given
\begin{equation*}
z_1, \ldots, z_n \in \C^+
\quad\text{and}\quad
\varphi_1, \ldots, \varphi_n \in \C,
\end{equation*}
find all Nevanlinna functions $\varphi \in \N$
satisfying interpolation conditions
\begin{equation}
\varphi(z_i) = \varphi_i \quad\text{for}\quad i = 1, \ldots, n.
\label{r1.2}
\end{equation}
\label{r1.2.1}
\end{Pb}
The complete solution of the problem consists
of two parts: (1) necessary and sufficient conditions under which
at least one solution exists;
(2) a description of all solutions, provided that
the necessary conditions are met.
Necessary conditions for the existence of a solution
can be obtained
through using the Schwarz-Pick Lemma, which may
be formulated as follows:
\begin{La}
A function $\varphi$ defined on $\C^+$
belongs to the class $\N$ if and only if
the kernel
\begin{equation*}
K_\varphi(z, \zeta) := \frac{\varphi(z) - \varphi(\zeta)^*}{
z - \bar{\zeta}} \quad \quad (z, \zeta \in \C^+)
\end{equation*}
is positive on $\C^+ \times \C^+$: in notation,
$K_\varphi \succeq 0$.
\label{r1.2.2}
\end{La}
In more detail,
$K_\varphi \succeq 0$ means that for each choice
of an integer $n$ and of an $n$-tuple
$\mathbf{z} = (z_1, \ldots, z_n)$
of points in $\C^+$,
the Hermitian matrix
\begin{equation*}
P^\varphi(\mathbf{z}) :=
\left[K_\varphi(z_j, z_i)\right]_{i,j=1}^n 
= \left[ \frac{\varphi(z_j) - \varphi(z_i)^*}{z_j -
\bar{z}_i}\right]_{i,j=1}^n
\end{equation*}
is positive semidefinite.
\begin{Tm}
Problem \ref{r1.2.1} has a solution
if and only if the 
\emph{Pick matrix} 
\begin{equation*}
P := \left[ \frac{\varphi_j - \bar{\varphi}_i}{
z_j - \bar{z}_i} \right]_{i,j=1}^n
\end{equation*}
is positive semidefinite.
\end{Tm}
Suppose that $\varphi \in \N$ satisfies conditions (\ref{r1.2}).  Then
the Pick matrix $P$ is equal to $P^\varphi(\mathbf{z})$ with $\mathbf{z} = 
(z_1, \ldots, z_n)$, and by Lemma \ref{r1.2.2}, 
$P \geq 0$.
This establishes the necessity part of the Theorem.  Sufficiency was
established in \cite{Nev} using induction arguments.
 
In the case where $P$ is not positive semidefinite,
one may ask whether there exists
a function $\w$ that satisfies interpolation conditions
(\ref{r1.2}) and is  ``nearly'' a Nevanlinna
function in some sense.  This generalization of the problem
calls for a suitable extension of the class
of Nevanlinna functions.

\section{Interpolation For Generalized Nevanlinna Functions}
\begin{Dn}
Let $\ka$ be a nonnegative integer.
A function $\w$ belongs to the \emph{generalized Nevanlinna
class} $\N_\ka$ provided that $\w$ is meromorphic on 
$\C^+$ and is such that the kernel
\begin{equation*}
K_\w(z, \zeta) := \frac{\w(z) - \w(\zeta)^*}{
z - \bar{\zeta}} \quad \quad (z,\zeta \in \rho(\w))
\end{equation*}
has $\ka$ negative squares on
$\rho(\w) \times \rho(\w)$: in notation,
\begin{equation}
{\rm sq}_-K_\w = \ka.
\label{1.8}
\end{equation}
\end{Dn}
More explicitly, condition (\ref{1.8})
means that for every choice of an integer $n$ and
of an $n$-tuple $\mathbf{z} = (z_1, \ldots, z_n)$ of points in
$\rho(\w)$, the Hermitian
matrix
\begin{equation}
P^\w(\mathbf{z}) := \left[K_\w(z_j,z_i)\right]_{i,j=1}^n =
\left[\frac{\w(z_j) - \w(z_i)^*}{z_j - \bar{z_i}}
\right]_{i,j=1}^n
\label{1.9}
\end{equation}
has at most $\ka$ negative eigenvalues,
\begin{equation*}
{\rm sq}_-\left[\frac{\w(z_j) - \w(z_i)^*}{z_j - \bar{z_i}}
\right]_{i,j=1}^n \leq \ka,
\end{equation*}
and has exactly $\ka$ negative eigenvalues
for at least one choice of $n$ and $n$ points
$z_1, \ldots, z_n \in \rho(\w)$.
Note that the class $\N_0$
coincides with the class $\N$
(more precisely, if $K_\w \succeq 0$
on a domain with nonempty interior,
then $\w$ admits
an (unique) analytic continuation to all of $\C^+$).
It will sometimes be convenient to consider an extended
class of functions constructed by adjoining
$\w \equiv \infty$ to the class $\N_\ka$:
\begin{equation*}
\tN_\ka := \N_\ka \cup \{\infty\}.
\end{equation*}
In the setting of generalized Nevanlinna functions, 
the Nevanlinna-Pick interpolation problem becomes:
\begin{Pb}
Given
\begin{equation*}
\ka \in \mathbb{Z}_+, \quad
z_1, \ldots, z_n \in \C^+,
\quad\text{and}\quad
\w_1, \ldots, \w_n \in \C,
\end{equation*}
find all functions $\w \in \N_\ka$ such that
\begin{equation*}
\w(z_i) = \w_i \quad\text{for}\quad i = 1, \ldots, n.
\end{equation*}
\label{aboveProb}
\end{Pb}

\section{Boundary Nevanlinna-Pick Interpolation}

Problem \ref{aboveProb} was studied in \cite{ADL2},
where a parameterization of the solution set was obtained.
A description of excluded parameters (i.e., parameters
that do not lead to a solution) was provided; however,
this description was not explicit.

In this thesis we study the boundary analogues of
Problem \ref{aboveProb},
where the interpolation nodes are taken to be on the boundary of
$\C^+$ (that is, on $\mathbb{R}$), and the prescribed
values of the interpolant $\w$ are replaced by the
prescribed values of its nontangential boundary limits.
The first difficulty arising in this setting is that the
nontangential boundary limits are only guaranteed to exist
at almost every boundary point (but not at every boundary point).  
Thus we first are interested
in conditions that guarantee the existence of the nontangential
boundary limits at the interpolation nodes.
Such conditions are established in the following two versions of the 
Carath\'eodory-Julia theorem for
generalized Nevanlinna functions.  
The proofs will be presented in Chapter III.
\begin{Tm}
Let $\ka \in \mathbb{Z}_+$,
$\w \in \N_\ka$, and $x_0 \in \mathbb{R}$.  Then the following
statements are equivalent:
\begin{enumerate}
\item{\begin{equation*}
\liminf_{z \to x_0} |\w(z)| < \infty
\quad\text{and}\quad
-\infty < d := \liminf_{z \to x_0}K_\w(z,z) < \infty.
\end{equation*}}
\item{The nontangential limits
\begin{equation*}
\w_0 := \lim_{z \to x_0}\w(z)
\quad\text{and}\quad
\tilde{d} := \lim_{z \to x_0}K_\w(z,z)
\end{equation*}
exist and are real.}
\item{The nontangential limits
\begin{equation*}
\w_0 := \lim_{z \to x_0}\w(z)
\quad\text{and}\quad
\w_1 := \lim_{z \to x_0}\w^\prime(z)
\end{equation*}
exist and are real.}
\item{The nontangential limits
\begin{equation*}
\w_0 := \lim_{z \to x_0}\w(z)
\quad\text{and}\quad
\tilde{\w}_1 := \lim_{z \to x_0}
\frac{\w(z) - \w_0}{z - x_0}
\end{equation*}
exist and are real.}
\end{enumerate}
In this case, $d = \tilde{d} = \w_1 = \tilde{\w}_1$.
\label{Tm1}
\end{Tm}
When the nontangential limit ${\displaystyle \lim_{z \to x_0}K_\w(z,z)}$ exists,
we sometimes choose to abbreviate it to 
$K_\w(x_0,x_0)$.
Theorem \ref{Tm1} is supplemented by the following result.
\begin{Tm}
Let $\ka \in \mathbb{Z}_+$,
$\w \in \N_\ka$ and $x_0 \in \mathbb{R}$.  Suppose that
the limit
\begin{equation*}
\lim_{z \to x_i}|\w(z)|
\end{equation*}
exists and equals infinity.  Then the following
statements are equivalent:
\begin{enumerate}
\item{\begin{equation*}
d := \liminf_{z \to x_0}K_{-1/\w}(z,z) \in \mathbb{R} \setminus
\{0\}.
\end{equation*}}
\item{The nontangential limit
${\displaystyle \tilde{d} := \lim_{z \to x_0}K_{-1/\w}(z,z)}$ exists
and belongs to $\mathbb{R} \setminus \{0\}$.}
\item{The nontangential limit
\begin{equation*}
\w_{-1} := \lim_{z \to x_0}(z-x_i)\w(z)
\end{equation*}
exists and belongs to $\mathbb{R} \setminus \{0\}$.}
\item{The nontangential limit
\begin{equation*}
\tilde{\w}_{-1} := -\lim_{z \to x_0}
(z - x_i)^2\w^\prime(z)
\end{equation*}
exists and belongs to $\mathbb{R} \setminus \{0\}$.}
\end{enumerate}
In this case, $d = \tilde{d} = \w_{-1} = \tilde{\w}_{-1}$.
\label{Tm2}
\end{Tm}
Theorems \ref{Tm1} and \ref{Tm2} suggest boundary interpolation
problems with the data set
\begin{equation}
\Omega = \left\{ \{x_i\}_{i = 1}^n, \{\w_j\}_{j = 1}^\ell,
\{\gamma_j\}_{j = 1}^\ell, \{\xi_k\}_{k = \ell + 1}^n \right\}, 
\quad 0 \leq \ell \leq n,
\label{1.18}
\end{equation}
where $x_i \in \mathbb{R}$ are the interpolation nodes,
$\w_j, \gamma_j$ are real numbers, and $\xi_k$ are nonzero
real numbers .
With the data set $\Omega$ we associate a certain matrix (the Pick matrix)
$P$ that will play an important
role in the subsequent analysis.  
For the sake of clarity, we display the Pick matrix 
in block matrix form and define
each block individually:
\begin{equation}
P = \begin{bmatrix}P_{11} & P_{12} \\ P_{21} & P_{22}\end{bmatrix}.
\label{1.19}
\end{equation}
The block entry $P_{11}$ is associated only to the data
 \begin{equation*}
\left\{\{x_i\}_{i = 1}^\ell, \{\w_i\}_{i = 1}^\ell,
\{\gamma_i\}_{i = 1}^\ell\right\}
\end{equation*}
and is given by:
\begin{equation}
P_{11} = [p_{ij}]_{i,j=1}^\ell \qquad p_{ij} = \left\{
\begin{array}{ll}
\frac{\w_j - \w_i}{x_j - x_i} & \text{ for } i\neq j\\
\gamma_i &\text{ for } i = j.
\end{array}\right.
\label{117}
\end{equation}
The block entry $P_{22}$ is the diagonal matrix:
\begin{equation}
P_{22} = [p_{ij}]_{i,j=\ell+1}^n = 
\begin{bmatrix}-\xi_{\ell+1} & & 0 \\ & \ddots & \\ 0 & & -\xi_n\end{bmatrix}.
\label{118}
\end{equation}
We define the block entry $P_{12}$ as follows,
\begin{equation}
P_{12} = [p_{ij}]_{i = 1, \ldots, \ell}^{j = \ell+1, \ldots, n}
\quad\text{where}\quad 
p_{ij} = \frac{\xi_j}{x_j - x_i},
\label{119}
\end{equation}
and take $P_{21}$ to be the adjoint of $P_{12}$:
\begin{equation}
P_{21} = \left( P_{12} \right)^*.
\label{120}
\end{equation}
Note that the Pick matrix $P$ is Hermitian and therefore has only real
eigenvalues. 
The number of negative eigenvalues,
\begin{equation}
\ka := {\rm sq}_-P,
\label{1.25}
\end{equation}
will also play an important role in the subsequent analysis.

\section{The Three Main Problems}
We now introduce the interpolation problems that are the main object
of our study.

\begin{Pb}
Given the data $\Omega$ as in (\ref{1.18}) and the integer
$\ka$ defined by (\ref{1.25}), 
find all functions $\w \in \N_\ka$ such that
\begin{align}
\w(x_i)&:=\lim_{z \to x_i}\w(z) = \w_i &(i = 1, \ldots, \ell), \label{1.26}\\
\w^\prime(x_i)&:=\lim_{z \to x_i}\w^\prime(z) = \gamma_i &(i = 1,
\ldots, \ell), \label{1.26b}\\
\w_{-1}(x_j)&:=\lim_{z \to x_j}(z - x_j)\w(z) = \xi_j &(j = \ell+1,
\ldots, n), \label{1.27c}
\end{align}
where all of the limits are nontangential.
\label{Prob1}
\end{Pb}
We shall sometimes call the interpolation nodes
$x_i$ with $i \in \{1, \ldots, \ell\}$ {\it regular interpolation
  nodes}, 
and similarly we shall sometimes call the interpolation nodes
$x_j$ with $j \in \{\ell + 1, \ldots, n\}$
{\it singular interpolation nodes}.  
Sometimes we shall
refer to
Problem \ref{Prob1} as
the {\it Boundary Interpolation Problem with equality}.
A closely related problem is the following:
\begin{Pb}
Given the data set $\Omega$ as in (\ref{1.18}) and the integer
$\ka$ defined by (\ref{1.25}), 
find all functions $\w \in \N_\ka$ such that
\begin{align}
\w(x_i) &= \w_i \quad\text{and}\quad -\infty < \w^\prime(x_i) \leq
\gamma_i &(i = 1, \ldots, \ell), \label{1.29}\\
-\infty &< -\frac{1}{\w_{-1}(x_j)} \leq -\frac{1}{\xi_j} &(j = \ell+1, \ldots, n).
\label{1.30}
\end{align}
\label{Prob2}
\end{Pb}
We shall sometimes call Problem \ref{Prob2} the 
{\it Boundary Interpolation Problem with inequalities}.
\begin{Rk}
If a generalized Nevanlinna function $\w$ satisfies interpolation
conditions (\ref{1.29}) and (\ref{1.30}), then its kernel $K_\w$ 
has at least $\ka = {\rm sq}_-P$ negative squares.
\label{Russia}
\end{Rk}
To see this, let $\w \in \N_{\tilde{\ka}}$ 
for some $\tilde{\ka} \in \mathbb{Z}_+$ and 
let us assume that the following nontangential limits
exist and satisfy
\begin{align}
\w(x_i) &\in \mathbb{R} \quad\text{and}\quad
K_\w(x_i,x_i) \in \mathbb{R} && \text{for} \quad i = 1, \ldots, \ell, \label{1.31}\\
\w(x_j) &= \infty \quad\text{and}\quad
K_{-1/\w}(x_i, x_i) \in \mathbb{R} \setminus \{0\} &&\text{for} \quad
j = \ell+1, \ldots, n \label{1.32}.
\end{align}
Note that, by Theorems \ref{Tm1} and \ref{Tm2}, 
conditions (\ref{1.31}) and (\ref{1.32}) are nonrestrictive:
if $\w$ does not satisfy these conditions, then it cannot satisfy all
of the conditions in (\ref{1.29}) and (\ref{1.30}).
Conditions (\ref{1.31}) and (\ref{1.32})
guarantee that the limits ${\displaystyle \lim_{z \to x_i} \w(z)}$ and 
${\displaystyle \lim_{z \to x_i}\w^\prime(z)}$ exist and are real for $i = 1, \ldots,
\ell$ and that the limit 
${\displaystyle \lim_{z \to x_j}(z - x_j)\w(z)}$ exists and is real and 
nonzero for $j = \ell+1, \ldots, n$.
Let us introduce
the diagonal matrix
\begin{equation*}
B(\mathbf{z}) = 
\begin{bmatrix}
I_\ell & & &\\
 & z_{\ell + 1} - x_{\ell + 1}& &\\
 & & \ddots & \\
 & & & z_n - x_n
\end{bmatrix},
\quad\text{where}\quad
\mathbf{z} = (z_1, \ldots, z_n).
\end{equation*}
The nontangential limits of the entries in $P^\w(\mathbf{z})$
and $B(\mathbf{z})$ exist as $z_i \to x_i$ for $i = 1, \ldots, n$;
the nontangential limit of the product
\begin{equation*}
P^\w(\mathbf{x}) :=
\lim_{\mathbf{z} \to \mathbf{x}} 
B(\mathbf{z})^* \cdot P^\w(\mathbf{z}) \cdot B(\mathbf{z}),
\qquad \text{where }\mathbf{x} = (x_1, \ldots, x_n),
\end{equation*}
exists as well, since conditions (\ref{1.31}) and (\ref{1.32}) hold.
The straightforward calculations reveal that the matrix 
$P^\w(\mathbf{x})$ is of the
form (\ref{1.19})--(\ref{120}),
but with $\w_i$ and $\gamma_i$ replaced
respectively by $\w(x_i)$ and $\w^\prime(x_i)$ for $i = 1, \ldots, \ell$,
and with $\xi_j$ replaced by $\w_{-1}(x_j)$ for $j = \ell+1, \ldots, n$.
Since $P^\w(\mathbf{x})$ is the limit of a sequence of matrices
each with at most $\tilde{\ka}$ negative squares (since $\w \in
\N_{\tilde{\ka}}$), 
it holds that
\begin{equation}
{\rm sq}_-P^\w(\mathbf{x}) \leq \tilde{\ka}.
\label{1.36}
\end{equation}
Define the two invertible diagonal matrices $A_\Omega$ and 
$A_\w$ by
\begin{equation*}
A_\Omega = \begin{bmatrix}
I_\ell & & & \\ 
 & \xi_{\ell + 1}^{-1} & &\\
& & \ddots & \\
& & & \xi_{n}^{-1}
\end{bmatrix}
\quad\text{and}\quad
A_\w = \begin{bmatrix}
I_\ell & & & \\
 & \w_{-1}^{-1}(x_{\ell + 1}) & &\\
& & \ddots & \\
& & & \w_{-1}^{-1}(x_n)
\end{bmatrix}.
\end{equation*}
If $\w$ satisfies conditions 
(\ref{1.26})--(\ref{1.27c}),
then the off-diagonal
entries
of the matrices
$A_\Omega \cdot P \cdot A_\Omega$ and $A_\w \cdot P^\w(\mathbf{x})
\cdot A_\w$ coincide, and we have
\begin{equation*}
A_\Omega \cdot P \cdot A_\Omega - A_\w \cdot P^\w(\mathbf{x}) \cdot
A_\w =
\begin{bmatrix}
D_1 & 0\\
0 & D_2
\end{bmatrix}
\end{equation*}
where $D_1$ and $D_2$ are diagonal matrices defined as
\begin{align*}
D_1 &=
\begin{bmatrix}
\gamma_1 - \w^\prime(x_1) & & 0\\
 & \ddots & \\
0 & & \gamma_\ell - \w^\prime(x_\ell)
\end{bmatrix},
\\
D_2 &=
\begin{bmatrix}
\w_{-1}^{-1}(x_{\ell + 1}) - \xi_{\ell + 1}^{-1} & & 0\\
 & \ddots & \\
0 & & \w_{-1}^{-1}(x_{n + 1}) - \xi_n^{-1} \\
\end{bmatrix}.
\end{align*}
This shows that the conditions in Problem \ref{Prob1} and Problem
\ref{Prob2} respectively imply
\begin{equation}
A_\w \cdot P^\w(\mathbf{x}) \cdot A_\w = A_\Omega \cdot P
\cdot A_\Omega \quad\text{and}\quad
A_\w \cdot P^\w(\mathbf{x}) \cdot A_\w \leq A_\Omega \cdot
P \cdot A_\Omega.
\label{1.40}
\end{equation}
Since the matrices $A_\w$ and $A_\Omega$ are Hermitian and invertible,
we conclude by (\ref{1.36}) that each
relation in (\ref{1.40}) implies
\begin{equation*}
{\rm sq}_-P \leq \tilde{\ka}.
\end{equation*}
Thus in Problems \ref{Prob1} and \ref{Prob2} we select the minimal
possible index $\ka$ for which solutions $\w \in \N_\ka$ may exist.
In other words, the above analysis shows that, if $\ka = {\rm sq}_-P$
and if $\hat{\ka}$ is integer such that $0 \leq \hat{\ka} < \ka$,
then neither Problem \ref{Prob1} nor \ref{Prob2} admits a solution
in the class $\N_{\hat{\ka}}$.

The formulation of the next problem may appear somewhat unusual.
However, as will be shown, this problem is in some sense a more
natural generalization of the classical Nevanlinna-Pick
interpolation problem than is either Problem \ref{Prob1} or
Problem \ref{Prob2}.
\begin{Pb}
Given the data $\Omega$ as in (\ref{1.18}) and the integer
$\ka$ defined by (\ref{1.25}), 
find all functions $\w \in \N_{\ka^\prime}$ 
for $\ka^\prime \leq \ka$
such that conditions (\ref{1.29}) and (\ref{1.30})
are satisfied at all but $\ka - \ka^\prime$ points
$x_1, \ldots, x_n$.
\label{Prob3}
\end{Pb}
We shall sometimes call Problem \ref{Prob3} the
{\it Master Boundary Interpolation Problem}.

\section{Our Study}
In this thesis, we provide, for Problems \ref{Prob1}, 
\ref{Prob2}, and \ref{Prob3},
necessary and sufficient conditions for the existence of at least
one solution, and we parameterize the set of all solutions
provided the necessary conditions are met.  In the next chapter,
we shall formulate the main results of the thesis and then
outline the structure of the thesis' main body.

\chapter{Main Results}

\section{Solution to the Master Boundary Interpolation Problem}

Our study of Problems \ref{Prob1}, \ref{Prob2}, and \ref{Prob3}
splits into two principally different cases: (1) the Pick matrix
$P$ is singular; (2) $P$ is invertible.
\begin{Tm}
Suppose that $P$ is singular.  Then Problem \ref{Prob3} has a unique
solution $\w$.  Furthermore, if
$\w \in \N_\ka$, then $\w$ is also a solution of Problem \ref{Prob2}.
\label{TmDeg}
\end{Tm}

When $P$ is invertible, the situation is more interesting and
intricate.  The parameterization of all solutions will be given
in terms of a linear fractional transformation with
rational coefficients.  The matrix of these coefficients is
introduced below exclusively in terms of the interpolation
data.  Let
\begin{equation}
\Theta(z) := 
\begin{bmatrix}\Theta_{11}(z) & \Theta_{12}(z) \\
\Theta_{21}(z) & \Theta_{22}(z)\end{bmatrix}
:= I_2 - i\begin{bmatrix}C \\ E\end{bmatrix}
(zI - X)^{-1}P^{-1}\begin{bmatrix}C^* & E^*\end{bmatrix}J,
\label{ThetaFormula}
\end{equation}
where
\begin{equation}
X = \begin{bmatrix}x_1 & & 0\\ & \ddots & \\0 & & x_n\end{bmatrix},
\quad
J = \begin{bmatrix}0 & -i \\ i & 0\end{bmatrix}, \label{r2.2}
\end{equation}
\begin{equation}
C = \begin{bmatrix}\w_1 & \ldots & \w_\ell & \xi_{\ell + 1} & \ldots
& \xi_n\end{bmatrix}, 
\quad
E = \begin{bmatrix}1 & \ldots & 1 & 0 & \ldots & 0\end{bmatrix}.
\label{r2.3}
\end{equation}
Note that
the function $\Theta$ has no poles except
possibly simple poles at $x_1, \ldots, x_n$.
It will be shown in Chapter IV that $\Theta$ is $J$-unitary
on the real line, i.e., 
\begin{equation*}
\Theta(x)J\Theta(x)^* = J \quad\text{for every}\; x\in \mathbb{R}
\setminus \{x_1, \ldots, x_n\},
\end{equation*}
and that the kernel
\begin{equation}
K_{\Theta,J}(z, \zeta) := \frac{J - \Theta(z)J\Theta(\zeta)^*}{
z - \bar{\zeta}}
\label{KTkernel}
\end{equation}
has $\ka := {\rm sq}_-P$ negative squares on $\C^+$:
\begin{equation}
{\rm sq}_-K_{\Theta,J}=\ka.
\label{con2}
\end{equation}
We shall use the symbol $\mathcal{W}_\ka$ to denote the class of
$J$-unitary
$2 \times 2$ meromorphic functions satisfying (\ref{con2}).
For every $\Theta \in \mathcal{W}_\ka$, the formula
\begin{equation}
\mathbf{T}_\Theta: \varphi \to
\frac{\Theta_{11}\varphi + \Theta_{12}}{
\Theta_{21}\varphi + \Theta_{22}}
\label{trans}
\end{equation}
defines a map from
$\widetilde{\N}_0$ into ${\displaystyle \bigcup_{\ka^\prime \leq
  \ka}\widetilde{\N}_{\ka^\prime}}$.  It turns out that for
$\Theta$ of the form (\ref{ThetaFormula}), the range of 
the transformation (\ref{trans}) is the set of solutions
of Problem \ref{Prob3}:
\begin{Tm}
Let the matrices $P$, $X$, $E$, and $C$ be associated to the data set
$\Omega$ as in (\ref{1.19}) -- (\ref{120}), (\ref{r2.2}), and (\ref{r2.3}),
and let $\w$ be a function meromorphic on $\C^+$.  If $P$ is invertible,
then $\w$ is a solution of Problem \ref{Prob3} if and only if it
is of the form
\begin{equation}
\w(z) = \mathbf{T}_\Theta[\varphi](z) :=
\frac{\Theta_{11}(z)\varphi(z) + \Theta_{12}(z)}{
\Theta_{21}(z)\varphi(z) + \Theta_{22}(z)}
\label{form1}
\end{equation}
for some function $\varphi \in \widetilde{\N}_0$.
\label{lessprecise}
\end{Tm}
\begin{Rk}
\emph{As a consequence of Theorems \ref{TmDeg} and
\ref{lessprecise}, we have the following.
In the context of generalized Nevanlinna functions,
Problem \ref{Prob3} plays the same role as does the
classical Nevanlinna-Pick interpolation in the context
of Nevanlinna functions.  In particular, Problem
\ref{Prob3} always admits a solution, and in the case
where the solution is not unique, the set of all
solutions is parameterized by a linear fractional
transformation with a free extended Nevanlinna
class parameter.}
\end{Rk}
The necessity part of Theorem \ref{lessprecise} will be
obtained in Chapter V using
V.~Potapov's method of the Fundamental Matrix Inequality
(FMI) appropriately adapted
to the context of generalized Nevanlinna functions.
The proof of the sufficiency part rests on Theorems \ref{Main1} --
\ref{Main4} below,
which are of a certain independent interest.

\section{Boundary Behavior of $\mathbf{T}_\Theta[\varphi]$}
Throughout this section we assume that the Pick matrix $P$
is invertible.
Let us introduce the numbers $\tilde{c}_1, \ldots, \tilde{c}_n$ and
$\tilde{e}_1, \ldots, \tilde{e}_n$ defined by
\begin{equation*}
\tilde{c}_i := CP^{-1}\mathbf{e}_i
\quad\text{and}\quad
\tilde{e}_i := EP^{-1}\mathbf{e}_i
\quad\text{for}\quad i = 1, \ldots, n.
\end{equation*}
In Chapter IV we shall show that $\tilde{e}_i, \tilde{c}_i \in \mathbb{R}$
and $|\tilde{c}_i| + |\tilde{e}_i| > 0$ for $i = 1, \ldots, n$,
and that
\begin{equation}
\eta_i := \frac{\tilde{c}_i}{\tilde{e}_i} = -\lim_{z \to x_i}
\frac{\Theta_{22}(z)}{\Theta_{21}(z)}
\label{inparticular}
\end{equation}
belongs to $\mathbb{R} \cup \{\infty\}$ for $i = 1, \ldots, n$.
Let $\tilde{p}_{ii}$ denote the $i$-th diagonal entry of the inverse
of the Pick matrix, $P^{-1}$.
If $\eta_i \neq \infty$ (i.e., $\tilde{e}_i \neq 0$), then any function
$\varphi \in \tN_0$
satisfies precisely one of the following
conditions:
\begin{align*}
\bC_1:&\quad
\varphi(x_i) \neq \eta_i.
 \\
\bC_2:&\quad
\varphi(x_i) = \eta_i \quad \text{and} \quad
K_\varphi(x_i,x_i) = \infty.
\\
\bC_3:&\quad
\varphi(x_i) = \eta_i
\quad \text{and} \quad
-\frac{\tilde{p}_{ii}}{\tilde{e}_i^2} < \varphi^\prime(x_i) <
\infty. \\
\bC_4:&\quad
\varphi(x_i) = \eta_i 
\quad \text{and} \quad
0 \leq \varphi^\prime(x_i) < -\frac{\tilde{p}_{ii}}{
\tilde{e}_i^2}.\\
\bC_5:&\quad
\varphi(x_i) = \eta_i 
\quad \text{and} \quad
\varphi^\prime(x_i) = -\frac{\tilde{p}_{ii}}{\tilde{e}_i^2}
> 0. \\
\bC_6:&\quad
\varphi(x_i) = \eta_i 
\quad \text{and} \quad
\varphi^\prime(x_i) = \tilde{p}_{ii} = 0.
\end{align*}
Condition $\bC_1$ means that either
the limit $\varphi(x_i)$ exists and does not equal $\eta_i$, or that the
limit $\varphi(x_i)$ does not exist.  In the remaining conditions, we
suppose the limits exist and satisfy their respective (in)equalities.
We shall denote by $\mathbf{C_{i-j}}$ the disjunction of conditions
$\mathbf{C_i}$ through $\mathbf{C_j}$.
Parallel to the preceding conditions are the following:
\begin{align*}
\btC_1:&\quad
\varphi(x_i) \neq \infty. \\
\btC_2:&\quad
\varphi(x_i) = \infty
\quad \text{and} \quad
\varphi_{-1}(x_i) = 0. \\
\btC_3:&\quad
\varphi(x_i) = \infty 
\quad \text{and} \quad
-\frac{\tilde{p}_{ii}}{\tilde{c}_i^2} < -\frac{1}{\varphi_{-1}(x_i)} < \infty. \\
\btC_4:&\quad
\varphi(x_i) = \infty 
\quad \text{and} \quad
0 \leq -\frac{1}{\varphi_{-1}(x_i)} <
 -\frac{\tilde{p}_{ii}}{\tilde{c}_i^2}. \\
\btC_5:&\quad
\varphi(x_i) = \infty 
\quad \text{and} \quad
\frac{1}{\varphi_{-1}(x_i)} = \frac{\tilde{p}_{ii}}{\tilde{c}_i^2} < 0. \\
\btC_6:&\quad
\varphi(x_i) = \infty 
\quad \text{and} \quad
\frac{1}{\varphi_{-1}(x_i)} = \tilde{p}_{ii} = 0.
\end{align*}
If $\eta_i = \infty$ ($\tilde{e}_i = 0$), then any function
$\varphi \in \tN_0$ satisfies precisely one of $\btC_j$,
$j = 1, \ldots, 6$.  For this second set of conditions,
we adopt the same conventions regarding limits as we did for the first.
In the sequel, if we say that a parameter $\varphi$ satisfies a condition
$\bC_j$ at $x_i$ for some $j = 1, \ldots, 6$, then we implicitly
assume that $\tilde{e}_i \neq 0$ also.  Similarly, if we say that
a parameter $\varphi$ satisfies $\btC_j$ at $x_i$, then we implicitly assume
that $\tilde{e}_i = 0$.
The next theorems give, in terms of the parameter $\varphi$,
 a classification of the interpolation conditions
that either are or are not satisfied by a function $\w$ of 
the form (\ref{form1}).

\begin{Tm}
Let the Pick matrix $P$ be invertible, let $\varphi \in \tN_0$,
let $\Theta$ be given by
(\ref{ThetaFormula}), let $\w = \mathbf{T}_\Theta[\varphi]$,
and let $x_i$ be an interpolation node with $i \in \{1, \ldots,
\ell\}$.  Suppose that $\tilde{e}_i = 0$.  Then
\begin{enumerate}
\item{ \text{The nontangential boundary limits $\w^\prime(x_i)$ and
$\w(x_i)$ exist and are subject to}
\begin{equation*}
\w^\prime(x_i) = \gamma_i
\quad\text{and}\quad
\w(x_i) = \w_i
\end{equation*}
if and only if the parameter $\varphi$ meets either
condition $\btC_1$ or $\btC_2$.}
\item{The nontangential boundary limits $\w^\prime(x_i)$ and
$\w(x_i)$ exist and are subject to
\begin{equation*}
-\infty < \w^\prime(x_i) < \gamma_i
\quad\text{and}\quad
\w(x_i) = \w_i
\end{equation*}
if and only if the parameter $\varphi$ meets condition
$\btC_3$.}
\item{The nontangential boundary limits $\w^\prime(x_i)$ and
$\w(x_i)$ exist and are subject to
\begin{equation*}
\gamma_i < \w^\prime(x_i) < \infty
\quad\text{and}\quad
\w(x_i) = \w_i
\end{equation*}
if and only if the parameter $\varphi$ meets condition
$\btC_4$.}
\item{If $\varphi$ meets condition $\btC_5$,
then the function $\w$ is subject to one of the following:
\begin{enumerate}
\item{The nontangential boundary limit $\w(x_i)$ fails
to exist.}
\item{The limit $\w(x_i)$ exists and $\w(x_i) \neq \w_i$.}
\item{$\w(x_i) = \w_i$ and $|K_\w(x_0,x_0)| = \infty$.}
\end{enumerate} }
\item{If $\varphi$ meets condition $\btC_6$,
then the nontangential boundary limit $\w(x_i)$ exists and
\begin{equation*}
\w(x_i) \neq \w_i.
\end{equation*} }
\end{enumerate}
\label{Main1}
\end{Tm}

\begin{Tm}
Let the Pick matrix $P$ be invertible, let $\varphi \in \tN_0$,
let $\Theta$ be given by
(\ref{ThetaFormula}), let $\w = \mathbf{T}_\Theta[\varphi]$,
and let $x_i$ be an interpolation node with $i \in \{\ell+1, \ldots,
n\}$.  Suppose that $\tilde{e}_i = 0$.  Then
\begin{enumerate}
\item{ \text{The nontangential boundary limit $\w_{-1}(x_i)$ exists
and is subject to}
\begin{equation*}
\w_{-1}(x_i) = \xi_i
\end{equation*} 
if and only if the parameter $\varphi$ meets either
condition $\btC_1$ or $\btC_2$.}
\item{The nontangential boundary limit $\w_{-1}(x_i)$ exists and is subject to
\begin{equation*}
-\infty < -\w_{-1}^{-1}(x_i) < -\xi_i^{-1}
\end{equation*}
if and only if the parameter $\varphi$ meets condition
$\btC_3$.}
\item{The nontangential boundary limit $\w_{-1}(x_i)$ exists
and is subject to
\begin{equation*}
-\xi_i^{-1} < -\w_{-1}^{-1}(x_i) < \infty
\end{equation*}
if and only if the parameter $\varphi$ meets condition
$\btC_4$.}
\item{The nontangential boundary limit $\w_{-1}(x_i)$ exists
and is subject to
\begin{equation*}
\w_{-1}(x_i) = 0
\end{equation*}
if and only if the parameter $\varphi$ meets either
condition $\btC_4$ or $\btC_5$.}
\end{enumerate}
\label{Main2}
\end{Tm}

\begin{Tm}
Let the Pick matrix $P$ be invertible, let $\varphi \in \tN_0$,
let $\Theta$ be given by
(\ref{ThetaFormula}), let $\w = \mathbf{T}_\Theta[\varphi]$,
and let $x_i$ be an interpolation node with $i \in \{1, \ldots,
\ell\}$.  Suppose that $\tilde{e}_i \neq 0$.  Then
\begin{enumerate}
\item{ \text{The nontangential boundary limits $\w^\prime(x_i)$ and
$\w(x_i)$ exist and are subject to}
\begin{equation*}
\w^\prime(x_i) = \gamma_i
\quad\text{and}\quad
\w(x_i) = \w_i
\end{equation*}
if and only if the parameter $\varphi$ meets either
condition $\bC_1$ or $\bC_2$.}
\item{The nontangential boundary limits $\w^\prime(x_i)$ and
$\w(x_i)$ exist and are subject to
\begin{equation*}
-\infty < \w^\prime(x_i) < \gamma_i
\quad\text{and}\quad
\w(x_i) = \w_i
\end{equation*}
if and only if the parameter $\varphi$ meets condition
$\bC_3$.}
\item{The nontangential boundary limits $\w^\prime(x_i)$ and
$\w(x_i)$ exist and are subject to
\begin{equation*}
\gamma_i < \w^\prime(x_i) < \infty
\quad\text{and}\quad
\w(x_i) = \w_i
\end{equation*}
if and only if the parameter $\varphi$ meets condition
$\bC_4$.}
\item{If $\varphi$ meets condition $\bC_5$,
then the function $\w$ is subject to one of the following:
\begin{enumerate}
\item{The nontangential boundary limit $\w(x_i)$ fails
to exist.}
\item{The limit $\w(x_i)$ exists and $\w(x_i) \neq \w_i$.}
\item{$\w(x_i) = \w_i$ and $|K_\w(x_0,x_0)| = \infty$.}
\end{enumerate} }
\item{If $\varphi$ meets condition $\bC_6$,
then the nontangential boundary limit $\w(x_i)$ exists and
\begin{equation*}
\w(x_i) \neq \w_i.
\end{equation*} }
\end{enumerate}
\label{Main3}
\end{Tm}

\begin{Tm}
Let the Pick matrix $P$ be invertible, let $\varphi \in \tN_0$,
let $\Theta$ be given by
(\ref{ThetaFormula}), let $\w = \mathbf{T}_\Theta[\varphi]$,
and let $x_i$ be an interpolation node with $i \in \{\ell+1, \ldots,
n\}$.  Suppose that $\tilde{e}_i \neq 0$.  Then
\begin{enumerate}
\item{ \text{The nontangential boundary limit $\w_{-1}(x_i)$ exists
and is subject to}
\begin{equation*}
\w_{-1}(x_i) = \xi_i
\end{equation*} 
if and only if the parameter $\varphi$ meets either
condition $\bC_1$ or $\bC_2$.}
\item{The nontangential boundary limit $\w_{-1}(x_i)$ exists and is subject to
\begin{equation*}
-\infty < -\w_{-1}^{-1}(x_i) < -\xi_i^{-1}
\end{equation*}
if and only if the parameter $\varphi$ meets condition
$\bC_3$.}
\item{The nontangential boundary limit $\w_{-1}(x_i)$ exists
and is subject to
\begin{equation*}
-\xi_i < -\w_{-1}^{-1}(x_i) < \infty
\end{equation*}
if and only if the parameter $\varphi$ meets condition
$\bC_4$.}
\item{The nontangential boundary limit $\w_{-1}(x_i)$ exists
and is subject to
\begin{equation*}
\w_{-1}(x_i) = 0
\end{equation*}
if and only if the parameter $\varphi$ meets either
condition $\bC_5$ or $\bC_6$.}
\end{enumerate}
\label{Main4}
\end{Tm}

\section{Solutions to the Boundary Interpolation Problem with
equality and the Boundary Interpolation Problem with inequalities}
As an immediate consequence of Theorems \ref{Main1} -- \ref{Main4}, we get
\begin{Cy}
A function $\w = \mathbf{T}_\Theta[\varphi]$ meets the
$i$-th interpolation conditions for Problem \ref{Prob2}:
\begin{align*}
\w(x_i) &= \w_i
\quad\text{and}\quad
-\infty < \w^\prime(x_i) \leq \w_i
&&\quad\text{for}\quad i \in \{1, \ldots, \ell\}, \\
-\infty &< -\w_{-1}^{-1}(x_i) \leq -\xi_i^{-1}
&&\quad\text{for}\quad i \in \{\ell + 1, \ldots, n\},
\end{align*}
if and only if the corresponding parameter $\varphi \in \tN_0$
meets the conditions $\bC_{1-3} \vee \btC_{1-3}$ at interpolation
node $x_i$.
\label{Cy24}
\end{Cy}
The next theorem relates the number of negative squares of the
kernel $K_\w$ for $\w = \mathbf{T}_\Theta[\varphi]$ to the
number of points at which the parameter $\varphi$ satisfies 
conditions $\bC_{4-6} \vee \btC_{4-6}$.
\begin{Tm}
If the Pick matrix $P$ is invertible and has $\ka$ negative
eigenvalues, then a function $\varphi \in \tN_0$ may satisfy
conditions $\bC_{4-6} \vee \btC_{4-6}$ at at most $\ka$
interpolation nodes.  Furthermore, if $\varphi$ meets conditions
 $\bC_{4-6} \vee \btC_{4-6}$ at exactly $k$ interpolation
nodes, then the function $\w = \mathbf{T}_\Theta[\varphi]$
belongs to the class $\N_{\ka - k}$. 
\label{Tm25}
\end{Tm}
Together Corollary \ref{Cy24} and Theorem \ref{Tm25}
imply the sufficiency part in Theorem \ref{lessprecise}.
Theorems \ref{lessprecise} and \ref{Tm25} also lead us
to the parameterizations of the solution sets of Problems
\ref{Prob1} and \ref{Prob2}.  
In particular, Theorem \ref{Tm25} indicates
which parameters lead to solutions and which do not.
\begin{Tm}
A function $\w$ of the form (\ref{form1}) is a solution of
Problem \ref{Prob1} if and only if the corresponding paramter
$\varphi \in \tN_0$ meets the conditions
$\bC_{1-2} \vee \btC_{1-2}$ at each interpolation node $x_i$.
\label{Tm26}
\end{Tm}
\begin{Tm}
A function $\w$ of the form (\ref{form1}) is a solution of
Problem \ref{Prob2} if and only if the corresponding paramter
$\varphi \in \tN_0$ meets the conditions
$\bC_{1-3} \vee \btC_{1-3}$ at each interpolation node $x_i$.
\label{Tm27}
\end{Tm}
Theorems \ref{lessprecise} and \ref{Tm27} also yield necessary
and sufficient conditions for Problems \ref{Prob2} and
\ref{Prob3} to be equivalent (i.e., to have the same set of
solutions).
\begin{Cy}
Problems \ref{Prob2} and \ref{Prob3} are equivalent if and only if
all of the diagonal entries of the inverse of the Pick matrix
$P$ are positive.
\label{Core2}
\end{Cy}
\begin{proof}
Indeed, there are no parameters $\varphi \in \tN_0$ satisfying
conditions $\bC_{4-6} \vee \btC_{4-6}$ at any interpolation node
$x_i$ if and only if each $\tilde{p}_{ii} > 0$, i.e.,
all the diagonal entries of $P^{-1}$ are positive.   The corollary
then follows from Theorems \ref{lessprecise} and \ref{Tm27}.
\end{proof}

In view of the preceding results, one may ask whether there
exists a parameter $\varphi \in \tN_0$
such that $\w = \mathbf{T}_\Theta[\varphi]$
``misses'' interpolation conditions (\ref{1.29}) and
(\ref{1.30}) at prescribed interpolation nodes
$x_{i_1}, \ldots, x_{i_k}$.  It turns out that this
problem admits a simple answer in terms of a certain
principal submatrix of the Pick
matrix $P$. 
\begin{Tm}
There exists a function $\varphi \in \tN_0$ satisfying
conditions $\bC_{4-6} \vee \tilde{\bC}_{4-6}$
at $x_{i_1}, \ldots, x_{i_k}$ if and only if the $k \times k$
matrix
\begin{equation*}
\mathcal{P} := \left[\tilde{p}_{i_\alpha,i_\beta}\right]_{\alpha,\beta=1}^k
\end{equation*}
is negative semidefinite.  If $\mathcal{P}$ is invertible, then there
are infinitely many such functions, and if $\mathcal{P}$ is singular,
then there is a unique such function.
\label{Tm29}
\end{Tm}

The rest of the thesis is organized as follows.
Chapters III and IV contain preliminaries needed for the subsequent
analysis.  Chapter V is dedicated to proving the necessity part
of Theorem \ref{lessprecise}.  In Chapter VI, we prove Theorems
\ref{Main1} and \ref{Main2}, and in Chapter VII we prove Theorems
\ref{Main3} and \ref{Main4}.  Note that this completes
the proof of Corollary \ref{Cy24}.
In Chapter VIII, we prove Theorems
\ref{Tm25} and \ref{Tm29}.  This will complete the proof
of Theorem \ref{lessprecise}, and hence the proofs
of Theorems \ref{Tm26} and \ref{Tm27} and
Corollary \ref{Core2} as well.
We prove Theorem \ref{TmDeg} in Chapter IX.  Finally,
in Chapter X we present some illustrative numerical examples,
to which a
wide range of the cases presented
in Theorems \ref{Main1} -- \ref{Main4}
is applicable, and 
in addition, we present a degenerate numerical interpolation problem
and find its unique solution.

\chapter{Generalized Nevanlinna and Schur Class Functions}

In this section we prove Theorems \ref{Tm1} and \ref{Tm2} along with some
needed auxiliary results.

\section{Generalized Schur Functions}
Related to $\N_\ka$
is the generalized Schur class $\Sch_\ka$, defined
below.
\begin{Dn}
Let $\ka$ be a nonnegative integer.
A function $S$ belongs to the \emph{generalized Schur class}
$\Sch_\ka$ provided that $S$ is meromorphic on
$\mathbb{D}$ and is such that the kernel
\begin{equation*}
K_S^{\Sch}(z, \zeta):=
\frac{1 - S(z)S(\zeta)^*}{1 - z\bar{\zeta}} \quad \quad
(z, \zeta \in \rho(S))
\end{equation*}
has $\ka$ negative squares on $\rho(S) \times \rho(S)$.
\end{Dn}
Generalized Schur functions appeared in \cite{Tak}
and \cite{Akh} and
were comprehensively studied in \cite{KL1} and \cite{KL2}.
It was shown in particular that
to each $S \in \Sch_\ka$ corresponds an essentially unique 
(i.e., unique up to a unimodular constant)
Krein-Langer
representation
\begin{equation}
S = \frac{f}{b}
\label{KL}
\end{equation}
where $f \in \Sch_0$ is a Schur function and $b$ is a finite
Blaschke product of degree $\ka$.
\begin{La}
Let $b$ be a finite Blaschke product of degree $\ka$:
\begin{equation}
b(z) = \prod_{j=1}^\ka\frac{z - c_j}{1 - z\bar{c}_j}
\quad \quad (c_j \in \mathbb{D}, \; j=1,\ldots,\ka)
\label{Blaschke}
\end{equation}
and let $t_0 \in \mathbb{T}$.
Then the limit
\begin{equation}
d_b := \lim_{z \to t_0} K_b^{\Sch}(z,z), \quad\quad (z \in \mathbb{D})
\label{dbKb}
\end{equation}
exists and is positive.
\label{BlaschkeLemma}
\end{La}
\begin{proof}
Observe that
\begin{align}
|b(z)|^2 = \prod_{j=1}^\ka \frac{|z - c_j|^2}{|1 - z\bar{c}_j|^2}
&= \prod_{j=1}^\ka \frac{|1 - z\bar{c}_j|^2 - (1 - |c_j|^2)(
1 - |z|^2)}{|1 - z\bar{c}_j|^2} \nonumber \\
&= \prod_{j=1}^\ka \left(1 - \frac{(1 - |c_j|^2)(1 - |z|^2)}{
|1 - z\bar{c}_j|^2}\right).
\label{bchain}
\end{align}
Since $c_j \in \mathbb{D}$, the latter equality
provides that the limit
\begin{equation}
d_b:= \lim_{z \to t_0}K_b^\Sch(z,z) =
\lim_{z \to t_0}\frac{1 - |b(z)|^2}{1 - |z|^2}
\label{bker}
\end{equation}
exists and is finite.  
Only the terms in $1 - |b(z)|^2$ 
with a factor $1 - |z|^2$ of multiplicity
exactly one have a nonzero limit, and it is easy to check by
expanding (\ref{bchain}) that the limit of each of these terms 
in (\ref{bker}) is positive, whence the limit (\ref{bker}) is
itself positive.  Note that even the tangential limit
${\displaystyle \lim_{z \to t_0} K_b^\Sch(z,z)}$ 
exists and equals $d_b$, since $b$ is rational.
\end{proof}
We shall now study the connection between the classes
$\tN_\ka$ and $\Sch_\ka$.
\begin{La}
Let $K_1$ and $K_2$ be two sesqui-analytic kernels
(i.e., analytic in $z$ and in $\bar{\zeta}$)
satisfying
\begin{equation*}
K_1(z,\zeta) = \varphi(\zeta)^* \cdot K_2(z,\zeta) \cdot \varphi(z),
\quad \quad\quad (z,\zeta \in \Delta)
\end{equation*}
where $\varphi$ is meromorphic on $\Delta$ and does not vanish
identically there.
Then
\begin{equation*}
{\rm sq}_-K_1 = {\rm sq}_-K_2 \quad \text{and} \quad
{\rm sq}_+K_1 = {\rm sq}_+K_2.
\end{equation*}
\label{Inertia}
\end{La}
\begin{proof}
For the proof, see \cite{AR}.
\end{proof}
\begin{La}
Let $\w \in \widetilde{\N}_\ka$.  Then
the function 
\begin{equation}
S:=\beta \circ \w \circ \beta^{-1},
\quad\text{where}\quad
\beta(z) = \frac{z-i}{z+i},
\label{Sformula}
\end{equation}
defined on
$D_S = \beta(\rho(\w))$
belongs to $\Sch_\ka$.
\label{NSLemma}
\end{La}
\begin{proof}
If $\w \equiv \infty$, then
(\ref{Sformula}) shows that $S \equiv 1$.  Let us henceforth
assume that $\w \not \equiv \infty$.
Formula (\ref{Sformula}) written in more detail as
\begin{equation}
S(\zeta) = \frac{\w(i\frac{\zeta+1}{\zeta-1}) - i}{
\w(i\frac{\zeta+1}{\zeta-1}) + i}
\label{Sdetail}
\end{equation}
shows that $S$ is meromorphic on $D_S$ provided that
\begin{equation}
w \not \equiv -i.
\label{wequivi}
\end{equation}
Let $\zeta_1, \zeta_2 \in \beta(\rho(\w))$ and set 
$z_j = \beta^{-1}(\zeta_j),
j=1,2$.  It is readily checked that
\begin{equation*}
1 - S(\zeta_1)S(\zeta_2)^* =
\frac{-2i \w(z_1) + 2i \w(z_2)^*}{(\w(z_1) + i)(\w(z_2)^* - i)}
\end{equation*}
and that
\begin{equation*}
1 - \zeta_1 \bar{\zeta}_2 = \frac{1}{2i}(1 - \beta(z_1))(
z_1 - \bar{z}_2)(1 - \beta(z_2)^*).
\end{equation*}
Combining these two relations yields
\begin{equation}
\frac{1 - S(\zeta_1)S(\zeta_2)^*}{1 - \zeta_1\bar{\zeta}_2} =
\varphi(\beta(z_1)) \cdot
\frac{\w(z_1) - \w(z_2)^*}{z_1 - \bar{z}_2} \cdot
\varphi(\beta(z_2))^*
\label{NSker}
\end{equation}
where
\begin{equation}
\varphi(\beta(z_j)) = \frac{2}{(\w(z_j)+i)(1 - \beta(z_j))} \quad
\quad (j=1,2).
\label{phi}
\end{equation}
From (\ref{phi}) it is clear that, if (\ref{wequivi}) holds,
then  $\varphi$ is meromorphic on $D_S$ and does not
vanish identically there.  
Thus if (\ref{wequivi}) is satisfied, then
Lemma \ref{Inertia} implies that,
\begin{equation*}
{\rm sq}_-K_\w = {\rm sq}_-K_S^{\Sch} = \ka,
\end{equation*}
from which the desired conclusion follows.  
The argument will be complete if we prove that
(\ref{wequivi}) always holds.  We do this
in the next lemma by showing that
$\w \equiv -i$ does not belong to $\N_\ka$ for any
$\ka \geq 0$.
\end{proof}
\begin{La}
The function $\w \equiv -i$ is not a generalized
Nevanlinna function, i.e., $\w \not \in \N_\ka$
for any $\ka \geq 0$.
\end{La}
\begin{proof}
Choose $n$ distinct points $\zeta_1, \ldots, \zeta_n \in \mathbb{D}$.
For each $k = 1, \ldots, n$, set
\begin{equation*}
z_k = \beta^{-1}(\zeta_k) = i\frac{1 + \zeta_k}{1 - \zeta_k} \in \C^+.
\end{equation*}
Since the following chain of equalities holds,
\begin{equation*}
\frac{-2i}{z_j - \bar{z}_k} = \frac{-2i}{
i\frac{1 + \zeta_j}{1 - \zeta_j} + i\frac{1 + \bar{\zeta}_k}{
1 - \bar{\zeta}_k}}
= -\frac{(1 - \zeta_j)(1 - \bar{\zeta}_k)}{
1 - \zeta_j \bar{\zeta}_k},
\end{equation*}
it holds that
\begin{equation*}
\left[ K_{-i}(z_j, z_k) \right]_{k,j = 1}^n =
 -\begin{bmatrix}1 - \zeta_1 & &0 \\
 & \ddots & \\
 0& & 1 - \zeta_n\end{bmatrix}
\left[ \frac{1}{1 - \zeta_j \bar{\zeta}_k} \right]_{k,j = 1}^n 
\begin{bmatrix}1 - \bar{\zeta}_1 & & 0 \\
 & \ddots & \\
0 & & 1 - \bar{\zeta}_n\end{bmatrix}.
\end{equation*}
To complete the proof we need only show that the matrix
\begin{equation*}
G := \left[ \frac{1}{1 - \zeta_j \bar{\zeta}_k} \right]_{k,j = 1}^n
\end{equation*}
is positive definite.
This follows from the fact that $G$ is
the Gram matrix of the linearly independent system of functions
\begin{equation*}
z \mapsto \frac{1}{1 - z\bar{\zeta}_k} \quad \quad (k = 1, \ldots, n)
\end{equation*}
considered as vectors in the Hardy space $H^2(\mathbb{D})$.
\end{proof}

\section{Carath\'eodory-Julia Theorem for Generalized Schur Functions}

We now state the Carath\'eodory-Julia Theorem for Schur functions.
\begin{Tm}
Let $f \in \Sch_0$ and $t_0 \in \mathbb{T}$.  Then the following
statements are equivalent:
\begin{enumerate}
\item{\begin{equation*}
d := \liminf_{z \to t_0}K_f^\Sch(z,z) < \infty.
\end{equation*}}
\item{The nontangential limit
${\displaystyle \tilde{d} := \lim_{z \to t_0}K_f^\Sch(z,z)}$ exists
and is finite.}
\item{The nontangential limits
\begin{equation*}
f_0 := \lim_{z \to t_0}f(z)
\quad\text{and}\quad
f_1 := \lim_{z \to t_0}f^\prime(z)
\end{equation*}
exist, $f_0 \in \mathbb{T}$, and
$f_1 \bar{f}_0 t_0 \geq 0$.}
\item{The nontangential limits
\begin{equation*}
f_0 := \lim_{z \to t_0}f(z)
\quad\text{and}\quad
\tilde{f}_1 := \lim_{z \to t_0}
\frac{f(z) - f_0}{z - t_0}
\end{equation*}
exist, $f_0 \in \mathbb{T}$, and
$\tilde{f}_1 \bar{f}_0 t_0 \geq 0
$.}
\end{enumerate}
In this case, $d = \tilde{d} = f_1\bar{f}_0t_0 = \tilde{f}_1 \bar{f}_0 t_0$.
\label{CJS}
\end{Tm}
\begin{proof}
For the proof, see \cite{Sar}.
\end{proof}
\begin{Tm}
Let $\ka \in \mathbb{Z}_\ka$,
$S \in \Sch_\ka$ and $t_0 \in \mathbb{T}$.  Then the following
statements are equivalent:
\begin{enumerate}
\item{\begin{equation}
d := \liminf_{z \to t_0}K_S^\Sch(z,z) < \infty.
\label{E323}
\end{equation}}
\item{The nontangential limit
${\displaystyle \tilde{d} := \lim_{z \to t_0}K_S^\Sch(z,z)}$ exists
and is real.}
\item{The nontangential limits
\begin{equation*}
S_0 := \lim_{z \to t_0}S(z)
\quad\text{and}\quad
S_1 := \lim_{z \to t_0}S^\prime(z)
\end{equation*}
exist, $S_0 \in \mathbb{T}$, and
$S_1 \bar{S}_0 t_0 \in \mathbb{R}$.}
\item{The nontangential limits
\begin{equation*}
S_0 := \lim_{z \to t_0}S(z)
\quad\text{and}\quad
\tilde{S}_1 := \lim_{z \to t_0}
\frac{S(z) - S_0}{z - t_0}
\end{equation*}
exist, $S_0 \in \mathbb{T}$, and
$\tilde{S}_1 \bar{S}_0 t_0 \in 
\mathbb{R}$.}
\end{enumerate}
In this case, $d = \tilde{d} = S_1\bar{S}_0t_0 = \tilde{S}_1\bar{S}_0t_0$.
\label{CJGS}
\end{Tm}
\begin{proof}
We first remark that, with $\ka = 0$, the above theorem reduces
to Theorem \ref{CJS}, but without the claim that $d \geq 0$.  We
shall derive the general result from Theorem \ref{CJS} by using the
Krein-Langer representation (\ref{KL}) of $S$.
We split the proof into several steps, showing that
$1. \implies 2. \implies 3. \implies 4. \implies 1.$

\nline
{\bf Step 1.}
\emph{Let $S$ be represented as in (\ref{KL}).  Then
\begin{equation}
\liminf_{z \to t_0}K_S^{\Sch}(z,z) < \infty \iff
\liminf_{z \to t_0}K_f^{\Sch}(z,z) < \infty.
\label{Sifff}
\end{equation} }
\nline

{\bf Proof of Step 1:}
Due to (\ref{KL}),
the kernel $K_S^{\Sch}$ takes the form
\begin{align}
K_S^{\Sch}(z,z) = \frac{1 - |S(z)|^2}{1 - |z|^2} 
&= \frac{1}{|b(z)|^2}\left(\frac{1 - |f(z)|^2}{
1 - |z|^2} - \frac{1 - |b(z)|^2}{1 - |z|^2}\right) 
\nonumber \\
&= \frac{1}{|b(z)|^2}\left(K_f^{\Sch}(z,z) -
K_b^{\Sch}(z,z)\right).
\label{Ksfb}
\end{align}
From (\ref{Blaschke}) it is clear that $b$ is rational, has $\ka$
poles, all in $\C \setminus \bar{\mathbb{D}}$, and
is unimodular on $\mathbb{T}$. 
For notational convenience, we set
$b_0 := b(t_0)$ and $b_1 := b^\prime(t_0)$.
By Lemma \ref{BlaschkeLemma},
\begin{equation}
0 < d_b := \lim_{z \to t_0}K_b^{\Sch}(z,z)
< \infty,
\label{db}
\end{equation}
and moreover ${\displaystyle d_b = \liminf_{z \to x_0}K_b^\Sch(z,z)}$
as well since $b(z)$ is rational.
By (\ref{db}), taking limits
in (\ref{Ksfb}) establishes (\ref{Sifff}).

\nline
{\bf Step 2.} \emph{Statements 1 -- 4 are equivalent.}
\nline

{\bf Proof of implication $1. \implies 2.$:}
Due
to Step 1 it holds that
\begin{equation*}
d_f := \liminf_{z \to t_0} K_f^{\Sch}(z,z) < \infty,
\end{equation*}
which implies, by Theorem \ref{CJS}, that the limit
\begin{equation*}
\tilde{d}_f := \lim_{z \to t_0} K_f^\Sch(z,z)
\end{equation*}
exists and that $\tilde{d}_f = d_f$.  Combining this with
(\ref{db}) and (\ref{Ksfb}) shows that
\begin{equation}
\tilde{d} := \lim_{z \to t_0} K_S^\Sch(z,z) = d_f - d_b.
\label{E328}
\end{equation}
Since $d_f - d_b \in \mathbb{R}$, the proof of Step 2
is complete.
\nline

{\bf Proof of implication $2. \implies 3.$:}
On account of (\ref{db}),
relation (\ref{Ksfb}) shows that
the limit 
\begin{equation*}
\tilde{d}_f := \lim_{z \to t_0}K_f(z,z)
\end{equation*}
exists
and is finite.
Theorem \ref{CJS} therefore
implies that the nontangential limits
\begin{equation*}
f_0 := \lim_{z \to t_0}f(z) 
\quad\text{and}\quad
f_1 := \lim_{z \to t_0}f^\prime(z)
\end{equation*}
exist, that $\tilde{d}_f \geq 0$, and that
\begin{equation}
\tilde{d}_f = f_1\bar{f}_0t_0.
\label{e3000}
\end{equation}
The existence of these limits in turn implies the
existence of the following nontangential limits:
\begin{align}
S_0 &:= \lim_{z \to t_0}S(z) = \lim_{z \to t_0}\frac{f(z)}{b(z)} =
\frac{f_0}{b_0}, \label{S0}\\
S_1 &:= \lim_{z \to t_0}S^\prime(z) = \lim_{z \to t_0}\frac{f^\prime(z)
b(z) - f(z)b^\prime(z)}{b(z)^2} = \frac{f_1b_0 - f_0b_1}{b_0^2}.
\label{S1}
\end{align}
Since both $f_0$ and $b_0$ are unimodular, $S_0$ is unimodular as well.
The nontangential limit of the right hand side in (\ref{Ksfb}) exists and equals
$\tilde{d}_f - d_b$, and therefore
\begin{equation}
\tilde{d} := \lim_{z \to t_0}K_S^{\Sch}(z,z) = \tilde{d}_f - d_b.
\label{dEq}
\end{equation}
Substituting (\ref{e3000})
into (\ref{dEq}) 
and using $d_b = b_1\bar{b}_0t_0$
yields
\begin{equation}
\tilde{d} = f_1\bar{f}_0t_0 - b_1\bar{b}_0t_0 
= \frac{1}{|b_0|^2}\left(f_1\bar{f}_0t_0 - \frac{b_1}{b_0}t_0\right)
= \left(\frac{f_1b_0-
    f_0b_1}{b_0^2}\right)\frac{\bar{f}_0}{\bar{b}_0}t_0.
\label{tdSp}
\end{equation}
Substituting (\ref{S0}) and (\ref{S1}) into (\ref{tdSp}) in turn
yields
\begin{equation}
\tilde{d} = S_1\bar{S}_0t_0.
\label{dS1S0t0}
\end{equation}
Thus $S_1 \bar{S}_0t_0 \in \mathbb{R}$ since $\tilde{d} \in \mathbb{R}$.
\nline

{\bf Proof of implication $3. \implies 4.$:}
The Krein-Langer representation (\ref{KL}) of $S$ and
the existence of the limit $S_0$ proves the existence of the limit
\begin{equation}
f_0 := \lim_{z \to t_0}f(z) = \lim_{z \to t_0}b(z)S(z) = b_0S_0.
\label{f0}
\end{equation}
Since both $b_0$ and $S_0$ are unimodular, $f_0$ is as well.
From (\ref{f0}) and the existence of the limit $S_1$, we
get that the limit
\begin{equation}
f_1 := \lim_{z \to t_0}f^\prime(z) = \lim_{z \to t_0}
\left(S^\prime(z)b(z) + 
f(z)\frac{b^\prime(z)}{b(z)}\right) = S_1b_0 + f_0\frac{b_1}{b_0}.
\label{139}
\end{equation}
exists.
Together (\ref{S0}) and (\ref{139}) show that
\begin{align*}
f_1 &= S_1 b_0 \cdot \bar{S}_0 S_0 +
f_0 \frac{b_1}{b_0}\\
&= S_1 \bar{S}_0 b_0 \cdot \frac{f_0}{b_0} + f_0b_1\bar{b}_0.
\end{align*}
Hence
\begin{equation*}
f_1 \bar{f}_0 t_0 = S_1\bar{S}_0 t_0 + b_1\bar{b}_0 t_0,
\end{equation*}
which shows that $f_1 \bar{f}_0 t_0 \in \mathbb{R}$.
Therefore, by Theorem \ref{CJS} applied to $f$, the limit
\begin{equation}
\tilde{f}_1 = \lim_{z \to t_0} \frac{f(z) - f_0}{z - t_0}
\label{fpbp}
\end{equation}
exists and $\tilde{f}_1 = f_1$.
Making use of (\ref{KL}), we have
\begin{equation}
\frac{S(z) - S_0}{z - t_0} =
\frac{1}{b_0 b(z)} \cdot \frac{f(z)b_0 - f_0b(z)}{z - t_0} = 
\frac{1}{b_0 b(z)}\left(\frac{f(z) - f_0}{z - t_0}b_0 - f_0\frac{b(z)
    - b_0}{z - t_0}\right).
\label{sslim}
\end{equation}
Due to (\ref{fpbp}) and the analyticity of $b$ on $\mathbb{T}$, the
limit in (\ref{sslim}) exists and satisfies
\begin{equation}
\tilde{S}_1 := 
\lim_{z \to t_0}\frac{S(z) - S_0}{z - t_0} =
\frac{f_1b_0 - f_0b_1}{b_0^2} = S_1.
\label{tSlim}
\end{equation}
Since $S_1 \bar{S_0} t_0 \in \mathbb{R}$, we have
$\tilde{S}_1 \bar{S_0} t_0 \in \mathbb{R}$ as well
due to (\ref{tSlim}).
\nline

{\bf Proof of implication $4. \implies 1.$:}
The Krein-Langer representation (\ref{KL}) of $S$ and
the existence of the limit $S_0$ prove the existence of the limit
\begin{equation}
f_0 := \lim_{z \to t_0}f(z) = \lim_{z \to t_0}b(z)S(z) = b_0S_0,
\label{f0step5}
\end{equation}
and since both $b_0$ and $S_0$ are unimodular, $f_0$ is as well.
Therefore equality (\ref{sslim}) holds.
The existence of the limit $\tilde{S}_1$ in conjunction with
(\ref{sslim}) implies that the limit
\begin{equation*}
\tilde{f}_1 := \lim_{z \to t_0}\frac{f(z) - f_0}{z - t_0}
\end{equation*}
exists and is finite.  Applying Theorem \ref{CJS} to
$f$ proves that the limit
\begin{equation*}
d_f := \lim_{z \to t_0}K_f^\Sch(z,z)
\end{equation*}
exists and is finite.
Therefore, on account of (\ref{Ksfb}), the limit
\begin{equation*}
\tilde{d} := \lim_{z \to t_0}K_S^\Sch(z,z)
\end{equation*}
exists and is real.
Since the limit $\tilde{d} := \lim_{z \to t_0}K_S^\Sch(z,z)$ 
is nontagential and the limit inferior
${\displaystyle d := \liminf_{z \to t_0}K_S^\Sch(z,z)}$
is arbitrary, it holds that
\begin{equation}
d \leq \tilde{d}.
\label{dleqtd}
\end{equation}
In particular, since $\tilde{d} < \infty$, it follows that
$d < \infty$.

\nline
{\bf Step 3.} \emph{If Statements 1--4 hold, then $d = \tilde{d} = S_1\bar{S}_0t_0 =
\tilde{S}_1\bar{S}_0t_0$.}
\nline

{\bf Proof of Step 3:}
Combining (\ref{dS1S0t0}), (\ref{tSlim}), 
and (\ref{dleqtd}), we have
\begin{equation*}
d \leq \tilde{d} = S_1\bar{S}_0t_0 = \tilde{S}_1\bar{S}_0t_0.
\end{equation*}
To complete the proof, we show that
$\tilde{d} \leq d$
holds as well.
Towards this end, take liminf of each side of equation (\ref{Ksfb}).  The left
hand side is $d$ by definition.  Let $\left\{z_n\right\}$ be a
sequence of points in $\mathbb{D}$ tending to $t_0$ such that
the right hand side in (\ref{Ksfb}) tends to $d$ as $n \to \infty$:
\begin{equation*}
d = \lim_{n \to \infty} \frac{1}{|b(z_n)|^2}\left(K_f^{\Sch}(z_n,z_n)
- K_b^{\Sch}(z_n,z_n)\right).
\end{equation*}
We have already established in Lemma \ref{BlaschkeLemma} that 
$d_b = \lim_{z \to t} K_b^{\Sch}(z,z)$ 
always exists whenever $t \in \mathbb{T}$, and 
that $\lim_{z \to t}b(z) = b(t)$ is unimodular.  Therefore,
\begin{equation}
d + d_b = \lim_{n \to \infty} K_f^{\Sch}(z_n,z_n) < \infty.
\label{ddbdf}
\end{equation}
The right hand side of (\ref{ddbdf})
is not less than $d_f$.
Hence
\begin{equation*}
d + d_b \geq d_f,
\end{equation*}
and so by (\ref{E328}),
\begin{equation*}
\tilde{d} = d_f - d_b \leq d,
\end{equation*}
as required.
\end{proof}

\section{Carath\'eodory-Julia Theorem for Generalized Nevanlinna Functions}

\begin{Tm}
(Theorem \ref{Tm1})
Let $\ka \in \mathbb{Z}_+$,
$\w \in \N_\ka$ ,and $x_0 \in \mathbb{R}$.  Then the following
statements are equivalent:
\begin{enumerate}
\item{\begin{equation*}
\liminf_{z \to x_0} |\w(z)| < \infty
\quad\text{and}\quad
-\infty < d := \liminf_{z \to x_0}K_\w(z,z) < \infty.
\end{equation*}}
\item{The nontangential limits
\begin{equation*}
\w_0 := \lim_{z \to x_0}\w(z)
\quad\text{and}\quad
\tilde{d} := \lim_{z \to x_0}K_\w(z,z)
\end{equation*}
exist and are real.}
\item{The nontangential limits
\begin{equation*}
\w_0 := \lim_{z \to x_0}\w(z)
\quad\text{and}\quad
\w_1 := \lim_{z \to x_0}\w^\prime(z)
\end{equation*}
exist and are real.}
\item{The nontangential limits
\begin{equation*}
\w_0 := \lim_{z \to x_0}\w(z)
\quad\text{and}\quad
\tilde{\w}_1 := \lim_{z \to x_0}
\frac{\w(z) - \w_0}{z - x_0}
\end{equation*}
exist and are real.}
\end{enumerate}
In this case, $d = \tilde{d} = \w_1 = \tilde{\w}_1$.
\label{CJGN1}
\end{Tm}

\begin{proof}
We shall prove
that $1. \implies 2. \implies 3. \implies 4.
\implies 1.$
\nline

{\bf Proof of implication $1. \implies 2.$:}
The function $S$ associated to $\w$
via formula (\ref{Sformula}) belongs to $\Sch_\ka$ by Lemma \ref{NSLemma}.  
Let $\{z_n\}$ be a sequence of points in
$\rho(\w)$ tending to $x_0$ such that
$K_\w(z_n,z_n) \to d$.  Since $\beta:\mathbb{R} \to \mathbb{T}
\setminus \{1\}$, 
the limit $t_0 := \lim_{n \to \infty}\beta(z_n)$
exists, is unimodular, and is not equal to one.  
Thus in setting $\zeta_n = \beta(z_n),
n = 1,2,\ldots$, we have $\zeta_n \to t_0$, 
$\zeta_n \in \mathbb{D}$ for $n = 1, 2, \ldots$.
Since $K_\w$ is bounded
along $\{z_n\}$ for $n$ large enough, the function
$\w(z_n)$ does not tend to $-i$ as $n \to \infty$, 
and so there exists an
$\epsilon > 0$ along with a subsequence $\{z_{n_j}\}$ 
such that $|w(z_{n_j}) + i| > \epsilon$ for $j = 1,2,\ldots$.
The function $\varphi \circ \beta$, as defined in (\ref{phi}), 
is therefore bounded along the subsequence $\{z_{n_j}\}$ for
$j$ large enough.
On account of (\ref{NSker}) and the boundedness
of $\varphi \circ \beta$ and $K_\w$, it follows that 
$K_S^{\Sch}(\zeta_{n_j}, \zeta_{n_j})$
is also bounded for $j$ large enough. In particular,
\begin{equation*}
d_S := \liminf_{\zeta \to t_0} K_S^{\Sch}(\zeta,\zeta) < \infty.
\end{equation*}
By Theorem \ref{CJGS}, the nontangential limit 
${\displaystyle \lim_{\zeta \to t_0} K_S^{\Sch}(\zeta,\zeta)}$ exists and equals
$d_S$, and the limit ${\displaystyle S_0 = \lim_{\zeta \to t_0}S(\zeta)}$ exists
and is unimodular.  Hence (\ref{Sformula}) shows that the
limit ${\displaystyle \w_0 = \lim_{z \to x_0} \w(z)}$ exists belongs
to $\mathbb{R} \cup \{-\infty, \infty\}$.  Since
${\displaystyle \liminf_{z \to x_0}|\w(z)| < \infty}$, in fact
$\w_0 \in \mathbb{R}$.  
The fact that
the limits $\w_0$ and $d_S$ exist and are real
shows via (\ref{NSker})
that the limit
\begin{equation*}
\tilde{d} = \lim_{z \to x_0}K_\w(z,z)
\end{equation*}
exists and is real.
\nline

{\bf Proof of implication $2. \implies 3.$:}
Since $\w_0$ is real, it follows from
(\ref{phi}) that
the limit ${\displaystyle \lim_{z \to x_i}\varphi(\beta(z))}$
exists, is nonzero, and is finite in absolute value.
Therefore, by formula (\ref{NSker}), the limit
\begin{equation*}
\tilde{d}_S = \lim_{\zeta \to t_0}K_S^\Sch(\zeta,\zeta)
\end{equation*}
exists and is real, which, by
Theorem \ref{CJGS}, implies that the limit
$S_0$ exists and is unimodular, and that the limit
$S_1$ exists and is such that
\begin{equation}
\tilde{d}_S = S_1 \bar{S}_0 t_0 \in
\mathbb{R}.
\label{st}
\end{equation}
Let us rewrite the relation 
(\ref{Sformula}) between $S$ and
$\w$ as follows
\begin{equation*}
(S \circ \beta)(z) = (\beta \circ \w)(z) \quad \quad (z \in \rho(\w))
\end{equation*}
and differentiate each side with respect to $z$:
\begin{equation}
S^\prime(\beta(z)) \cdot \beta^\prime(z) =
\beta^\prime(\w(z)) \cdot \w^\prime(z).
\label{154}
\end{equation}
Making use of the relation $\beta^\prime(z) = 2i(z+i)^{-2}$, 
we may rewrite (\ref{154}) as
\begin{equation}
\w^\prime(z) = \left(\frac{\w(z) + i}{z + i}\right)^2
S^\prime(\zeta)
\label{wDerivative}
\end{equation}
which shows that $\w_1 := \lim_{z \to x_0}\w^\prime(z)$ exists
since both $\w_0$ and $S_1$ exist.
Using (\ref{st}), writing $t_0$ in terms of $x_0$ via
$\beta$, and taking limits in (\ref{wDerivative}),
we have
\begin{equation}
\w_1 = \left(\frac{\w_0 + i}{x_0 + i}\right)^2S_1
= \left(\frac{w_0 + i}{x_0 + i}\right)^2 \cdot
\frac{\tilde{d}_S}{\bar{S}_0} \cdot
\frac{x_0 + i}{x_0 - i}
=\frac{(\w_0 + i)^2}{x_0^2 + 1} \cdot S_0 \cdot \tilde{d}_S
\label{158}
\end{equation}
Since $S_0 = \beta(\w_0)$, 
(\ref{158}) yields
\begin{equation}
\w_1 = \frac{(\w_0 + i)^2}{x_0^2 + 1} \cdot
\frac{\w_0 - i}{\w_0 + i} \cdot \tilde{d}_S
= \frac{\w_0^2 + 1}{x_0^2 + 1} \cdot \tilde{d}_S
\label{w1dS}
\end{equation}
Since both $\w_0$ and $\tilde{d}_S$ are real, it follows that $\w_1$
is real as well.
\nline

{\bf Proof of implication $3. \implies 4.$:}
The associated function $S$ given by
(\ref{Sformula}) and presented in detail in
(\ref{Sdetail}) is such that the limit
\begin{equation}
S_0 := \lim_{\zeta \to t_0}S(\zeta) = \lim_{z \to x_0}
\frac{\w(z) - i}{\w(z) + i} = \frac{w_0 - i}{w_0 + i}
\label{E362}
\end{equation}
exists, is unimodular, and not equal to one.  
By relation (\ref{wDerivative}), the limit
\begin{equation}
S_1 := \lim_{\zeta \to t_0}S^\prime(\zeta) = \lim_{z \to x_0}
\left(\frac{z+i}{\w(z) + i}\right)^2\w^\prime(z) =
\left(\frac{x_0 + i}{\w_0 + i}\right)^2\w_1
\label{E363}
\end{equation}
exists.  
Together (\ref{E362}) and (\ref{E363})
show that $S_1\bar{S}_0t_0 \in \mathbb{R}$:
\begin{equation*}
S_1\bar{S}_0t_0 = \left(\frac{x_0 + i}{\w_0 + i}\right)^2\w_1 \cdot
\frac{\w_0 + i}{\w_0 - i}t_0
= \frac{(x_0 + i)^2}{\w_0^2 + 1}\w_1 \cdot \frac{x_0-i}{x_0+i}
= \frac{x_0^2 + 1}{\w_0^2 + 1}\w_1.
\end{equation*}
Hence, by Theorem \ref{CJGS}, the limit
\begin{equation}
\bar{S}_1 := \lim_{\zeta \to t_0}\frac{S(\zeta) - S_0}{\zeta - t_0}
\label{dCsense}
\end{equation}
exists and equals $S_1$.
Let us express
$
\frac{\w(z) - \w_0}{z - x_0}
$ 
in terms of the generalized
Schur function $S$ associated to $\w$ via (\ref{Sformula}),
and in terms of
$\zeta = \beta(z)$, $t_0 = \beta(x_0)$, and
$S_0 = \beta(\w_0)$.
We get
\begin{equation}
\frac{\w(z) - \w_0}{z - x_0} =
\frac{ i\frac{1 + S(\zeta)}{1 - S(\zeta)} -
i\frac{1 + S_0}{1 - S_0}}{
i\frac{1 + \zeta}{1 - \zeta}-
i\frac{1 + t_0}{1 - t_0}} =
\frac{(1 - \zeta)(1 - t_0)}{(1 - S(\zeta))(1 - S_0)}
\cdot \frac{S(\zeta) - S_0}{\zeta - t_0}.
\label{rearrangement}
\end{equation}
Sending $z \in \C^+$ to $x_0$ nontangentially results
in $\zeta \in \mathbb{D}$ being sent to $t_0$ nontangentially. 
Making use of (\ref{dCsense}),
we pass to limit in (\ref{rearrangement})
\begin{equation*}
\lim_{z \to x_0}
\frac{\w(z) - \w_0}{z - x_0} =
\left(\frac{1 - t_0}{1 - S_0}\right)^2 S_1,
\end{equation*}
which, upon using $x_0 = \beta^{-1}(t_0)$ and
$w_0 = \beta^{-1}(S_0)$,
may be rewritten as
\begin{equation}
\tilde{\w}_1 =
\frac{\left(1 - \frac{x_0 - i}{x_0 + i}\right)^2}{\left(
1 - \frac{\w_0 - i}{\w_0 + i}\right)^2}
S_1 =
\left( \frac{\w_0 + i}{x_0 + i} \right)^2 S_1.
\label{E370}
\end{equation}
Together (\ref{E370}) and
(\ref{158}) show that
\begin{equation}
\tilde{\w}_1 = \w_1 \in \mathbb{R}.
\label{tw1w1}
\end{equation}
\nline

{\bf Proof of implication $4. \implies 1.$:}
The associated function $S$ given by
(\ref{Sdetail}) is such that the limit
\begin{equation}
S_0 := \lim_{\zeta \to t_0}S(\zeta) = \lim_{z \to x_0}
\frac{\w(z) - i}{\w(z) + i} = \frac{w_0 - i}{w_0 + i}
\label{E373}
\end{equation}
exists, is unimodular, and not equal to one, 
since the limit $w_0$ exists and
is real.  By (\ref{rearrangement}),
\begin{equation}
\frac{\w(z) - \w_0}{z - x_0} =
\frac{(1 - \zeta)(1 - t_0)}{(1 - S(\zeta))(1 - S_0)}
\cdot \frac{S(\zeta) - S_0}{\zeta - t_0}.
\label{E374}
\end{equation}
Since the limit $\tilde{w}_1$ exists and is real,
(\ref{E373}) and (\ref{E374}) show that the limit
\begin{equation*}
\tilde{S}_1 = \lim_{\zeta \to t_0}
\frac{S(\zeta) - S_0}{\zeta - t_0}
\end{equation*}
exists and that
\begin{equation}
\left(\frac{1 - t_0}{1 - S_0}\right)^2
\tilde{S}_1 \in \mathbb{R}.
\label{E376}
\end{equation}
Using (\ref{E376}) and the fact that
$S_0 = \beta(\w_0)$ and $t_0 = \beta(x_0)$,
it follows that
\begin{equation*}
\tilde{S}_1\bar{S}_0 t_0 =
\left(\frac{x_0 + i}{\w_0 + i}\right)^2
\tilde{\w}_1 \cdot
\frac{\w_0 + i}{\w_0 - i} \cdot
\frac{x_0 + i}{w_0 - i}
= \frac{x_0 + i}{\w_0 + i} \tilde{\w}_1
\cdot \frac{x_0 - i}{\w_0 - i} \\
= \frac{x_0^2 + 1}{\w_0^2 + 1}\tilde{\w}_1,
\end{equation*}
which shows that $\tilde{S}_1 \bar{S}_0 t_0 \in
\mathbb{R}$.  By Theorem \ref{CJGS}, the limit
\begin{equation*}
d_S := \liminf_{\zeta \to t_0} K_S^\Sch(\zeta,\zeta)
\end{equation*}
satisfies $(-\infty < )\; d_S < \infty$.  
In view of this and the fact that $\w(z)$ has a finite
limit at $x_i$, formula (\ref{NSker}) shows that
\begin{equation*}
d := \liminf_{z \to x_0}K_\w(z,z)
\end{equation*}
is such that $-\infty < d < \infty$.
\nline

Finally, we show that when Statements 1--4 hold,
then $d = \tilde{d} = \w_1 = \tilde{\w}_1$.
From formula (\ref{tw1w1}) we have that
that $\tilde{w}_1 = \w_1$.  Since 
Statements 1--4 hold, the function $S$
associated to $\w$ satisfies the statements
in Theorem \ref{CJGS}, and so
$\tilde{d}_S = d_S$.  This combined with
(\ref{NSker}) and the fact that 
$0 < |\varphi(\beta(\w_0))| < \infty$ shows
that $\tilde{d} = d$.  To complete the
proof, we show that $\tilde{d} = \w_1$.
Towards this end, note that, by (\ref{NSker}),
it holds that
\begin{equation}
\tilde{d}_S = \frac{4}{\left|\w_0 + i\right|^2
\left|1 - \frac{x_0 - i}{x_0 + i}\right|^2}
\tilde{d}
= \frac{|x_0 + i|^2}{|\w_0 + i|^2}\tilde{d} 
= \frac{x_0^2 + 1}{\w_0^2 + 1} \tilde{d}.
\label{fstep}
\end{equation}
Comparing (\ref{w1dS}) and (\ref{fstep})
shows that indeed $\tilde{d} = \w_1$.
\end{proof}

\begin{Tm}
(Theorem \ref{Tm2})
Let $\ka \in \mathbb{Z}_+$,
$\w \in \N_\ka$ and $x_0 \in \mathbb{R}$.  Suppose that
the limit
\begin{equation*}
\lim_{z \to x_i}|\w(z)|
\end{equation*}
exists and equals infinity.  Then the following
statements are equivalent:
\begin{enumerate}
\item{\begin{equation*}
d := \liminf_{z \to x_0}K_{-1/\w}(z,z) \in \mathbb{R} \setminus
\{0\}
\end{equation*}}
\item{The nontangential limit
${\displaystyle \tilde{d} := \lim_{z \to x_0}K_{-1/\w}(z,z)}$ exists
and belongs to $\mathbb{R} \setminus \{0\}$.}
\item{The nontangential limit
\begin{equation*}
\w_{-1} := \lim_{z \to x_0}(z-x_i)\w(z)
\end{equation*}
exist and belongs to $\mathbb{R} \setminus \{0\}$.}
\item{The nontangential limit
\begin{equation*}
\tilde{\w}_{-1} := -\lim_{z \to x_0}
(z - x_i)^2\w^\prime(z)
\end{equation*}
exists and belongs to $\mathbb{R} \setminus \{0\}$.}
\end{enumerate}
In this case, $d = \tilde{d} = \w_{-1} = \tilde{\w}_{-1}$.
\label{CJGN2}
\end{Tm}
\begin{proof}
Set $\hat{\w}(z) := -\frac{1}{\w(z)}$.  Note that the limit
\begin{equation}
\lim_{z \to x_0} \hat{\w}(z) = 0
\label{hatw}
\end{equation}
exists and is zero.  The theorem then follows from
applying Theorem \ref{CJGN1} to $\hat{\w}$ and taking
into account (\ref{hatw}).
\end{proof}

For functions $\varphi \in \N_0$, slightly stronger versions of
the Carath\'eodory-Julia Theorems hold. 
To begin with, note that due to 
the fact that $\Im\; \varphi \geq 0$ on $\C^+$,
we have
\begin{equation*}
K_\varphi(z,z) = 
\frac{\varphi(z) - \varphi(z)^*}{z - \bar{z}} =
\frac{\Im\; \varphi}{\Im\; z} \geq 0.
\end{equation*}
Hence, for any point $x_0 \in \mathbb{R}$,
\begin{equation*}
\liminf_{z \to x_0} K_\varphi(z,z) \geq 0.
\end{equation*}
In view of Theorem \ref{CJGN1}, the preceding
argument proves the following
\begin{La}
Let $\varphi \in \N_0$ be a Nevanlinna function and
let $x_0 \in \mathbb{R}$.
Then the limit
\begin{equation*}
K_\varphi(x_0,x_0) := \lim_{z \to x_0} K_\varphi(z,z)
\end{equation*}
exists and, moreover,
\begin{equation*}
0 \leq K_\varphi(x_0,x_0) \leq \infty.
\end{equation*}
\label{CJN1}
\end{La}
Slight modifications to Theorems \ref{CJGN1} and
\ref{CJGN2} may be performed to prove the
following result.
\begin{La}
Let $\varphi \in \N_0$ be a Nevanlinna function and
let $x_0 \in \mathbb{R}$.
Then the limit
\begin{equation*}
\varphi_{-1}(x_0) := \lim_{z \to x_0} (z - x_0)\varphi(z)
\end{equation*}
exists and satisfies
\begin{equation}
-\infty < \varphi_{-1}(x_0) \leq 0.
\end{equation}
Moreover, the limit
\begin{equation*}
\lim_{z \to x_0} |z - x_0|^2 K_\varphi(z,z)
\end{equation*}
exists and equals $-\varphi_{-1}(x_0)$.
\label{CJN2}
\end{La}

\begin{Rk}
Lemmas \ref{CJN1} and \ref{CJN2} together with Theorems
\ref{CJGN1} and \ref{CJGN2} justify the claim in Section
2.3 that any function $\varphi \in \tN_0$
satisfies precisely one of $\bC_1$ -- $\bC_6$
if $\tilde{e}_i \neq 0$, and precisely one of
$\btC_1$ -- $\btC_6$ if $\tilde{e}_i = 0$.
\end{Rk}


\chapter{The Lyaponov Identity and Properties of $\Theta$}


\section{The Lyaponov Identity}

\begin{La}
Let the data $\Omega$, defined in (\ref{1.18}), be given.
Let $P$ be the Pick matrix associated to $\Omega$ as in
(\ref{1.19}) -- (\ref{120}) and
set $\ka = {\rm sq}_-P$.  Let $X$, $E$, and $C$ be defined as
in (\ref{r2.2}) and (\ref{r2.3}).  Then
the Lyaponov identity holds: 
\begin{equation}
PX - XP = E^*C - C^*E.
\label{LI}
\end{equation}
\end{La}
\begin{proof}
Let us introduce
the block partitionings
\begin{equation}
P^{-1} = \begin{bmatrix}\tilde{P}_{11} & \tilde{P}_{12} \\
 \tilde{P}_{21} & \tilde{P}_{22}\end{bmatrix}, \quad
X = \begin{bmatrix}X_{1} & 0 \\ 0 & X_{2}\end{bmatrix}, \quad
E = \begin{bmatrix}E_1 & E_2\end{bmatrix}, \quad
C = \begin{bmatrix}C_1 & C_2\end{bmatrix}
\label{e2.14}
\end{equation}
conformal with the partitioning of $P$ given in (\ref{1.19}).  
Note that $E_2 = 0$.
Substituting (\ref{e2.14}) into (\ref{LI}), we see that 
(\ref{LI}) reduces to essentially three equalities:
\begin{align}
P_{11}X_1 - X_1P_{11} &= E_1^*C_1 - C_1^*E_1,
\label{rLI}\\
P_{12}X_2 - X_1P_{12} &= E_1^*C_2 \label{sLI} \\
P_{22}X_2 - X_2P_{22} &= 0. \label{tLI}
\end{align}
The $(i,j)$ offdiagonal entry of left hand side of (\ref{rLI})
is given by
\begin{equation*}
\frac{\w_j - \w_i}{x_j - x_i}\cdot x_j - x_i \cdot
\frac{\w_j - \w_i}{x_j - x_i} = \w_j - \w_i,
\end{equation*}
and so is thus seen to be equal to the $(i,j)$ offdiagonal
entry of the right hand side.
The diagonal entries of both sides vanish, and 
hence (\ref{rLI}) is proved.
Identities (\ref{sLI}) and (\ref{tLI}) are similarly easily verified.
\end{proof}

\begin{La} 
Suppose that $P$ is invertible. Then
\begin{enumerate}
\item \text{The row vectors }
\begin{equation}
\tE := EP^{-1} = 
\begin{bmatrix}\tilde{e}_1 & \ldots & \tilde{e}_n\end{bmatrix}
\quad\text{and}\quad 
\tC := CP^{-1} =
\begin{bmatrix}\tilde{c}_1 & \ldots & \tilde{c}_n\end{bmatrix}
\label{e3.2}
\end{equation}
\text{ satisfy the Lyaponov identity}
\begin{equation}
XP^{-1} - P^{-1}X = \tE^*\tC -
\tC^*\tE.
\label{e3.4}
\end{equation}

\item \text{The offdiagonal entries $\tilde{p}_{ij}$ of $P^{-1}$
 are given by }
\begin{equation}
\tilde{p}_{ij} = \frac{\tilde{e}_i\tilde{c}_j -
\tilde{c}_i\tilde{e}_j}{x_i - x_j}.
\label{4.10a}
\end{equation}

\item \text{$\tilde{e}_i, \tilde{c}_i \in \mathbb{R}$ and
    $|\tilde{e}_i| + |\tilde{c}_i| > 0$, for $i = 1, \ldots, n$.}
\end{enumerate}
\label{L2-1}
\end{La}
\begin{proof}
To prove the first statement, multiply both sides 
in (\ref{LI}) by $P^{-1}$ on the left and on the right.
Setting $[\tilde{p}_{ij}]_{i,j=1}^n = P^{-1}$ allows equality
(\ref{e3.4}) to be rewritten entrywise as
\begin{equation*}
x_i\tilde{p}_{ij} - \tilde{p}_{ij}x_j =
\tilde{e}_i^*\tilde{c}_j -
\tilde{c}_i^*\tilde{e}_j
\end{equation*}
which implies
(\ref{4.10a}).
Since the entries of $E$ and $P$ are real, $\tE$
and $\tC$ have real entries as well.
Suppose that 
\begin{equation}
\tilde{e}_i = \tilde{c}_i = 0.
\label{onestar}
\end{equation}
Let ${\bf e}_i$
be the $i$-th column of the identity matrix $I_n$.  Then,
due to (\ref{e3.4}),
\begin{equation*}
(XP^{-1} - P^{-1}X){\bf e_i} = (\tE^*\tC -
\tC^*\tE){\bf e}_i,
\end{equation*}
or equivalently,
\begin{equation*}
-(x_iI-X)P^{-1}{\bf e}_i = \tE^*\tilde{c}_i - 
\tC^*\tilde{e}_i = 0.
\end{equation*}
Since all of the diagonal entries in the diagonal
matrix $x_iI - X$ are not zeros except for the $i$-th diagonal
entry, it follows that
\begin{equation}
P^{-1}{\bf e}_i = \alpha{\bf e}_i \quad \text{for some scalar }
\alpha \neq 0.
\label{e3.10}
\end{equation}
On the one hand, by (\ref{onestar}),
\begin{equation*}
\begin{bmatrix}\tC \\ \tE\end{bmatrix}{\bf e}_i = 
\begin{bmatrix}\tilde{c}_i \\ \tilde{e}_i\end{bmatrix} =
\begin{bmatrix}0 \\ 0\end{bmatrix},
\end{equation*}
while on the other hand, (\ref{e3.10}) and (\ref{e3.2}) imply
\begin{equation*}
\begin{bmatrix}\tC \\ \tE\end{bmatrix}
{\bf e}_i =
\begin{bmatrix}C \\ E\end{bmatrix}P^{-1}{\bf e}_i = 
\alpha \begin{bmatrix}C \\ E\end{bmatrix}{\bf e}_i =
\alpha \begin{bmatrix}c_i \\ e_i \end{bmatrix}.
\end{equation*}
If $e_i = 1$, then the right hand side is nonzero, and so
gives the desired contradiction.  If $e_i = 0$,
then $c_i = \xi_i \neq 0$,
and so in this case we also arrive at a contradiction, thus
proving that $|\tilde{e}_i| + |\tilde{c}_i| > 0$.
\end{proof}

\section{The Function $\Theta$}
In this section we study in more detail the function
\begin{equation}
\Theta(z) := 
\begin{bmatrix}\Theta_{11}(z) & \Theta_{12}(z) \\
\Theta_{21}(z) & \Theta_{22}(z)\end{bmatrix}
:= I_2 - i\begin{bmatrix}C \\ E\end{bmatrix}
(zI - X)^{-1}P^{-1}\begin{bmatrix}C^* & E^*\end{bmatrix}J.
\label{ThetaFormula2}
\end{equation}
\begin{La}
Let $\Theta$ be defined by (\ref{ThetaFormula2}).
If $i \in \{1, \ldots, \ell\}$ (i.e., if $x_i$ is a 
regular interpolation node),
then
\begin{equation*}
\lim_{z \to x_i}(z - x_i)\Theta(z) =
\lim_{z \to x_i}\begin{bmatrix}(z - x_i)\Theta_{11}(z) &
(z - x_i)\Theta_{12}(z) \\ (z - x_i)\Theta_{21}(z) &
(z - x_i)\Theta_{22}(z)\end{bmatrix} =
\begin{bmatrix}\w_i\tilde{e}_i & -\w_i\tilde{c}_i \\
\tilde{e}_i & -\tilde{c}_i\end{bmatrix}.
\end{equation*}
If $i \in \{\ell + 1, \ldots, n\}$ (i.e., if
$x_i$ is a singular interpolation node), then
\begin{equation*}
\lim_{z \to x_i}(z - x_i)\Theta(z) =
\lim_{z \to x_i}\begin{bmatrix}(z - x_i)\Theta_{11}(z) &
(z - x_i)\Theta_{12}(z) \\ (z - x_i)\Theta_{21}(z) &
(z - x_i)\Theta_{22}(z)\end{bmatrix} =
\begin{bmatrix}\xi_i\tilde{e}_i & -\xi_i\tilde{c}_i \\
0 & 0\end{bmatrix}.
\end{equation*}
\label{L42}
\end{La}
\begin{proof}
Making use of (\ref{e3.2}), we may write
(\ref{ThetaFormula2}) as
\begin{equation}
\Theta(z) = I_2 - i\begin{bmatrix}C \\ E\end{bmatrix}(zI - X)^{-1}
\begin{bmatrix}\tC^* & \tE^* \end{bmatrix}J.
\label{matlim}
\end{equation}
By (\ref{r2.2}),
\begin{equation}
\lim_{z \to x_i} (z - x_i)(zI - X)^{-1} = {\bf e}_i{\bf e}_i^*.
\label{6LIM}
\end{equation}
Hence
\begin{align}
\lim_{z \to x_i}(z - x_i)\Theta(z) &= -\lim_{z \to x_i}
i \begin{bmatrix}C \\ E\end{bmatrix}{\bf e}_i{\bf e}_i^*
\begin{bmatrix}\tC & \tE^*\end{bmatrix}J \nonumber \\
&=
\begin{bmatrix}C \\ E\end{bmatrix}{\bf e}_i {\bf e}_i^*
\begin{bmatrix}\tE^* & -\tC^*\end{bmatrix}
\label{r2.2a}
\end{align}
If $i \in \{1, \ldots, \ell\}$, then
$E\mathbf{e}_i = 1$ and $C\mathbf{e}_i = \w_i$, and so,
in view of (\ref{r2.2a}), we have
\begin{equation*}
\lim_{z \to x_i}(z - x_i)\Theta(z)
=\begin{bmatrix}\w_i \\ 1\end{bmatrix}
\begin{bmatrix}\tilde{e}_i & -\tilde{c}_i\end{bmatrix} =
\begin{bmatrix}\w_i \tilde{e}_i &
-\w_i \tilde{c}_i \\ \tilde{e}_i &
-\tilde{c}_i\end{bmatrix}.
\end{equation*}
On the other hand, if $i \in \{\ell + 1, \ldots, n\}$, then
$E\mathbf{e}_i = 0$ and $C\mathbf{e}_i = \xi_i$, and so,
in view of (\ref{r2.2a}), we have
\begin{equation*}
\lim_{z \to x_i}(z - x_i)\Theta(z)
=\begin{bmatrix}\xi_i \\ 0\end{bmatrix}
\begin{bmatrix}\tilde{e}_i & -\tilde{c}_i\end{bmatrix} =
\begin{bmatrix}\xi_i \tilde{e}_i &
-\xi_i \tilde{c}_i \\ 0 &
0\end{bmatrix}.
\end{equation*}
\end{proof}

\begin{Rk}
It holds that
$\lim_{z \to x_i} \frac{\Theta_{22}(z)}{\Theta_{21}(z)} = -\eta_i$,
where $\eta_i := \frac{\tilde{c}_i}{\tilde{e}_i}$ (as defined in
(\ref{inparticular})).
\end{Rk}

\begin{La}
If $P$ is invertible, then $\Theta$ belongs to the class
$\mathcal{W}_\ka$.
\label{LemmaRemark} 
\end{La}
\begin{proof}
First we shall show that
\begin{equation}
K_{\Theta,J}(z, \zeta) = \begin{bmatrix}C \\ E\end{bmatrix}
(zI - X)^{-1}P^{-1}(\bar{\zeta}I - X)^{-1}
\begin{bmatrix}C^* & E^*\end{bmatrix},
\label{KTJ}
\end{equation}
where the kernel $K_{\Theta,J}$ is given by (\ref{KTkernel}).
Indeed, by (\ref{ThetaFormula2}),
\begin{align}
J - \Theta(z)J\Theta(\zeta)^* =& J -
\left(I - i\begin{bmatrix}C \\ E\end{bmatrix}(zI - X)^{-1}
P^{-1}\begin{bmatrix}C^* & E^*\end{bmatrix}J\right) J \nonumber \\
& \times \left(I + iJ\begin{bmatrix}C \\ E\end{bmatrix}
P^{-1}(\bar{\zeta}I - X)^{-1}\begin{bmatrix}C^* & E^*\end{bmatrix}
\right) \nonumber \\
=& -i\begin{bmatrix}C \\ E\end{bmatrix}P^{-1}
(\bar{\zeta}I - X)^{-1}\begin{bmatrix}C^* & E^*\end{bmatrix} \nonumber
\\
& + i\begin{bmatrix}C \\ E\end{bmatrix}(zI - X)^{-1}P^{-1}
\begin{bmatrix}C^* & E^*\end{bmatrix} \nonumber \\
& - \begin{bmatrix}C \\ E\end{bmatrix}(zI - X)^{-1}P^{-1}
\begin{bmatrix}C^* & E^*\end{bmatrix}J
\begin{bmatrix}C \\ E\end{bmatrix}P^{-1}
(\bar{\zeta}I - X)^{-1}\begin{bmatrix}C^* & E^*\end{bmatrix} \nonumber
\\
=& -i\begin{bmatrix}C \\ E\end{bmatrix}(zI - X)^{-1}
(-P^{-1}(\bar{\zeta}I - X) + (zI - X)P^{-1} \nonumber \\
& + iP^{-1}\begin{bmatrix}C^* & E^*\end{bmatrix}J
\begin{bmatrix}C \\ E\end{bmatrix}P^{-1})
(\bar{\zeta}I - X)^{-1}\begin{bmatrix}C^* & E^*\end{bmatrix}.
\label{obs1}
\end{align}
Our next observation is that
\begin{align}
&-P^{-1}(\bar{\zeta}I - X) + (zI - X)P^{-1} + iP^{-1}
\begin{bmatrix}C^* & E^*\end{bmatrix}J\begin{bmatrix}C \\ 
E\end{bmatrix}P^{-1} \nonumber \\
&= -P^{-1}(\bar{\zeta}I - X) + (zI = X)P^{-1} +
P^{-1}(C^*E - E^*C)P^{-1} \nonumber \\
&= -P^{-1}(\bar{\zeta}I - X) + (zI - X)P^{-1} +
P^{-1}(PX - XP)P^{-1} \nonumber \\
&=(z - \bar{\zeta})P^{-1}.
\label{obs2}
\end{align}
Substituting (\ref{obs2}) into (\ref{obs1}) yields
\begin{equation}
J - \Theta(z)J\Theta(\zeta)^* = -i(z - \bar{\zeta})
\begin{bmatrix}C \\ E\end{bmatrix}(zI - X)^{-1}
P^{-1}(\bar{\zeta}I - X)^{-1}\begin{bmatrix}C^* & E^*
\end{bmatrix},
\label{Junitary}
\end{equation}
which, upon being divided by $i(z - \bar{\zeta})$,
implies (\ref{KTJ}).

Equality (\ref{Junitary}) implies that $\Theta$ is $J$-unitary on $\mathbb{R}$,
and (\ref{KTJ}) implies that
\begin{equation*}
{\rm sq}_-K_{\Theta,J} \leq {\rm sq}_-P = \ka.
\end{equation*}
To show that in fact ${\rm sq}_-K_{\Theta,J} = \ka$,
it suffices to 
check that the pair $(\begin{bmatrix}C \\ E\end{bmatrix}, X)$
is observable, i.e., that
\begin{equation}
\bigcap_{j = 0}^\infty {\rm Ker}\; \begin{bmatrix}C \\ E\end{bmatrix}
X^j = \{0\}.
\label{observability}
\end{equation}
Let us suppose that
(\ref{observability}) does not hold.  Then there exists
a nonzero column vector $y \in \C^n$ such that
\begin{equation}
\begin{bmatrix}C \\ E\end{bmatrix} X^j y =
\begin{bmatrix}0 \\ 0\end{bmatrix}
\quad\quad \text{for all $j \in \mathbb{Z}_+$.}
\label{obsA}
\end{equation}
Let $\begin{bmatrix}y_1^* & y_2^*\end{bmatrix}^* = y$
be a conformal partitioning of $y$ so that
$y_1 \in \C^\ell$ and $y_2 \in
\C^{n - \ell}$, and let $X$, $E$, and
$C$ be partitioned as in (\ref{e2.14}).
We may then
rewrite (\ref{obsA}) as
\begin{equation}
\begin{bmatrix}C_1 & C_2 \\ E_1 & 0\end{bmatrix}
\begin{bmatrix}X_1 & 0 \\ 0 & X_2\end{bmatrix}^j
\begin{bmatrix}y_1 \\ y_2\end{bmatrix} =
\begin{bmatrix}0 \\ 0\end{bmatrix}.
\label{obsA2}
\end{equation}
It is known (e.g., see \cite{BGR}) that for any $k \times k$ diagonal
matrix $A$ and for any $1 \times k$ row vector $v$ with
no zero entries, the pair
$(v, A)$ is observable.  Since neither $E_1$ nor $C_2$
has any entries equal to zero, the pairs
$(\begin{bmatrix}C_1 \\ E_1\end{bmatrix}, X_1)$ and 
$(\begin{bmatrix}C_2 \\ 0 \end{bmatrix}, X_2)$ are observable.
In view of (\ref{obs2}), the observability of the pairs 
implies, respectively, that
$y_1 = 0$ and $y_2 = 0$.
Thus the implication is that $y = 0$, which is a contradiction.
Therefore, the pair $(\begin{bmatrix}C \\ E\end{bmatrix}, X)$ is
indeed
observable.
\end{proof}
\begin{Rk}
The function $\Theta$ is $J$-unitary
on $\mathbb{R}$, and therefore, by the Schwarz reflection
principle, $\Theta(z)^{-1} = J\Theta(\bar{z})^*J$, from which it follows that
\begin{equation}
\Theta(z)^{-1} = I_2 + i\begin{bmatrix}C \\ E\end{bmatrix}
P^{-1}(zI - X)^{-1}\begin{bmatrix}C^* & E^*\end{bmatrix}J.
\label{e3.20}
\end{equation}
\end{Rk}
We shall have need of the following realization
formulas in the sequel.  The verification is
straightforward and is based on Lyaponov
identities (\ref{LI}) and (\ref{e3.4}), respectively.
\begin{Rk}
The following equalities hold for every choice of $z, \zeta \in \C \setminus \{x_1, \ldots, x_n\}:$
\begin{align*}
\frac{\Theta(\zeta)^{-*}J\Theta(z)^{-1} - J}{-i(z - \bar{\zeta})} &=
J\begin{bmatrix}C \\ E\end{bmatrix}(\bar{\zeta}I - X)^{-1}P^{-1}
(zI - X)^{-1}\begin{bmatrix}C^* & E^*\end{bmatrix}J, \\
\frac{J - \Theta(\zeta)^*J\Theta(z)}{-i(z - \bar{\zeta})} &=
J\begin{bmatrix}\tC \\ \tE\end{bmatrix}
(\bar{\zeta}I - X)^{-1}P(zI - X)^{-1}
\begin{bmatrix}\tC^* & \tE^*\end{bmatrix}J.
\end{align*}
\label{ThetaRemark}
\end{Rk}

\begin{Rk}
We now adopt some new conventions regarding
block matrix partitionings.  In the sequel, it will
sometimes be convenient to consider rearrangements of the interpolation
nodes $x_i$.  Thus the matrices $P$, $P^{-1}$, and $X$ undergo
permutation similiarities, and the columns of $E$ and $C$ are
permuted.  Nevertheless, we shall still use the symbols $P$,
$P^{-1}$, $X$, $E$, and $C$ to denote the transformed matrices.
Let $k \in \{1, \ldots, n-1\}$.  Note that the Lyaponov
Identity (\ref{LI}) holds regardless of how the interpolation
nodes are rearranged.
After rearrangment, we may thus consider the
following block partitionings:
\begin{equation}
P = \begin{bmatrix}P_{11} & P_{12} \\ P_{21} & P_{22}\end{bmatrix},
\quad
P^{-1} = \begin{bmatrix}\tilde{P}_{11} & \tilde{P}_{12} \\
 \tilde{P}_{21} & \tilde{P}_{22}\end{bmatrix},
\label{PP}
\end{equation}
and
\begin{equation}
X = \begin{bmatrix}X_{1} & 0 \\ 0 & X_{2}\end{bmatrix}, \quad
E = \begin{bmatrix}E_1 & E_2\end{bmatrix}, \quad
C = \begin{bmatrix}C_1 & C_2\end{bmatrix},
\label{XEC}
\end{equation}
where $P_{11}, \tilde{P}_{11}, X_1 \in \C^{k \times k}$, and
$E_1, C_1 \in \C^{1 \times k}$.
\label{BlockRemark}
\end{Rk}
In the remainder of this chapter, we assume that the
block partitionings are given by (\ref{PP}) and (\ref{XEC}),
with $k$ arbitrary in $\{1, \ldots, n-1\}$.
\begin{La}
Let us suppose that $P_{11}$ is invertible and that
${\rm sq}_-P_{11} = \ka_1 \leq \ka$.  Then $\tilde{P}_{22}$
is invertible, ${\rm sq}_-\tilde{P}_{22} = \ka - \ka_1$,
and the functions
\begin{align}
\Theta^{(1)}(z) = I_2 - i\begin{bmatrix}C_1 \\ E_1\end{bmatrix}
(zI - X_1)^{-1}P_{11}^{-1}\begin{bmatrix}C_1^* & E_1^*\end{bmatrix}J,
\label{e3.23} \\
\tilde{\Theta}^{(2)}(z) = I_2 - i\begin{bmatrix}\tC_2 \\
\tE_2\end{bmatrix}\tilde{P}_{22}^{-1}(zI - X_2)^{-1}
\begin{bmatrix}\tC_2^* & \tE_2^*\end{bmatrix}J
\label{e3.24}
\end{align}
belong to $\mathcal{W}_{\ka_1}$ and $\mathcal{W}_{\ka - \ka_1}$, respectively.
Furthermore, $\Theta$ admits the factorization
\begin{equation}
\Theta(z) = \Theta^{(1)}(z)\tilde{\Theta}^{(2)}(z).
\label{e3.25}
\end{equation}
\label{InvLemma}
\end{La}
\begin{proof}
To prove the first statement we use the standard
Schur complement argument.  Since both $P$ and $P_{11}$ are
invertible, the Schur complement of $P_{11}$ in $P$ is invertible
and has $\ka - \ka_1$ negative eigenvalues.  This matrix
is precisely $P_{22} - P_{21}P_{11}^{-1}P_{12}$.  The block
$\tilde{P}_{22}$ in $P^{-1}$ equals
$(P_{22} - P_{21}P_{11}^{-1}P_{12})^{-1}$, and so also has
$\ka - \ka_1$ negative eigenvalues.  Let us set for short
\begin{equation*}
R(z) = -\begin{bmatrix}C_1 \\ E_1\end{bmatrix}(zI - X_1)^{-1}, \quad
\tilde{R}(z) = -\begin{bmatrix}\tC_2 \\ \tE_2\end{bmatrix}
\tilde{P}_{22}^{-1}(zI - X_2)^{-1}.
\end{equation*}
The Lyaponov identities
\begin{equation*}
P_{11}X_{1} - X_{1}P_{11} = E_1^*C_1 - C_1^*E_1 
\quad\text{and}\quad
X_2\tilde{P}_{22}^{-1} - \tilde{P}_{22}^{-1}X_2 =
\tE_2^*\tC_2 - \tC_2^*\tE_2
\end{equation*}
are contained in Lyaponov identities (\ref{LI})
and (\ref{e3.4}).  These
allow us to establish  the following
realization formulas:
\begin{equation}
K_{\Theta^{(1)},J}(z, \zeta) = R(z)P_{11}^{-1}R(\zeta)^*
\quad\text{and}\quad
K_{\tilde{\Theta}^{(2)},J} = \tilde{R}(z)\tilde{P}_{22}
\tilde{R}(\zeta)^*.
\label{realizationformulas}
\end{equation}
The calculations that establish these formulas are
straightforward and parallel those performed in Lemma
\ref{LemmaRemark}.
Formulas (\ref{realizationformulas}) show that the rational functions
$\Theta^{(1)}$ and $\tilde{\Theta}^{(2)}$ are $J$-unitary on
$\mathbb{R}$ and that
\begin{equation}
{\rm sq}_-K_{\Theta^{(1)},J} \leq {\rm sq}_-P_{11} = \ka_1
\quad\text{and}\quad
{\rm sq}_-K_{\tilde{\Theta}^{(2)},J} \leq {\rm sq}_-
\tilde{P}_{22} = \ka - \ka_1.
\label{e3.29}
\end{equation}
Under the assumption that the factorization formula (\ref{e3.25}) holds, 
it follows that
\begin{equation*}
K_{\Theta, J}(z, \zeta) = K_{\Theta^{(1)},J}(z, \zeta) +
\Theta^{(1)}(z) K_{\tilde{\Theta}^{(2)},J}(z, \zeta)
\Theta^{(1)}(\zeta)^*
\end{equation*}
and so
\begin{equation*}
\ka = {\rm sq}_-K_{\Theta, J} \leq
{\rm sq}_-K_{\Theta^{(1)},J} +
{\rm sq}_-K_{\tilde{\Theta}^{(2)},J},
\end{equation*}
which, combined with inequalities (\ref{e3.29}), implies
\begin{equation*}
{\rm sq}_-K_{\Theta^{(1)},J} = \ka_1 
\quad\text{and}\quad
{\rm sq}_-K_{\tilde{\Theta}^{(2)},J} = \ka - \ka_1.
\end{equation*}
Only the factorization remains to be proven.  Towards this end
we shall systematically make use of the well-known equality
\begin{equation}
P^{-1} = \begin{bmatrix}P_{11}^{-1} & 0 \\ 0 & 0\end{bmatrix}
+ \begin{bmatrix}-P_{11}^{-1}P_{12} \\ I_{k}\end{bmatrix}
\tilde{P}_{22}
\begin{bmatrix}-P_{21}P_{11}^{-1} & I_k\end{bmatrix}.
\label{wellknown}
\end{equation}
This formula allows us to write
\begin{equation*}
\begin{bmatrix}\tC_2 \\ \tE_2\end{bmatrix} =
\begin{bmatrix}C \\ E\end{bmatrix}P^{-1}\begin{bmatrix}0 \\
 I_k\end{bmatrix} =
\begin{bmatrix}C \\ E\end{bmatrix}
\begin{bmatrix}-P_{11}^{-1}P_{12} \\ I\end{bmatrix}
\tilde{P}_{22}.
\end{equation*}
Then (\ref{e3.24}) may be expressed as
\begin{equation}
\tilde{\Theta}^{(2)}(z) = I_2 - i
\begin{bmatrix}C \\ E\end{bmatrix}
\begin{bmatrix}-P_{11}^{-1}P_{12} \\ I\end{bmatrix}
(zI - X_2)^{-1} \begin{bmatrix}\tC_2^* &
\tE_2^*\end{bmatrix} J.
\label{tildeT2}
\end{equation}
On the other hand, substituting (\ref{wellknown}) into
(\ref{ThetaFormula2}) and taking into account (\ref{e3.23})
yields
\begin{align*}
\Theta(z) &= I_2 -i\begin{bmatrix}C \\ E\end{bmatrix}
(zI-X)^{-1}P^{-1}\begin{bmatrix}C^* & E^*\end{bmatrix}J \\
&= \Theta^{(1)}(z) -i\begin{bmatrix}C \\ E\end{bmatrix}
(zI-X)^{-1}\begin{bmatrix}-P_{11}^{-1}P_{12} \\ I_k\end{bmatrix}
\tilde{P}_{22}\begin{bmatrix}-P_{21}P_{11}^{-1} & I_k\end{bmatrix}
\begin{bmatrix}C^* & E^*\end{bmatrix}J \\
&= \Theta^{(1)}(z) + i\begin{bmatrix}C \\ E\end{bmatrix}
\begin{bmatrix}-P_{11}^{-1} \\ I_k\end{bmatrix}
\begin{bmatrix}\tC_2^* & \tE_2^*\end{bmatrix}J.
\end{align*}
So (\ref{e3.25}) is equivalent to
\begin{equation}
\tilde{\Theta}^{(2)}(z)
= I_2 - i\Theta^{(1)}(z)^{-1}\begin{bmatrix}C \\ E\end{bmatrix}
(zI-X)^{-1}\begin{bmatrix}-P_{11}^{-1}P_{12} \\ I_k\end{bmatrix}
\begin{bmatrix}\tC_2^* & \tE_2^*\end{bmatrix}J. 
\label{tildeT3}
\end{equation}
To complete the proof, it suffices to show that
\begin{equation}
\Theta^{(1)}(z)^{-1}\begin{bmatrix}C \\ E\end{bmatrix}
(zI-X)^{-1}\begin{bmatrix}-P_{11}^{-1}P_{12} \\ I_k\end{bmatrix} =
\begin{bmatrix}C \\ E\end{bmatrix}\begin{bmatrix}
-P_{11}^{-1}P_{12} \\ I_k\end{bmatrix}(zI - X_2)^{-1},
\label{e3.36}
\end{equation}
since substituting (\ref{e3.36}) into (\ref{tildeT3}) yields
(\ref{tildeT2}).
The formula for $\Theta^{(1)}(z)^{-1}$ follows from the
reflection principle:
\begin{equation}
\Theta^{(1)}(z)^{-1} = J\Theta^{(1)}(\bar{z})^*J
= I_2 + i\begin{bmatrix}C_1 \\ E_1\end{bmatrix}
P_{11}^{-1}(zI - X_1)^{-1}\begin{bmatrix}C_1^* & E_1^*\end{bmatrix}J.
\label{E423}
\end{equation}
Now compare the top block entries in the Lyaponov equation (\ref{LI}):
\begin{equation*}
\begin{bmatrix}P_{11}X_1 & P_{12}X_2\end{bmatrix} -
\begin{bmatrix}X_1P_{11} & X_1P_{12}\end{bmatrix} =
E_1^*C - C_1^*E.
\end{equation*}
Multiplying each side by $(zI - X_1)^{-1}$ on the left and
by $(zI - X)^{-1}$
on the right yields
\begin{equation*}
(zI - X_1)^{-1}(E_1^*C - C^*E)(zI - X)^{-1}
=(zI - X_1)^{-1}\begin{bmatrix}P_{11} & P_{12}\end{bmatrix} -
\begin{bmatrix}P_{11} & P_{12}\end{bmatrix}(zI - X_1)^{-1}.
\end{equation*}
Therefore,
\begin{align*}
&\Theta^{(1)}(z)^{-1}\begin{bmatrix}C \\ E\end{bmatrix}
(zI - X)^{-1}\begin{bmatrix}-P_{11}^{-1}P_{12} \\ I\end{bmatrix}
\\
=& \begin{bmatrix}C \\ E\end{bmatrix}(zI - X)^{-1}
\begin{bmatrix}-P_{11}^{-1}P_{12} \\ I\end{bmatrix}
-\begin{bmatrix}C_1 \\ E_1\end{bmatrix}
\begin{bmatrix}I & P_{11}^{-1}P_{12}\end{bmatrix}(zI - X)^{-1}
\begin{bmatrix}-P_{11}^{-1}P_{12} \\ I\end{bmatrix}
\\
=& -\begin{bmatrix}C_1 \\ E_1\end{bmatrix}(zI - X_1)^{-1}
P_{11}^{-1}P_{12} +\begin{bmatrix}C_2 \\ E_2\end{bmatrix}
(zI - X_2)^{-1} \\
&-\begin{bmatrix}C_1 \\ E_1\end{bmatrix}\left(
P_{11}^{-1}P_{12}(zI - X_2)^{-1} - (zI - X_1)^{-1}P_{11}^{-1}P_{12}
\right)
\\
=& \begin{bmatrix}C \\ E\end{bmatrix}
\begin{bmatrix}-P_{11}^{-1}P_{12} \\ I\end{bmatrix}
(zI - X_2)^{-1},
\end{align*}
which proves (\ref{e3.36}) and thus completes the proof.
\end{proof}

\chapter{Fundamental Matrix Inequality}
In the 1970's, V.~Potapov suggested a general method
to solve classical interpolation problems, the method
of the Fundamental Matrix Inequality (FMI)
(see \cite{Ame}).  This
method was subsequently developed by other mathematicians
and applied to very general modern interpolation problems
that had been
advanced far beyond their classical roots
(e.g., see \cite{And}, \cite{Kup}).  The first attempt
to adapt this method to generalized Schur functions was
undertaken in \cite{Gol}.  Here
we make an appropriate adaptation of Potapov's 
method of the Fundamental Matrix
Inequality in order to characterize the solution set 
of Problem \ref{Prob3}.

We start with some
some simple observations.
For the proof, see \cite{BK}.
\begin{Pn}
Let $K(z, \zeta)$ be a sesqui-analytic kernel defined on 
$\Omega$ 
such that ${\rm sq}_-K = \ka$.
Then
\begin{enumerate}
\item \text{For every choice of an integer $p$, of a 
Hermitian $p \times p$ matrix $A$,}\\
\text{
and of a $p \times 1$ vector-valued function $G$,}
\begin{equation*}
{\rm sq}_-
\begin{bmatrix}A & G(z) \\ G(\zeta)^* & K(z, \zeta)\end{bmatrix}
\leq \ka + p.
\end{equation*}
\item \text{If $\lambda_1, \ldots, \lambda_n$ are points
in $\Omega$ and if}
\begin{equation*}
A = \left[K(\lambda_j, \lambda_i)\right]_{i,j=1}^p
\quad\text{and}\quad G(z) = 
\begin{bmatrix}K(z, \lambda_1) \\ \vdots \\ K(z, \lambda_p)
\end{bmatrix},
\end{equation*}
\text{then}
\begin{equation*}
{\rm sq}_-
\begin{bmatrix}A & G(z) \\ G(\zeta)^* & K(z, \zeta)\end{bmatrix}
= \ka.
\end{equation*}
\end{enumerate}
\label{p2.1}
\end{Pn}

\section{The Fundamental Matrix Inequality}
\begin{Tm}
Let $\w$
be a function meromorphic on $\mathbb{C}^+$.
Then for $\w$ to be a solution
of Problem \ref{Prob3}, it is necessary and sufficient that the kernel
\begin{equation}
\mathbf{K}_\w(z, \zeta) :=
\begin{bmatrix}
P &
(z I_n - X)^{-1}(\w(z) E^* - C^*) \\
(\w(\zeta)^* E - C)(\bar{\zeta} I_n - X)^{-1} &
K_\w(z, \zeta)
\end{bmatrix}
\label{r5.1} 
\end{equation}
have $\ka$ negative squares on $\rho(\w) \times \rho(\w)$:
\begin{equation}
{\rm sq}_-\mathbf{K}_\w(z, \zeta) = \ka.
\label{r5.2}
\end{equation}
\end{Tm}
\begin{proof}
We prove the necessity part of the theorem here, postponing
the proof of the sufficienty part until Chapter VII.
For $\w$ to be a solution
of Problem \ref{Prob3}, $\w$ must belong to the
class $\N_{\ka^\prime}$ for some $\ka^\prime \leq \ka$
and satisfy interpolation conditions 
(\ref{1.19}) and (\ref{1.30})
at all but
precisely $\ka - \ka^\prime$ interpolation nodes.
Let us first consider the case where $\w$ belongs to
$\N_\ka$ and so satisfies all of the interpolation
conditions.  In this case,
${\rm sq}_-K_\w = \ka$ and so by the second statement
in Proposition \ref{p2.1}, the kernel
\begin{equation*}
\mathbf{K}_\w^{(1)}(z, \zeta):=
\begin{bmatrix}
K_\w(z_1, z_1) & \ldots & K_\w(z_n, z_1) & K_\w(z, z_1) \\
\vdots & & \vdots & \vdots \\
K_\w(z_1, z_n) & \ldots & K_\w(z_n, z_n) & K_\w(z, z_n) \\
K_\w(z_1, \zeta) & \ldots & K_\w(z_n, \zeta) & 
 K_\w(z, \zeta)
\end{bmatrix}
\end{equation*}
has $\ka$ negative squares on 
$\rho(\w) \times \rho(\w)$
for every choice of points $z_1, \ldots, z_n \in \rho(\w)$.
Let
$\mathbf{z}, \mathbf{x}, P^\w(\mathbf{z}),
B(\mathbf{z}), P^\w(\mathbf{x}), A_\Omega$, and $A_\w$ be defined
as they were in Section 1.5.  Define
\begin{equation*}
\tilde{B}(\mathbf{z}) = \begin{bmatrix}B(\mathbf{z}) & 0 \\ 0 &
  1\end{bmatrix}, \quad
\tA_\Omega = \begin{bmatrix}A_\Omega & 0 \\ 0 & 1\end{bmatrix}, \quad
\tA_\w = \begin{bmatrix}A_\w & 0 \\ 0 & 1\end{bmatrix}.
\end{equation*}
The kernel $\tilde{B}({\bf z})^* \cdot \mathbf{K}_\w^{(1)}(z,
\zeta) \cdot \tilde{B}(\mathbf{z})$ then has $\ka$ negative squares
on $\rho(\w) \times \rho(\w)$.  The limit of these kernels
(as $\mathbf{z} \to \mathbf{x}$) is the
kernel of the form
\begin{equation*}
\mathbf{K}_\w^{(2)}(z, \zeta) = \begin{bmatrix}
P^\w(\mathbf{x}) & (zI - X)^{-1}(\w(z)E^* - C^*) \\
(\w(\zeta)^*E - C)(\bar{\zeta}I - X)^{-1} & K_\w(z, \zeta)
\end{bmatrix}.
\end{equation*}
Since the kernel $\mathbf{K}_\w^{(2)}$ is the limit of a family of kernels, each
of which has $\ka$ negative squares, it follows that ${\rm
  sq}_-\mathbf{K}_\w^{(2)} \leq \ka$.  Note that
\begin{equation}
\tA_\Omega \cdot \mathbf{K}_\w(z, \zeta) \cdot \tA_\Omega
= \tA_\w \cdot \mathbf{K}_\w^{(2)}(z, \zeta) \cdot 
\tA_\w + \begin{bmatrix}A_\Omega \cdot P \cdot A_\Omega - 
A_\w \cdot P^\w(\mathbf{x}) \cdot A_\w & 0 \\ 0 & 0\end{bmatrix},
\label{r5.3}
\end{equation}
where, as was shown in Section 1.5,
\begin{equation*}
A_\Omega \cdot P \cdot A_\Omega -
A_\w \cdot P^\w(\mathbf{x}) \cdot A_\w =
\begin{bmatrix}
D_1 & 0\\
0 & D_2
\end{bmatrix}
\end{equation*}
with
\begin{align*}
D_1 &=
\begin{bmatrix}
\gamma_1 - \w^\prime(x_1) & & 0\\
 & \ddots & \\
0 & & \gamma_\ell - \w^\prime(x_\ell)
\end{bmatrix},
\\
D_2 &=
\begin{bmatrix}
\w_{-1}^{-1}(x_{\ell + 1}) - \xi_{\ell + 1}^{-1} & & 0\\
 & \ddots & \\
0 & & \w_{-1}^{-1}(x_{n + 1}) - \xi_n^{-1} \\
\end{bmatrix}.
\end{align*}
Since, by assumption, $\w$ satisfies all the interpolation
conditions, it follows that the second term in the right
hand side of (\ref{r5.3}) is positive semidefinite.
Therefore
\begin{equation}
{\rm sq}_-\mathbf{K}_\w \leq {\rm sq}_-\mathbf{K}_\w^{(2)} \leq \ka.
\label{r5.4}
\end{equation}
Since $\mathbf{K}_\w$ contains the kernel $K_\w$ as a
principal submatrix, ${\rm sq}_-\mathbf{K}_\w \geq 
{\rm sq}_-K_\w = \ka$, which combined with (\ref{r5.4}) yields
(\ref{r5.2}).  
We now turn to the general case, supposing that $\w \in
\N_{\ka^\prime}$ for some $\ka^\prime \leq \ka$.
The interpolation conditions 
(\ref{1.29}) and (\ref{1.30})
are satisfied at all but
$k = \ka - \ka^\prime$ nodes $x_i$.  Let us denote
the set of indices associated to these $k$ nodes
by $J_k$.
Explicitly, only the following interpolation conditions are satisfied:
\begin{align}
\w(x_i) = \w_i 
\quad&\text{and}\quad
\w^\prime(x_i) \leq \gamma_i && \text{for}\quad i \in
\mathbb{N}_\ell \cap J_k \label{r5.5} \\
\frac{-1}{\w_{-1}(x_j)} &\leq \frac{-1}{\w_j} && \text{for}\quad j \in (\mathbb{N}_n
\setminus
\mathbb{N}_\ell) \cap J_k.
\label{r5.6}
\end{align}
Let us rearrange the interpolation nodes so that the indices
of the first $k$ interpolation nodes belong to $J_k$.
Recall the conformal partitionings of $P$, $E$, and $C$ introduced in
(\ref{PP}) and (\ref{XEC}).
Let us set
\begin{equation}
F_i(z) = (z I - X_i)^{-1}(\w(z)E_i^* - C_i^*) \qquad
\text{for}\quad i = 1,2,
\label{r5.7}
\end{equation}
so that
\begin{equation*}
\begin{bmatrix} F_1(z) \\ F_2(z)\end{bmatrix} =
(z I_n - X)^{-1}(\w(z) E^* - C^*).
\end{equation*}
Now the matrix $P_{11}$ is the Pick matrix of the truncated
interpolation problem (Problem \ref{Prob3}) with
interpolation conditions (\ref{r5.5}) and (\ref{r5.6}).  
Applying the first part of the proof to this truncated
problem shows that the kernel
\begin{equation}
\widetilde{\mathbf{K}}_\w(z, \zeta) :=
\begin{bmatrix}
P_{11} & F_1(z) \\
F_1(\zeta)^* & K_\w(z, \zeta)
\end{bmatrix}
\label{r5.8}
\end{equation}
has $\ka^\prime$ negative squares on $\rho(\w)$.
Applying the first statement in Proposition \ref{p2.1} to
\begin{equation}
K(z, \zeta) = \tilde{\mathbf{K}}_\w(z, \zeta), \quad
G(z) = \begin{bmatrix}P_{21} & F_2(z)\end{bmatrix},
\quad\text{and}\quad
A = P_{22},
\label{r5.9}
\end{equation}
we conclude that
\begin{equation}
{\rm sq}_-\begin{bmatrix}P_{22} & G(z) \\
G(\zeta)^* & \tilde{\mathbf{K}}_\w(z, \zeta)
 \end{bmatrix}
\leq {\rm sq}_-\widetilde{\mathbf{K}}_\w + k =
\ka^\prime + (\ka - \ka^\prime) =
\ka.
\label{r5.10}
\end{equation}
On account of (\ref{r5.7}) and (\ref{r5.9}), 
\begin{equation}
\begin{bmatrix}
P_{22} & G(z)\\
G(\zeta)^* & \tilde{\mathbf{K}}_\w(z, \zeta)
\end{bmatrix} =
\begin{bmatrix}
P_{22} & P_{21} & F_2(z)\\
P_{12} & P_{11} & F_1(z)\\
F_2(\zeta)^* & F_1(\zeta)^* & K_\w(z, \zeta)
\end{bmatrix}.
\label{r5.11}
\end{equation}
Together (\ref{r5.8}) and (\ref{r5.11}) show that
\begin{equation*}
\mathbf{K}_\w(z, \zeta) = U
\begin{bmatrix}
P_{22} & G(z) \\
G(\zeta)^* & \tilde{\mathbf{K}}_\w(z, \zeta)
\end{bmatrix}
U^*,
\text{ where }
U =
\begin{bmatrix}
0 & I_{\ka - k} & 0 \\
I_k & 0 & 0\\
0 & 0 & 1
\end{bmatrix},
\end{equation*}
which, due to (\ref{r5.10}), implies that ${\rm sq}_-\mathbf{K}_\w
\leq \ka$.  On the other hand, $P$ rests as a principal
submatrix in $\mathbf{K}_\w$, and so
${\rm sq}_-\mathbf{K}_\w \geq {\rm sq}_-P = \ka$.
This completes the proof of the necessity part of the theorem.
For the proof of sufficiency, see Chapter VII.
\end{proof}

\section{Parameterization of the Solution Set}
When $P$ is invertible, all functions for which (\ref{r5.2})
is satisfied
can be described in terms of a linear fractional transformation.
\begin{Tm}
Suppose that the Pick matrix P is invertible and let 
$\Theta = [\Theta_{ij}]$ be the $2 \times 2$ matrix
valued function defined in (\ref{ThetaFormula2}).  
A function $\w$
meromorphic in $\mathbb{C}^+$ is subject to (\ref{r5.2})
if and only if it is of the form
\begin{equation}
\w(z) = \mathbf{T}_\Theta[\varphi](z):=
\frac{\Theta_{11}(z)\varphi(z) + \Theta_{12}(z)}{
\Theta_{21}(z)\varphi(z) + \Theta_{22}(z)}
\label{r5.12}
\end{equation}
for some function $\varphi \in \tN_0$.
\end{Tm}
\begin{proof}
Let $\mathbf{S}$ be the Schur complement of $P$
in the kernel $\mathbf{K}_\w$ defined in (\ref{r5.1}):
\begin{equation}
\mathbf{S}(z, \zeta) := K_\w(z, \zeta) - 
(\w(\zeta)^*E - C)(\bar{\zeta}I - X)^{-1} P^{-1}
(zI - X)^{-1}(\w(z)E^* - C^*).
\label{r5.13}
\end{equation}
The equalities
\begin{equation*}
K_\w(z, \zeta) := \frac{\w(z) - \w(\zeta)^*}{z - \zeta} =
\begin{bmatrix}\w(\zeta)^* & 1\end{bmatrix}
\frac{J}{i(z - \bar{\zeta})}
\begin{bmatrix}\w(z) \\ 1\end{bmatrix}
\end{equation*}
and
\begin{equation*}
\w(\zeta)^*E - C =
i\begin{bmatrix}\w(\zeta)^* & 1\end{bmatrix}J
\begin{bmatrix}C \\ E\end{bmatrix},
\end{equation*}
where $J$ is the signature matrix defined in (\ref{r2.2}),
allow us to write
(\ref{r5.13}) as
\begin{align*}
\mathbf{S}(z, \zeta) =&
\begin{bmatrix}\w(\zeta)^* & 1\end{bmatrix}\left\{
\frac{J}{i(z - \bar{\zeta})} - J\begin{bmatrix}C \\ E\end{bmatrix}
(\bar{\zeta}I - X)^{-1}P^{-1}\right. \nonumber \\
& \left. \times (zI - X)^{-1}\begin{bmatrix}C^* & E^*\end{bmatrix}
J\right\}\begin{bmatrix}\w(z) \\ 1\end{bmatrix},
\end{align*}
which, due to Remark \ref{ThetaRemark}, is the same as
\begin{equation*}
\mathbf{S}(z, \zeta) = \begin{bmatrix}\w(\zeta)^* & 1\end{bmatrix}
\frac{\Theta(\zeta)^{-*}J\Theta(z)^{-1}}{i(z - \bar{\zeta})}
\begin{bmatrix}\w(z) \\ 1\end{bmatrix}.
\end{equation*}
The standard Schur complement argument reveals
\begin{equation*}
{\rm sq}_-\mathbf{K}_\w = {\rm sq}_-P +
{\rm sq}_-\mathbf{S} =
\ka + {\rm sq}_-\mathbf{S},
\end{equation*}
and so (\ref{r5.2}) holds
if and only if the kernel
$\mathbf{S}$ is positive definite on
$\rho(\w) \times
\rho(\w)$, i.e.,
\begin{equation}
\begin{bmatrix}\w(\zeta)^* & 1\end{bmatrix}
\frac{\Theta(\zeta)^{-*}J\Theta(z)^{-1}}{
i(z - \bar{\zeta})}
\begin{bmatrix}\w(z) \\ 1\end{bmatrix} \succeq 0.
\label{r5.14}
\end{equation}
It remains to show that (\ref{r5.14}) holds 
for meromorphic $\w$ if and only if
$\w$ is of the form (\ref{r5.12}).  
To prove this, let us first assume that
$\w$ satisfies (\ref{r5.14}),
and introduce meromorphic functions $u$ and $v$:
\begin{equation}
\begin{bmatrix}u(z) \\ v(z)\end{bmatrix} :=
\Theta(z)^{-1} \begin{bmatrix}\w(z) \\ 1\end{bmatrix}.
\label{r5.15}
\end{equation}
Inequality (\ref{r5.14})
may then be written in terms of $u$ and $v$ as
\begin{equation}
\begin{bmatrix}u(\zeta)^* & v(\zeta)^*\end{bmatrix}
\frac{J}{i(z - \bar{\zeta})}
\begin{bmatrix}u(z) \\ v(z)\end{bmatrix} =
\frac{u(z)v(\zeta)^* - u(\zeta)^*v(z)}{z - \bar{\zeta}} \succeq 0.
\label{r5.16}
\end{equation}
By (\ref{r5.15}), $u$ and $v$ are analytic on
$\rho(\w)$.  We shall show that
\begin{equation}
v(z) \neq 0 \quad\text{for every}\quad z \in \rho(\w),
\label{r5.17}
\end{equation}
unless $v \equiv 0$.
To this end, let us suppose that there are points $z, \zeta \in \rho(\w)$
such that $v(z) \neq 0$ and $v(\zeta) = 0$.  Positivity
of the kernel in (\ref{r5.16}) implies that
the $2 \times 2$ matrix
\begin{equation*}
\begin{bmatrix}
\frac{u(\zeta)v(\zeta)^* - u(\zeta)^*v(\zeta)}{\zeta - \bar{\zeta}} &
\frac{u(z)v(\zeta)^* - u(\zeta)^*v(z)}{z - \bar{\zeta}} \\
\frac{u(\zeta)v(z)^* - u(z)^*v(\zeta)}{\zeta - \bar{z}} &
\frac{u(z)v(z)^* - u(z)^*v(z)}{z - \bar{z}}
\end{bmatrix} =
\begin{bmatrix}
0 & -\frac{u(\zeta)^*v(z)}{z - \bar{\zeta}} \\
\frac{u(\zeta)v(z)^*}{\zeta - \bar{z}} & 
\frac{u(z)v(z)^* - u(z)^*v(z)}{z - \bar{z}}
\end{bmatrix}
\end{equation*}
is positive semidefinite.  
Then $u(\zeta)^*v(z) = 0$,
which in turn implies $u(\zeta) = 0$.  This forces
${\rm det}\;\Theta(\zeta)^{-1} = 0$ on account of (\ref{r5.15}),
which is a contradiction, since $\Theta^{-1}$ is invertible,
except on $\{x_1, \ldots, x_n\}$.  This proves (\ref{r5.17}).
Now set
\begin{equation}
\varphi(z) = \frac{u(z)}{v(z)}
\label{r5.18}
\end{equation}
and suppose that $v(z) \not \equiv 0$ (i.e., $v(z) \neq 0$ for
every $z \in \rho(\w)$).
Substituting (\ref{r5.18}) into (\ref{r5.16}) yields
\begin{equation*}
\frac{v(z)\varphi(z)v(\zeta)^* - v(\zeta)^*\varphi(\zeta)^*v(z)}{
z - \bar{\zeta}} =
v(z) \frac{\varphi(z) - \varphi(\zeta)^*}{z - \bar{\zeta}}v(\zeta)^*
\succeq 0 \qquad (z, \zeta \in \rho(\w)).
\end{equation*}
From (\ref{r5.17})
we conclude that
\begin{equation*}
K_\varphi(z,\zeta):=
\frac{\varphi(z) - \varphi(\zeta)^*}{z - \bar{\zeta}} \succeq 0 \qquad
(z, \zeta \in \rho(\w)).
\end{equation*}
This latter equality means that (after an analytic continuation to
all of $\mathbb{C}^+$) $\varphi$ is a Nevanlinna function.
Finally, it follows from (\ref{r5.15}) that
\begin{equation*}
\Theta \begin{bmatrix}u \\ v\end{bmatrix} =
\begin{bmatrix}\Theta_{11}u + \Theta_{12}v \\
\Theta_{21}u + \Theta_{22}v\end{bmatrix}
\end{equation*}
which implies
\begin{equation*}
\w = \frac{\Theta_{11}u + \Theta_{12}v}{\Theta_{21}u + 
\Theta_{22}v} =
\frac{\Theta_{11}\varphi + \Theta_{12}}{\Theta_{21}\varphi +
\Theta_{22}} = {\bf T}_\Theta[\varphi].
\end{equation*}

Now suppose that $v(z) \equiv 0$ on $\rho(\w)$.  Then
\begin{equation*}
\Theta(z)\begin{bmatrix}u(z) \\ 0\end{bmatrix} =
\begin{bmatrix}\Theta_{11}(z)u(z) \\ \Theta_{21}(z)u(z)\end{bmatrix} =
\begin{bmatrix}\w(z) \\ 1\end{bmatrix},
\end{equation*}
from which it follows that $u(z) \equiv \Theta_{21}(z)^{-1}$, which
in turn implies $\w(z) \equiv \frac{\Theta_{11}(z)}{\Theta_{21}(z)}$.
Such $\w$ indeed satisfies (\ref{r5.14}) since it results in
$\mathbf{S} \equiv 0$, as may be seen by using (\ref{r5.16}).  
Note that this situation corresponds to
$\varphi(z) \equiv \infty$ and so $\varphi \in \tN_0$ as required.

Thus we have shown that if a meromorphic $\w$ satisfies
(\ref{r5.14}), then $\w$ is of the form (\ref{r5.12}).
Let us now assume we have a meromorphic function $\w$ given by
\begin{equation*}
\w(z) = \frac{\Theta_{11}(z)\varphi(z) + \Theta_{12}(z)}{
\Theta_{21}(z)\varphi(z) + \Theta_{22}(z)}
\end{equation*}  
for some Nevanlinna function $\varphi$.
Set
\begin{equation*}
V(z) := \Theta_{21}(z)\varphi(z) + \Theta_{22}(z).
\end{equation*}
If $V(z) \equiv 0$, then $\w(z) \equiv \infty$ 
and so is not a meromorphic function.
Therefore $V(z) \not \equiv 0$.  This permits us to express $\w$ in the following way:
\begin{equation*}
\begin{bmatrix}\w(z) \\ 1\end{bmatrix} = 
\Theta(z)\begin{bmatrix}\varphi(z) \\ 1\end{bmatrix}
\frac{1}{V(z)},
\end{equation*}
which leads us to
\begin{align*}
\frac{\varphi(z) - \varphi(\zeta)^*}{z - \bar{\zeta}}
&= \begin{bmatrix}\varphi(\zeta)^* & 1\end{bmatrix}
\frac{J}{i(z - \bar{\zeta})}
\begin{bmatrix}\varphi(z) \\ 1 \end{bmatrix} \\
&= V(\zeta)^*
\begin{bmatrix}\w(\zeta)^* & 1\end{bmatrix}
\frac{\Theta(\zeta)^{-*}J\Theta(z)^{-1}}{i(z - \bar{\zeta})}
\begin{bmatrix}\w(z) \\ 1\end{bmatrix}V(z)
\end{align*}
Since
$\varphi$ is a Nevanlinna function, the latter kernel is positive
on $\rho(\w) \times \rho(\w)$, whence $\mathbf{S}(z, \zeta) \succeq
0$.
\end{proof}

\chapter{Parameters and Interpolation Conditions I}
In this section we characterize
the boundary behavior of $\w = \mathbf{T}_\Theta[\varphi]$ 
around interpolation nodes $x_i$ for which $\tilde{e}_i = 0$.  In
particular, we are interested in under what circumstances
$\w$ satisfies some or all of the interpolation conditions 
stated in Problem \ref{Prob1} or in Problem \ref{Prob2}, i.e.,
the conditions
(\ref{1.26})--(\ref{1.27c}) or
the conditions (\ref{1.29}) and (\ref{1.30}), respectively.

Let us define the functions
\begin{equation}
U_\varphi(z) = \Theta_{11}(z)\varphi(z) + \Theta_{12}(z)
\quad\text{and}\quad
V_\varphi(z) = \Theta_{21}(z)\varphi(z) + \Theta_{22}(z)
\label{r6.1}
\end{equation}
for a fixed function $\varphi \in \N_0$ so that
\begin{equation}
\begin{bmatrix}U_\varphi(z) \\ V_\varphi(z)\end{bmatrix} =
\Theta(z)\begin{bmatrix}\varphi(z) \\ 1\end{bmatrix}.
\label{r6.2}
\end{equation}
This permits the following representation of $\w$:
\begin{equation}
\w(z) := {\bf T}_\Theta[\varphi](z) = \frac{U_\varphi(z)}{V_\varphi(z)}.
\label{r6.3}
\end{equation}
Set
\begin{align}
\Psi_\varphi(z) &= i(zI-X)^{-1}
\begin{bmatrix}\tC^* & \tE^*\end{bmatrix}J
\begin{bmatrix}\varphi(z) \\ 1\end{bmatrix}
\nonumber \\
&= (zI - X)^{-1}\left(\tC^* - \varphi(z)\tE^*
\right).
\label{r6.4}
\end{align}
Substituting (\ref{r6.4}) and (\ref{ThetaFormula2})
into (\ref{r6.2})
yields
\begin{equation}
U_\varphi(z) = \varphi(z) - C\Psi_\varphi(z) 
\quad\text{and}\quad
V_\varphi(z) = 1 - E\Psi_\varphi(z).
\label{r6.4a}
\end{equation}
Due to equality (\ref{r6.3}), we have
\begin{equation}
K_\w(z,\zeta) =
\frac{\w(z) - \w(\zeta)^*}{z - \bar{\zeta}} =
\frac{1}{V_\varphi(\zeta)^* V_\varphi(z)} 
\frac{V_\varphi(\zeta)^*U_\varphi(z) - V_\varphi(z)U_\varphi(\zeta)^*}
{z - \bar{\zeta}}.
\label{r6.5}
\end{equation}
Note that by (\ref{r6.2}), 
\begin{align*}
V_\varphi(\zeta)^*U_\varphi(z) - V_\varphi(z)U_\varphi(\zeta)^* &=
-i\begin{bmatrix}U_\varphi(\zeta)^* & V_\varphi(\zeta)^*\end{bmatrix}J
\begin{bmatrix}U_\varphi(z) \\ V_\varphi(z)\end{bmatrix} \\
&= -i\begin{bmatrix}\varphi(\zeta)^* & 1\end{bmatrix}
\Theta(\zeta)^*J\Theta(z)
\begin{bmatrix}\varphi(z) \\ 1\end{bmatrix}.
\end{align*}
By (\ref{Junitary}),
the right hand side of the latter equality admits the representation
\begin{equation*}
 -\begin{bmatrix}\varphi(z)^* & 1\end{bmatrix}
\left(-iJ + (z - \bar{\zeta)}J
\begin{bmatrix}\tC \\ \tE\end{bmatrix}
(\bar{\zeta}I - X)^{-1}P(zI - X)^{-1}
\begin{bmatrix}\tC^* & \tE^*\end{bmatrix}J\right)
\begin{bmatrix}\varphi(z) \\ 1\end{bmatrix},
\end{equation*}
which simplifies, yielding,
\begin{equation*}
V_\varphi(\zeta)^*U_\varphi(z) - V_\varphi(z)U_\varphi(\zeta)^*
= \varphi(z) - \varphi(\zeta)^* + (z - \bar{\zeta})\Psi_\varphi(\zeta)^*P\Psi_\varphi(z),
\end{equation*}
where $\Psi_\varphi$ is given in (\ref{r6.4}).
Thus (\ref{r6.5}) may be expressed as follows:
\begin{equation}
K_\w(z,\zeta) =
\frac{1}{V_\varphi(\zeta)^*V_\varphi(z)}
\left( K_\varphi(z,\zeta) +
\Psi_\varphi(\zeta)^*P\Psi_\varphi(z)\right).
\label{r6.6}
\end{equation}
\begin{Rk}
\emph{
Equality (\ref{r6.6}) implies that for every function
$\varphi \in \N_0$ and $\Theta \in \mathcal{W}_\ka$, the function
$\w := {\bf T}_\Theta[\varphi]$ belongs to the generalized Nevanlinna
class $\tN_{\ka^\prime}$ for some $\ka^\prime \leq \ka$.
Indeed, it follows from (\ref{r6.6}) that 
${\rm sq}_-K_\w \leq {\rm sq}_-K_\varphi + {\rm sq}_-P = 0 + \ka$.}
\label{Remark6.1}
\end{Rk}
Upon evaluating equality (\ref{r6.6}) at $\zeta = z$, we obtain
\begin{equation}
K_\w(z,z) = 
\frac{1}{|V_\varphi(z)|^2}\left(
K_\varphi(z,z) +
\Psi_\varphi(z)^*P\Psi_\varphi(z)\right).
\label{r6.7}
\end{equation}

For $\varphi \equiv \infty$, we interpret the above expressions
projectively, which means that we take
\begin{equation}
U_\infty(z) = \Theta_{11}(z)
\quad\text{and}\quad
V_\infty(z) = \Theta_{21}(z)
\label{r6.7a}
\end{equation}
which become
\begin{equation}
U_\infty(z) = 1 - C\Psi_\infty(z)
\quad\text{and}\quad
V_\infty(z) = -E\Psi_\infty(z)
\label{r6.7b}
\end{equation}
upon setting
\begin{equation}
\Psi_\infty(z) = -(zI-X)^{-1}\tE^*.
\label{r6.7c}
\end{equation}
Then
\begin{equation}
K_\w(z,z) = \frac{\Psi_\infty(z)^*P\Psi_\infty(z)}{|V_\infty(z)|^2}.
\label{r6.8}
\end{equation}

\section{Interpolation Conditions at Singular Interpolation Nodes}

Throughout this section we let $i \in \{\ell + 1, \ldots, n\}$
be fixed.  We use the conformal partitionings introduced
in Remark \ref{BlockRemark}. 
We take $k = 1$ and, in particular, we rearrange interpolation nodes
({\it after} first choosing $i \in \{\ell + 1, \ldots, n\}$) so that
$P_{11} = \xi_i$.  The arrangement of the remaining points
may be arbitrary.
Since $\xi_i \neq 0$, we may define the functions $\Theta^{(1)}$ and 
$\tilde{\Theta}^{(2)}$ using (\ref{e3.23}) and (\ref{e3.24}):
\begin{align} 
\Theta^{(1)}(z) &= I_2 - i\begin{bmatrix}\xi_i \\ 0\end{bmatrix}\cdot
\frac{-\xi_i^{-1}}{z - x_i} \cdot \begin{bmatrix}\xi_i & 0\end{bmatrix} J
= \begin{bmatrix}1 & \frac{\xi_i}{z - x_i} \\ 0 & 1\end{bmatrix},
\label{r6.9}
\\
\tilde{\Theta}^{(2)}(z) &:=
\begin{bmatrix}\tilde{\Theta}_{11}^{(2)}(z) & \Tilde{\Theta}_{12}^{(2)}(z) \\
\tilde{\Theta}_{21}^{(2)}(z) & \Tilde{\Theta}_{22}^{(2)}(z)\end{bmatrix} =
I_2 - i\begin{bmatrix}\tC_2 \\
  \tE_2\end{bmatrix}
\tilde{P}_{22}^{-1}(zI - X_2)^{-1}\begin{bmatrix}\tC_2^* &
  \tE_2^*\end{bmatrix} J. \nonumber
\end{align}
Note that by Lemma \ref{InvLemma},
\begin{equation}
\Theta(z) = \Theta^{(1)}(z)\tilde{\Theta}^{(2)}(z).
\label{r6.10}
\end{equation}
\begin{Rk}
Since the poles of $\tilde{\Theta}^{(2)}$ belong to
$\sigma(X_2) = \{x_1, \ldots, x_n\} \setminus \{x_i\}$,
the function $\tilde{\Theta}^{(2)}$ is analytic at $x_i$.
\label{PoleRemark}
\end{Rk}
Substituting (\ref{r6.9}) into (\ref{r6.10}) yields
\begin{equation}
\Theta(z) = \begin{bmatrix}\tilde{\Theta}_{11}^{(2)}(z) + \frac{\xi_i}{z -
    x_i}\tilde{\Theta}_{21}^{(2)}(z) &
\tilde{\Theta}_{12}^{(2)}(z) + \frac{\xi_i}{z -
    x_i}\tilde{\Theta}_{22}^{(2)}(z) \\
\tilde{\Theta}_{21}^{(2)}(z) &
\tilde{\Theta}_{22}^{(2)}(z)\end{bmatrix}.
\label{r6.11}
\end{equation}
Comparing entries in (\ref{r6.11}) with those in
\begin{equation*}
\Theta(z) :=
\begin{bmatrix}
\Theta_{11}(z) & \Theta_{12}(z)\\
\Theta_{21}(z) & \Theta_{22}(z)
\end{bmatrix},
\end{equation*}
shows that
\begin{align}
\tilde{\Theta}_{11}^{(2)}(z) + \frac{\xi_i}{z -
    x_i}\tilde{\Theta}_{21}^{(2)}(z) &\equiv \Theta_{11}(z),
&\tilde{\Theta}_{21}^{(2)}(z) \equiv \Theta_{21}(z),
\label{r6.12}
\\
\tilde{\Theta}_{12}^{(2)}(z) + \frac{\xi_i}{z -
    x_i}\tilde{\Theta}_{22}^{(2)}(z)
&\equiv \Theta_{12}(z),
&\tilde{\Theta}_{22}^{(2)}(z) \equiv \Theta_{22}(z).
\label{r6.12a}
\end{align}
By combining relations (\ref{r6.12}), (\ref{r6.12a}), and (\ref{r6.11}), we may express
in the following way
the functions $U_\varphi$, defined in (\ref{r6.1}):
\begin{equation}
U_\varphi(z) = \left( \tilde{\Theta}_{11}^{(2)}(z) + \frac{\xi_i}{z -
    x_i}\Theta_{21}(z) \right)\varphi(z) +
\tilde{\Theta}_{12}^{(2)}(z) + \frac{\xi_i}{z -
    x_i}\Theta_{22}(z).
\label{r6.13}
\end{equation}
Multiplying the Lyaponov equation (\ref{e2.14})
by $\mathbf{e}_i^*$ on the left and using the fact that
$e_i = 0$ shows that
\begin{equation*}
\mathbf{e}_i^*P(x_i I - X) = \xi_iE
\end{equation*}
and hence that
\begin{equation}
E = \frac{1}{\xi_i}\mathbf{e}_i^*P(x_i I - X).
\label{r6.14}
\end{equation}
Note that
\begin{align}
\lim_{z \to x_i} (x_i I - X)(zI - X)^{-1} &=
I - \mathbf{e}_i\mathbf{e}_i^*,
\label{r72lim1}\\
\lim_{z \to x_i} (z - x_i)(zI - X)^{-1} &=
\mathbf{e}_i\mathbf{e}_i^*
\label{r72lim2}
\end{align}
by the definition (\ref{r2.2})
of $X$.
By (\ref{ThetaFormula2}) and (\ref{r6.14}), we get
\begin{align}
\Theta_{21}(x_i) &= \lim_{z \to x_i}\Theta_{21}(z) 
= \lim_{z \to x_i}\frac{1}{\xi_i}\mathbf{e}_i^*P(x_i I -
X)(zI - X)^{-1}P^{-1}E^* \nonumber \\
&= \frac{1}{\xi_i}\mathbf{e}_i^*P(I - \mathbf{e}_i
\mathbf{e}_i^*)P^{-1}E^*
= \tilde{e}_i.
\label{r6.15}
\end{align}
where we have made use of the limit (\ref{6LIM})
and the fact that $\mathbf{e}_i^*P\mathbf{e}_i = -\xi_i$.
Similarly,
\begin{equation}
\Theta_{22}(x_i) = \lim_{z \to x_i}\left(1 - \frac{1}{\xi_i}
\mathbf{e}_i^*P(x_i I - X)(zI - X)^{-1}P^{-1}C^*\right)
= -\tilde{c}_i.
\label{r6.16}
\end{equation}

Now we prove Theorem \ref{Main2}.  We repeat the formulation,
but in an abbreviated form.
\begin{Tm}
Let $i \in \{\ell + 1, \ldots, n\}$ and suppose that
$\tilde{e}_i = 0$ (i.e., $\eta_i = \infty$).
Let
$\varphi \in \tN_0$ and let $\w:= {\bf T}_\Theta[\varphi]$.  Then
\begin{align*}
&\text{$\varphi$ satisfies $\btC_{1-2}$} &\iff &\quad 
&&\w_{-1}(x_i) = \xi_i,\\
&\text{$\varphi$ satisfies $\btC_3$} &\iff &\quad
&&-\infty < \frac{-1}{\w_{-1}(x_i)} < \frac{-1}{\xi_i},\\
&\text{$\varphi$ satisfies $\btC_4$} &\iff &\quad
&&\frac{-1}{\xi_i} < \frac{-1}{\w_{-1}(x_i)} < \infty,\\
&\text{$\varphi$ satisfies $\btC_{5-6}$} &\iff &\quad
&&\w_{-1}(x_i) = 0.
\end{align*}
\end{Tm}
\begin{proof}
Our objective is to prove the existence of the limit
\begin{equation}
\lim_{z \to x_i} (z-x_i)\frac{U_\varphi(z)}{V_\varphi(z)} =
\lim_{z \to x_i} (z-x_i)\w(z) =: \w_{-1}(x_i)
\label{r6.17}
\end{equation}
and express its value in terms of the data and the parameter $\varphi$.
Existence will follow from the evaluations of the limits
\begin{equation*}
\lim_{z \to x_i} (z-x_i)U_\varphi(z) 
\quad\text{and}\quad
\lim_{z \to x_i} V_\varphi(z).
\end{equation*}

By multiplying the identity (\ref{e3.4}) by $\mathbf{e}_i$ on the right, 
and by using the fact that $\tilde{e}_i = 0$, we obtain
\begin{equation}
\tE^* = \frac{-1}{\tilde{c}_i}(x_iI-X)P^{-1}\mathbf{e}_i.
\label{r6.18}
\end{equation}
By (\ref{r6.15}),
\begin{equation}
\Theta_{21}(x_i) = \tilde{e}_i = 0.
\label{r6.19}
\end{equation}
Using (\ref{ThetaFormula}), we may
calcuate the derivatives with respect to
$z$ of the functions $\Theta_{11}$ 
and $\Theta_{21}$:
\begin{equation}
\Theta_{11}^\prime(z) = -C(zI-X)^{-2}\tE^*
\quad\text{and}\quad
\Theta_{21}^\prime(z) = -E(zI-X)^{-2}\tE^*.
\label{r6.20}
\end{equation}
Substituting (\ref{r6.14}) and (\ref{r6.18}) into
(\ref{r6.20}), we have
\begin{equation*}
\Theta_{21}^\prime(z) =
\frac{1}{\xi_i \tilde{c}_i}\mathbf{e}_i^*P(x_iI-X)(zI-X)^{-2}
(x_iI-X)P^{-1}\mathbf{e}_i.
\end{equation*}
Taking limits as $z \to x_i$ yields
\begin{equation}
\Theta_{21}^\prime(x_i) = \frac{1}{\xi_i \tilde{c}_i}
\mathbf{e}_i^*P(I - \mathbf{e}_i\mathbf{e}_i^*)P^{-1}\mathbf{e}_i
= \frac{1}{\xi_i \tilde{c}_i} + \frac{\tilde{p}_{ii}}{\tilde{c}_i},
\label{r6.21}
\end{equation}
where the last equality is due to the fact that
\begin{equation*}
\mathbf{e}_i^*P\mathbf{e}_i = -\xi_i
\quad\text{and}\quad
\mathbf{e}_i^*P\mathbf{e}_i = \tilde{p}_{ii}.
\end{equation*}
By (\ref{r6.12}), the function $\Theta_{21}(z)$ is identified with the function
$\tilde{\Theta}_{21}^{(2)}(z)$, which
is analytic at $x_i$ by Remark \ref{PoleRemark}. 
Therefore,
by (\ref{r6.19}) and (\ref{r6.20}), $\Theta_{21}(z)$ 
admits the Taylor series expansion
\begin{equation*}
\Theta_{21}(z) = (z-x_i)\Theta_{21}^\prime(x_i) + O(|z - x_i|^2).
\end{equation*}
Hence
\begin{equation*}
\lim_{z \to x_i}\Theta_{21}(z)\varphi(z) =
\lim_{z \to x_i}\left((z - x_i)\varphi(z)\Theta_{21}^\prime(x_i) +
O(|z - x_i|^2)\varphi(z)\right).
\end{equation*}
The last term on the right hand side tends to zero
by Lemma \ref{CJN2}, and thus,
by (\ref{r6.21}),
\begin{equation}
\lim_{z \to x_i}\Theta_{21}(z)\varphi(z) =
\varphi_{-1}(x_i)\left(\frac{1}{\xi_i \tilde{c}_i} + \frac{\tilde{p}_{ii}}{\tilde{c}_i}
\right).
\label{r6.22}
\end{equation}
It now follows from (\ref{r6.16}) and (\ref{r6.22}) that
\begin{equation}
\lim_{z \to x_i}V_\varphi(z) = 
\varphi_{-1}(x_i) \left( 
 \frac{1}{\xi_i \tilde{c}_i} + \frac{\tilde{p}_{ii}}{\tilde{c}_i}
\right) - \tilde{c}_i.
\label{r6.23}
\end{equation}
Furthermore, in view of (\ref{r6.13}),
\begin{align}
\lim_{z \to x_i} (z - x_i)U_\varphi(z) =&
\lim_{z \to x_i} \left(\xi_i (\Theta_{21}(z)\varphi(z) +
  \Theta_{22}(z)) \right. \nonumber \\
&\left. + \tilde{\Theta}_{11}^{(2)}(z)(z - x_i)\varphi(z)
+ (z - x_i)\tilde{\Theta}_{12}^{(2)}(z)\right).
\label{r6.24}
\end{align}
The functions $\tilde{\Theta}_{11}^{(2)}$ and $\tilde{\Theta}_{12}^{(2)}$ are analytic
at $x_i$ and so in particular
$(z - x_i)\tilde{\Theta}_{12}^{(2)}(z)$ tends to zero and
$\tilde{\Theta}_{11}^{(2)}(z)$ has a finite limit at $x_i$.  Hence
$\tilde{\Theta}_{11}^{(2)}(z)(z - x_i)\varphi(z)$ tends to
$\tilde{\Theta}_{11}^{(2)}(x_i)\varphi_{-1}(x_i)$.
The term $\xi_i \left(\Theta_{21}(z)\varphi(z) + \Theta_{22}(z)\right)$
is just $\xi_i V(z)$.
Therefore, to find the limit in (\ref{r6.24}), we must determine
the limit $\tilde{\Theta}_{11}^{(2)}(x_i)$.
By (\ref{r6.12a}) and (\ref{ThetaFormula2}),
\begin{equation*}
\tilde{\Theta}_{11}^{(2)}(z) + \frac{\xi_i}{z -
  x_i}E(zI-X)^{-1}\tE^*
= \Theta_{11}(z) =
1 + C(zI - X)^{-1}\tE^*,
\end{equation*}
from which it follows that
\begin{equation}
\tilde{\Theta}_{11}^{(2)}(z) = 1 + \left(C - \frac{\xi_i}{z - x_i}
E\right)(zI-X)^{-1}\tE^*.
\label{r6.25}
\end{equation}
Due to (\ref{r6.14}) and (\ref{r6.18}), taking limits in (\ref{r6.25}) 
as $z \to x_i$ yields
\begin{align}
\tilde{\Theta}_{11}^{(2)}(x_i) &= 1 - \frac{1}{\tilde{c}_i}C
(I - \mathbf{e}_i\mathbf{e}_i^*)P^{-1}\mathbf{e}_i -
\frac{1 + \xi_i \tilde{p}_{ii}}{\tilde{c}_i} \nonumber \\
&= 1 - \frac{\tilde{c}_i}{\tilde{c}_i} + \frac{\xi_i \tilde{p}_{ii}}{
\tilde{c}_i} - \frac{1 + \xi_i\tilde{p}_{ii}}{\tilde{c}_i} \nonumber \\
&= -\frac{1}{\tilde{c}_i}.
\label{r6.26}
\end{align}
Therefore
\begin{align}
\lim_{z \to x_i} (z - x_i)U(z) &=
\xi_i V(x_i) - \varphi_{-1}(x_i) \cdot \frac{1}{\tilde{c}_i} \nonumber
\\
&= \xi_i \left(\varphi_{-1}(x_i) \cdot \frac{1 + \xi_i\tilde{p}_{ii}}{
\xi_i \tilde{c}_i} - \tilde{c}_i\right) - \varphi_{-1}(x_i) \cdot
\frac{1}{\tilde{c}_i} \nonumber \\
&= \xi_i \left(\varphi_{-1}(x_i)\frac{\tilde{p}_{ii}}{\tilde{c}_i} - 
\tilde{c}_i\right),
 \label{r6.27}
\end{align}
where the second equality follows from (\ref{r6.23}).
Thus, by (\ref{r6.27}), (\ref{r6.23}), and (\ref{r6.17}),
\begin{equation}
\w_{-1}(x_i) =
\frac{\xi_i \left(\varphi_{-1}(x_i) \frac{\tilde{p}_{ii}}{\tilde{c}_i}
- \tilde{c}_i\right)}{
\left(\varphi_{-1}(x_i) \frac{\tilde{p}_{ii}}{\tilde{c}_i}
- \tilde{c}_i\right) + \varphi_{-1}(x_i) \frac{1}{\xi_i \tilde{c}_i}}.
\label{r6.28}
\end{equation}
If $\varphi$ satisfies condition $\btC_{1-2}$ at $x_i$, then
$\varphi_{-1}(x_i) = 0$.  Substituting this
into (\ref{r6.28}) shows that $\w_{-1}(x_i) = \xi_i$.
If $\varphi$ meets condition $\btC_{3-5}$ at $x_i$, then
\begin{equation}
-\infty < \varphi_{-1}(x_i) < 0.
\label{6.28a}
\end{equation}
When
$\btC_{5}$ is satisfied, we have
$\varphi_{-1}(x_i)\tilde{p}_{ii} - \tilde{c}_i^2 = 0$, and 
substituting this equality 
into (\ref{r6.28}) shows that $\w_i = 0$.  When
$\btC_{3-4}$ holds, $\varphi_{-1}(x_i)\tilde{p}_{ii} -
\tilde{c}_i^2 \neq 0$, and so $\w_{-1}^{-1}(x_i)$ exists
and is finite, by (\ref{r6.28}).  A simple calculation then
shows that
\begin{equation}
\frac{1}{\w_{-1}(x_i)} - \frac{1}{\xi_i} =
\frac{\varphi_{-1}(x_i)}{\xi_i^2(\tilde{p}_{ii}
\varphi_{-1}(x_i) - \tilde{c}_i^2)}.
\label{r6.29}
\end{equation}
Formula (\ref{r6.29}) combined with (\ref{6.28a})
implies that if $\varphi$ satisfies $\btC_{3}$ at $x_i$,
then $\w_{-1}^{-1}(x_i) < \xi_i^{-1}$, and if
$\varphi$ satisfies $\btC_4$ at $x_i$, then
$\w_{-1}^{-1}(x_i) > \xi_i^{-1}$.

Now suppose that $\varphi$ meets $\btC_6$ at $x_i$.
In order for this to occur, we must have $\varphi \equiv \infty$
and $\tilde{p}_{ii} = 0$.  The formulas for $U_\infty$ and $V_\infty$
then become, on acount of (\ref{r6.7a}), (\ref{r6.12}), and (\ref{r6.12a}),
\begin{equation}
U_\infty(z) = \tilde{\Theta}_{11}^{(2)}(z) + \frac{\xi_i}{z -
  x_i}\Theta_{21}(z) 
\quad\text{and}\quad
V_\infty(z) = \Theta_{21}(z),
\label{r6.30}
\end{equation}
and the Taylor series expansion of $\Theta_{21}(z)$ 
at $x_i$ simplifies to
\begin{equation}
\Theta_{21}(z) = \frac{z-x_i}{\xi_i \tilde{c}_i} + O(|z - x_i|^2).
\label{r6.31}
\end{equation}
Together (\ref{r6.30}) and (\ref{r6.31}), along
with (\ref{r6.26}), show that
\begin{align*}
\lim_{z \to x_i}U_\infty(z) &=
\lim_{z \to x_i} \left(\tilde{\Theta}_{11}^{(2)}(z) +
\frac{1}{\tilde{c}_i} + O(|z - x_i|)\right)
= -\frac{1}{\tilde{c}_i} + \frac{1}{\tilde{c}_i} = 0,\\
\lim_{z \to x_i} (z - x_i)^{-1}V_\infty(z) &=
\lim_{z \to x_i} \left(
\frac{1}{\xi_i \tilde{c}_i} + O(|z - x_i|)\right)
= \frac{1}{\xi_i \tilde{c}_i}.
\end{align*}
Therefore
\begin{equation*}
\w_{-1}(x_i) = \lim_{z \to x_i}(z -
x_i)\frac{U_\infty(z)}{V_\infty(z)} =
\frac{0}{(\xi_i \tilde{c}_i)^{-1}} = 0.
\end{equation*}

Since every function $\varphi \in \tN$ satisfies exactly one of
$\btC_1, \ldots, \btC_6$, the proof
is complete.
\end{proof}

We may now solve Problem \ref{Prob1} for invertible $P$ when 
$\ka := {\rm sq}_-P = 0$ and
the data $\Omega$
is such that $\ell = 0$ (i.e., all interpolation nodes
are singular interpolation nodes).
\begin{Tm}
Let $\Omega$ be such that $\ell = 0$:
\begin{equation*}
\Omega = \left\{ \{x_i\}_{i=1}^n, \{\xi_i\}_{i=1}^n \right\}.
\end{equation*}
Let $P$ be the
Pick matrix associated to $\Omega$ and suppose that $P$ is
positive definite.  Then all solutions of Problem \ref{Prob1}
are parameterized by $\mathbf{T}_\Theta[\varphi]$, where
$\varphi \in \tN_0$ and satisfies $\btC_{1-2}$
at each interpolation node.
\label{ellis0}
\end{Tm}
\begin{proof}
Since $\ell = 0$, the matrix $E$ is the zero vector, and hence
$\tilde{E} := EP^{-1} = 0$ also.  Since each $\tilde{e}_i = 0$,
we need only examine under what conditions parameters $\varphi$
satisfy one of conditions $\btC_1$ through $\btC_6$
(since conditions $\bC_1$ through $\bC_6$ only apply when
$\tilde{e}_i \neq 0$).  
Note that positive definiteness of $P$ (and
hence that of $P^{-1}$) precludes $\varphi$ from satisfying
$\btC_{4-6}$.  
Observe that Theorem \ref{Main2} has been proven,
the relevant parts of Theorem \ref{Main1} dealing with conditions
$\btC_1$ through $\btC_3$ have already been proven, and the necessity
part of Theorem \ref{lessprecise} has been proven.  Note also
that, by Remark \ref{Remark6.1}, 
a function $\w$ of the form $\mathbf{T}_\Theta[\varphi]$
for $P > 0$ belongs to the class $\N_0$. 
Therefore, due to 
Theorems \ref{lessprecise}, \ref{Main1}, and \ref{Main2}, a function
$\w = \mathbf{T}_\Theta[\varphi]$ is a solution of Problem
\ref{Prob1} (under assumption that $P > 0$) if and only if
$\varphi \in \tN_0$ satisfies $\btC_{1-2}$ at each interpolation
node.
\end{proof}

\section{Interpolation Conditions at Regular Interpolation Nodes}
In this section we prove Theorem \ref{Main1} by breaking it
into three separate Theorems.  In each case, the
relevant parts of Theorem \ref{Main1} will be repeated
in an abbreviated form.
\begin{Tm}
Let $i \in \{1, \ldots, \ell\}$ and suppose that $\tilde{e}_i = 0$
(i.e., that $\eta_i = \infty$).  Suppose also
that $\varphi \in \tN_0$ satisfies
condition $\btC_{1-4}$ at $x_i$.
Let $\w := \mathbf{T}_\Theta[\varphi]$.
Then
\begin{enumerate}
\item  If $\varphi$ satisfies $\btC_{1-2}$, then
$\w(x_i) = \w_i$ and $\w^\prime(x_i) = \gamma_i$.
\item If $\varphi$ satisfies $\btC_{3}$, then 
$\w(x_i) = \w_i$ and $-\infty < \w^\prime(x_i) < \gamma_i$.
\item If $\varphi$ satisfies $\btC_{4}$, then
$\w(x_i) = \w_i$ and $\gamma_i < \w^\prime(x_i) < \infty$.
\end{enumerate}
\label{T6.4}
\end{Tm}
\begin{proof}
By (\ref{ThetaFormula2}),
\begin{equation}
\Theta_{11}(z) = 1 + C(zI - X)^{-1}\tE^*
\quad\text{and}\quad
\Theta_{21}(z) = E(zI - X)^{-1}\tE^*.
\label{s6.1}
\end{equation}
Since $\tilde{e}_i = 0$, we may use the formula (\ref{r6.18})
for $E$ in the latter equalities.  We then get
\begin{align}
\Theta_{11}(z) &= 
1 - \frac{1}{\tilde{c}_i}C(zI - X)^{-1}(x_iI - X)P^{-1}\mathbf{e}_i,
\label{LimT11} \\
\Theta_{21}(z) &= 
-\frac{1}{\tilde{c}_i}E(zI - X)^{-1}(x_iI-X)P^{-1}\mathbf{e}_i.
\label{LimT21}
\end{align}
Taking limits in (\ref{LimT11}) yields
\begin{align}
\lim_{z \to x_i}\Theta_{11}(z) &= 1 - 
\frac{1}{\tilde{c}_i}C(I -
\mathbf{e}_i\mathbf{e}_i^*)P^{-1}\mathbf{e}_i
= 1 - \frac{1}{\tilde{c}_i}\left(CP^{-1}\mathbf{e}_i -
  \tilde{p}_{ii}C\mathbf{e}_i\right) \nonumber \\
&= \w_i \cdot \frac{\tilde{p}_{ii}}{\tilde{c}_i}.
\label{r6.35}
\end{align}
Similarly, taking limits in (\ref{LimT21}) gives
\begin{equation}
\lim_{z \to x_i}\Theta_{21}(z) = \frac{\tilde{p}_{ii}}{\tilde{c}_i}.
\label{r6.36}
\end{equation}
Due to Lemma \ref{L42}, relations (\ref{r6.35}) and (\ref{r6.36}) respectively
imply that
\begin{align}
\lim_{z \to x_i}\Theta_{11}(z)\varphi(z)(z - x_i) &= \w_i \cdot
\frac{\tilde{p}_{ii}}{\tilde{c}_i} \varphi_{-1}(x_i) \label{r6.37}\\
\lim_{z \to x_i}\Theta_{21}(z)\varphi(z)(z - x_i) &=
\frac{\tilde{p}_{ii}}{\tilde{c}_i} \varphi_{-1}(x_i). \label{r6.38}
\end{align}
By Lemma \ref{L42},
\begin{equation}
\lim_{z \to x_i}(z - x_i)\Theta_{12}(z) = -\w_i \tilde{c}_i
\quad\text{and}\quad
\lim_{z \to x_i}(z - x_i)\Theta_{22}(z) = -\tilde{c}_i.
\label{r6.39}
\end{equation}
Since $\varphi$ satisfies $\btC_{1-4}$ at $x_i$, it holds that
\begin{equation}
\tilde{p}_{ii} \varphi_{-1}(x_i) - \tilde{c}_i^2 \neq 0.
\label{r6.40}
\end{equation}
Together (\ref{r6.37}) -- (\ref{r6.40}) show that
\begin{align*}
\lim_{z \to x_i}\w(z) &= \lim_{z \to x_i}
\frac{\Theta_{11}(z)\varphi(z)(z - x_i) + (z - x_i)\Theta_{12}(z)}{
\Theta_{21}(z)\varphi(z)(z-x_i) + (z-x_i)\Theta_{22}(z)} \nonumber \\
&= \frac{\w_i \cdot \frac{\tilde{p}_{ii}}{\tilde{c}_i}\varphi_{-1}(x_i) -
  \w_i\tilde{c}_i}{\frac{\tilde{p}_{ii}}{\tilde{c}_i}\varphi_{-1}(x_i) - \tilde{c}_i}
 \nonumber \\
&= \frac{\w_i (\tilde{p}_{ii}\varphi_{-1}(x_i) - \tilde{c}_i^2)}{\tilde{p}_{ii}\varphi_{-1}(x_i) - \tilde{c}_i^2}
= \w_i.
\end{align*}
This proves that $\w(x_i) = \w_i$.
Note that if $\varphi \equiv \infty$ satisfies $\btC_{1-4}$,
then $\tilde{p}_{ii} \neq 0$, and so the desired conclusion
still follows with only minor modifications:
\begin{equation*}
\lim_{z \to x_i} \w(z) = \lim_{z \to x_i}
\frac{\Theta_{11}(z)}{\Theta_{21}(z)} =
\frac{\w_i (\tilde{p}_{ii} / \tilde{c}_i)}{
\tilde{p}_{ii} / \tilde{c}_i} = \w_i.
\end{equation*}
Let us now
compute $K_\w(z,z)$ using (\ref{r6.7}) under the assumption that
$\varphi \not \equiv \infty$.
By (\ref{r6.4}) and (\ref{r6.18}),
\begin{align*}
\lim_{z \to x_i}(z - x_i)\Psi_\varphi(z) &=
\lim_{z \to x_i}(z - x_i)(zI - X)^{-1}\left(\tC^* -
  \varphi(z)\tE^*\right) \\
&= \mathbf{e}_i\mathbf{e}_i^* + \lim_{z \to x_i}
\frac{1}{\tilde{c}_i}(z - x_i)\varphi(z)
(zI - X)^{-1}(x_i I - X) P^{-1} \mathbf{e}_i \\
&= \mathbf{e}_i\mathbf{e}_i^* + \frac{1}{\tilde{c}_i}
\varphi_{-1}(x_i)(I - \mathbf{e}_i\mathbf{e}_i^*)
P^{-1}\mathbf{e}_i \\
&=  \mathbf{e}_i\mathbf{e}_i^* + \frac{1}{\tilde{c}_i}
\varphi_{-1}(x_i)\left(P^{-1} - 
\tilde{p}_{ii}I\right)\mathbf{e}_i.
\end{align*}
Therefore
\begin{align}
\lim_{z \to x_i}|z-x_i|^2\Psi_\varphi(z)^*P\Psi_\varphi(z) =& 
\; \tilde{c}_i^2\gamma_i -
\tilde{c}_i\varphi_{-1}(x_i)\mathbf{e}_i^*P\cdot
\frac{-1}{\tilde{c}_i}(P^{-1}\mathbf{e}_i -
\tilde{p}_{ii}\mathbf{e}_i)  \nonumber \\
& -\tilde{c}_i\varphi_{-1}(x_i)\cdot
\frac{-1}{\tilde{c}_i}(\mathbf{e}_i^* P^{-1} -
\tilde{p}_{ii}\mathbf{e}_i^*)P\mathbf{e}_i \nonumber \\
&
+\varphi_{-1}(x_i)^2\cdot\frac{1}{\tilde{c}_i^2}\cdot(\mathbf{e}_i^*P^{-1}
- \tilde{p}_{ii}\mathbf{e}_i^*)P(P^{-1}\mathbf{e}_i -
\tilde{p}_{ii}\mathbf{e}_i) \nonumber \\
=& \;\tilde{c}_i^2\gamma_i + \varphi_{-1}(x_i)(1 -
\gamma_i\tilde{p}_{ii}) + \varphi_{-1}(x_i)(1 -
\gamma_i\tilde{p}_{ii})  \nonumber \\
& +\frac{\varphi_{-1}(x_i)^2}{\tilde{c}_i^2}(\tilde{p}_{ii}
-2\tilde{p}_{ii} + \tilde{p}_{ii}^2\gamma_i),
\label{r6.41}
\end{align}
which, after further simplification, becomes
\begin{equation}
\lim_{z \to x_i}|z - x_i|^2\Psi_\varphi(z)^*P\Psi_\varphi(z) =
\gamma_i\left(\varphi_{-1}(x_i)\frac{\tilde{p}_{ii}}{\tilde{c}_i} -
  \tilde{c}_i\right)^2 + \varphi_{-1}(x_i)\left(2 - \varphi_{-1}(x_i)
\frac{\tilde{p}_{ii}}{\tilde{c}_i^2}\right).
\label{r6.42}
\end{equation}
By Lemma \ref{CJN2},
\begin{equation}
\lim_{z \to x_i} |z - x_i|^2K_\varphi(z,z) = -\varphi_{-1}(x_i).
\label{Kwphi}
\end{equation}
By (\ref{r6.1}), (\ref{r6.38}), and (\ref{r6.39}), 
\begin{equation}
\lim_{z \to x_i} (z - x_i)V(z) =
\frac{\tilde{p}_{ii}}{\tilde{c}_i} \varphi_{-1}(x_i)
- \tilde{c}_i.
\label{vlim6}
\end{equation}
Therefore, by (\ref{r6.7}), (\ref{r6.41}), (\ref{Kwphi}), and
(\ref{vlim6}), it follows that
\begin{align}
\lim_{z \to x_i}K_\w(z,z) &= \lim_{z \to x_i} \frac{|z -
      x_i|^2\left(K_\varphi(z,z) + \Psi(z)^*P\Psi(z)\right)}{
|z - x_i|^2|V_\varphi(z)|^2} \nonumber \\
&=
      \frac{- \varphi_{-1}(x_i) + \gamma_i\left(\varphi_{-1}(x_i)\frac{\tilde{p}_{ii}}{\tilde{c}_i}
      - \tilde{c}_i\right)^2 + \varphi_{-1}(x_i)\left(2 -
      \varphi_{-1}(x_i)\frac{\tilde{p}_{ii}}{\tilde{c}}^2\right)}{\left(\varphi_{-1}(x_i)\frac{\tilde{p}_{ii}}{\tilde{c}_i}
      - \tilde{c}_i\right)^2} \nonumber \\
&= \gamma_i+ \varphi_{-1}(x_i)\frac{\tilde{c}_i^2 -
      \varphi_{-1}(x_i)\tilde{p}_{ii}}{(\varphi_{-1}(x_i)\tilde{p}_{ii}
      - \tilde{c}_i^2)^2} \nonumber \\
&= \gamma_i -\varphi_{-1}(x_i)\left(\tilde{c}_i^2 - \varphi_{-1}(x_i)\tilde{p}_{ii}\right).
\label{6670}
\end{align}
By Theorem \ref{CJGN1}, $w^\prime(x_i) = K_\w(x_i,x_i)$.  Thus, 
formula (\ref{6670}) together with the fact that
$-\infty < \varphi_{-1}(x_i) \leq 0$ completes the
proof for functions $\varphi \in \N_0$: if
$\varphi$ satisfies $\btC_{1-2}$, then $\varphi_{-1}(x_i) = 0$
and so $\w^\prime(x_i) = \gamma_i$; if $\varphi$ satisfies
$\btC_3$, then $\tilde{c}_i^2 - \varphi_{-1}(x_i)\tilde{p}_{ii} < 0$
and so $\w^\prime(x_i) < \gamma_i$; if $\varphi$ satisfies
$\btC_4$, then $w^\prime(x_i) > \gamma_i$.

It remains to consider the function $\varphi \equiv \infty$.
Note that $\varphi \equiv \infty$ can never satisfy 
$\btC_{1-2}$, and, if it satisfies
$\btC_{3-4}$, then $\tilde{p}_{ii} \neq 0$.
To calculate the limit of $K_\w(z,z)$ under the
assumption that $\varphi \equiv \infty$
satisfies $\btC_{3-4}$, we use formula (\ref{r6.8}):
\begin{align*}
\lim_{z \to x_i}K_\w(z,z) &= \lim_{z \to x_i}\frac{\tE
(\bar{z}I-X)^{-1}P(zI-X)^{-1}\tE^*}{\tE(\bar{z}I-X)^{-1}E^*
E(zI-X)^{-1}\tE^*} \\
&= \frac{\frac{1}{\tilde{c}_i^2}\left(\mathbf{e}_i^*P^{-1} - \tilde{p}_{ii}
\mathbf{e}_i^*\right)P\left(P\mathbf{e}_i -
\tilde{p}_{ii}\mathbf{e}_i\right)}{\frac{1}{\tilde{c}_i^2}\left(\mathbf{e}_i^*
P^{-1} -
\tilde{p}_{ii}\mathbf{e}_i^*\right)E^*E\left(P^{-1}\mathbf{e}_i -
\tilde{p}_{ij}\mathbf{e}_i\right)} \\
&= \frac{\tilde{p}_{ii}^2\gamma_i - \tilde{p}_{ii}}{\tilde{p}_{ii}^2}
= \gamma_i - \frac{1}{\tilde{p}_{ii}}.
\end{align*}
By Theorem \ref{CJGN1}, $\w^\prime(x_i) = K_\w(x_i,x_i)$.  If
$\varphi$ meets condition $\btC_3$, then
$\tilde{p}_{ii} > 0$ and so $\w^\prime(x_i) < \gamma_i$.  On
the other had, if $\varphi$ meets condition $\btC_4$, then
$\tilde{p}_{iI} < 0$ and so $\w^\prime(x_i) > \gamma_i$.
\end{proof}

\begin{Tm}
Let $i \in \{1, \ldots, \ell\}$ and suppose that
$\tilde{e}_i = 0$.  Suppose that $\varphi \in \tN_0$
satisfies condition $\tilde{\bC}_5$ at $x_i$.  Then
the function $\w := \mathbf{T}_\Theta[\varphi]$ is subject to
one of the following:
\begin{enumerate}
\item \text{The nontangential limit $\w(x_i)$ does not exist.}
\item \text{The nontangential limit $\w(x_i)$ exists and $\w(x_i) \neq
  \w_i$.}
\item \text{The nontangential limit $\w(x_i)$ exists and equals
  $\w_i$, and ${\displaystyle \left|\lim_{z \to x_i}K_\w(z,z)\right| = \infty.}$}
\end{enumerate}
\label{T6.5}
\end{Tm}
\begin{proof}
The problem of finding
all functions in $\varphi \in \tN_0$ that
meet condition
$\btC_5$ at $x_i$ is an
interpolation problem of the same type as Problem
\ref{Prob1} with one interpolation node, $\ell = 0$,
and with positive definite Pick matrix,
which is the the scalar $P = -\tilde{c}_i^2 /
\tilde{p}_{ii} > 0$.  
By Theorem \ref{ellis0}, 
any solution $\varphi$ of this problem admits a representation
\begin{equation*}
\varphi = \mathbf{T}_{\hat{\Theta}}[\hat{\varphi}]
\end{equation*}
where $\hat{\Theta}$ is given by
\begin{equation}
\hat{\Theta}(z) = I_2 + \frac{i}{z - x_i}
\begin{bmatrix}\frac{\tilde{c}_i^2}{\tilde{p}_{ii}} \\ 0\end{bmatrix}
\cdot \frac{\tilde{p}_{ii}}{\tilde{c}_i^2}\cdot
\begin{bmatrix}\frac{\tilde{c}_i^2}{\tilde{p}_{ii}} & 0\end{bmatrix}
J
= \begin{bmatrix}1 & \frac{1}{z - x_i}\cdot \frac{\tilde{c}_i^2}{
\tilde{p}_{ii}} \\ 0 & 1\end{bmatrix},
\label{hatTheta}
\end{equation}
and where
$\hat{\varphi}$ satisfies condition $\btC_{1-2}$ at
$x_i$ (note that this excludes $\hat{\varphi} \equiv \infty$).
From (\ref{hatTheta}) it follows that
\begin{equation}
\hat{\Theta}(z)^{-1} = \begin{bmatrix}1 &
\frac{-1}{z - x_i}\cdot \frac{\tilde{c}_i^2}{\tilde{p}_{ii}} \\
0 & 1\end{bmatrix} 
\label{e4.120}
\end{equation}
and
\begin{equation}
{\rm Res}_{z = x_i}\hat{\Theta}(z)^{-1} =
\begin{bmatrix}
0 & -\frac{\tilde{c}_i^2}{\tilde{p}_{ii}} \\
0 & 0
\end{bmatrix}.
\label{6RES}
\end{equation}
Observe that $\hat{\Theta}(z)^{-1}$ coincides with the function 
$\widetilde{\Theta}^{(2)}$ introduced
in (\ref{e3.24}), if we rearrange interpolation nodes
so that $x_i \mapsto x_n$ (this may be done
without loss of generality) and if we take
$E_2 = \tilde{c}_i^2/\tilde{p}_{ii}$, $C_2 = 0$, and 
$\widetilde{P}_{22}^{-1} = -\tilde{p}_{ii}/\tilde{c}_i^2$.  
Therefore, by (\ref{e4.120}) and (\ref{e3.25}),
\begin{equation}
\Theta(z) = \Theta^{(1)}(z)\hat{\Theta}(z)^{-1},
\label{tt1tzi}
\end{equation}
where $\Theta^{(1)}$ is given in (\ref{e3.23}).  
Due to (\ref{tt1tzi}), the following chain of equalities holds:
\begin{equation*}
w:= {\bf T}_\Theta[\varphi] = {\bf T}_\Theta[{\bf T}_{\hat\Theta}
[ \hat{\varphi}]]= {\bf T}_{\Theta \hat{\Theta}}[ \hat{\varphi}] =
{\bf T}_{\Theta^{(1)}}[\hat{\varphi}].
\end{equation*}
Thus, upon setting
\begin{equation*}
U_{\hat{\varphi}} = \Theta_{11}^{(1)}\hat{\varphi} +
\Theta_{12}^{(1)}
\quad\text{and}\quad
V_{\hat{\varphi}} = \Theta_{21}^{(1)}\hat{\varphi} +
\Theta_{22}^{(1)},
\end{equation*}
so that
\begin{equation}
\begin{bmatrix}U_{\hat{\varphi}} \\
V_{\hat{\varphi}}\end{bmatrix} =
\Theta^{(1)}
\begin{bmatrix}\hat{\varphi} \\ 1\end{bmatrix},
\label{6.49a}
\end{equation}
we arrive at
\begin{equation}
\w(z) = {\bf T}_{\Theta^{(1)}}[\hat{\varphi}](z) =
\frac{\Theta_{11}^{(1)}(z)\hat{\varphi}(z) + \Theta_{12}^{(1)}(z)}{
\Theta_{21}^{(1)}(z)\hat{\varphi}(z) + \Theta_{22}^{(1)}(z)} =
\frac{U_{\hat{\varphi}}(z)}{V_{\hat{\varphi}}(z)}.
\label{e4.125}
\end{equation}
Note that $\hat{\varphi}$ satisfies precisely one of the following:
\begin{align*}
&(a)\;\;\text{The limit $\hat{\varphi}(x_i)$ does not exist.}\\
&(b)\;\;\text{The limit $\hat{\varphi}(x_i)$ exists and is not
equal to $\eta_i = \infty$.}\\
&(c)\;\;\text{It holds that $\hat{\varphi}(x_i) = \eta_i = \infty$ and
$\hat{\varphi}_{-1}(x_i) = 0$.}
\end{align*}
The function $\Theta^{(1)}$ is
a rational function analytic and invertible at $x_i$, and 
so the limit $\w(x_i)$ exists if and only if
the limit $\hat{\varphi}(x_i)$ exists.  
Therefore $(1) \iff (a)$.
Henceforth,
let us assume that $\hat{\varphi}(x_i)$ exists.
We break the rest of the proof into two 
complementary cases:

\nline
{\bf Case I}: 
Suppose that $\hat{\varphi}(x_i) \neq \eta_i = \infty$, i.e.,
$\hat{\varphi}$ satisfies $(b)$.  
If $V_{\hat{\varphi}}(x_i) = 0$ and
if $U_{\hat{\varphi}}(x_i) \neq 0$, then $\w(x_i) = \infty \neq \w_i$.
On the other hand, if $V_{\hat{\varphi}}(x_i) = U_{\hat{\varphi}}(x_i) = 0$, then
\begin{equation*}
\Theta^{(1)}(x_i)\begin{bmatrix}\hat{\varphi}(x_i) \\ 1\end{bmatrix} =
\begin{bmatrix}U_{\hat{\varphi}}(x_i) \\
 V_{\hat{\varphi}}(x_i)\end{bmatrix} = 
\begin{bmatrix}0 \\ 0\end{bmatrix}, 
\end{equation*}
which implies that $\Theta^{(1)}(x_i)$ is singular, 
which is a contradiction.
Let us henceforth assume that $V_{\hat{\varphi}}(x_i) \neq 0$.

Using
(\ref{6RES}), Lemma \ref{L42}, and the analyticity of
$\Theta^{(1)}$ at $x_i$ to compare residues on
both sides of (\ref{tt1tzi}), we get
\begin{equation*}
\tilde{c}_i\begin{bmatrix}\w_i \\ 1\end{bmatrix} =
\Theta^{(1)}(x_i)\begin{bmatrix}\frac{\tilde{c}_i^2}{\tilde{p}_{ii}}
\\0 \end{bmatrix}.
\end{equation*}
By Lemma \ref{L2-1}, $\tilde{e}_i = 0$ implies $\tilde{c}_i \neq 0$,
and so after cancellation we have
\begin{equation}
\begin{bmatrix}\w_i \\ 1\end{bmatrix} =
\frac{1}{\tilde{p}_{ii}}\Theta^{(1)}(x_i)\begin{bmatrix}
\tilde{c}_i \\ 0\end{bmatrix}.
\label{e6.62}
\end{equation}
Upon multiplying (\ref{e6.62}) by $J$ on the left
and using $\Theta^{(1)}(x_i)^{-1} = J\Theta^{(1)}(x_i)J$,
we get
\begin{equation}
\begin{bmatrix}1 & -\w_i\end{bmatrix} = \frac{1}{\tilde{p}_{ii}}
\begin{bmatrix}0 & -\tilde{c}_i\end{bmatrix}
\Theta^{(1)}(x_i)^{-1}.
\label{arrive}
\end{equation}
Hence
\begin{equation}
\begin{bmatrix}1 & -\w_i\end{bmatrix}\Theta^{(1)}(x_i)
= \frac{1}{\tilde{p}_{ii}}\begin{bmatrix}0 &
  -\tilde{c}_i\end{bmatrix}.
\label{arrive2}
\end{equation}
Note that
\begin{equation*}
\w(x_i) - \w_i = \frac{U_{\hat{\varphi}}(x_i) - \w_i
  V_{\hat{\varphi}}(x_i)}{V_{\hat{\varphi}}(x_i)}
= \frac{1}{V_{\hat{\varphi}}(x_i)}
\begin{bmatrix}1 & -\w_i\end{bmatrix}
\begin{bmatrix}U_{\hat{\varphi}}(x_i) \\
V_{\hat{\varphi}}(x_i)\end{bmatrix}.
\end{equation*}
In view of (\ref{6.49a}), the latter equality becomes
\begin{equation*}
\w(x_i) - \w_i =
\frac{1}{V_{\hat{\varphi}}(x_i)} \begin{bmatrix}1 & -\w_i\end{bmatrix}
\Theta^{(1)}(x_i)\begin{bmatrix}\hat{\varphi}(x_i) \\ 1\end{bmatrix},
\end{equation*}
which, due to (\ref{arrive2}), yields
\begin{equation*}
\w(x_i) - \w_i =
\frac{-\tilde{c}_i}{\tilde{p}_{ii}V_{\hat{\varphi}}(x_i)} \neq 0.
\end{equation*}
We have thus shown that if $\hat{\varphi}(x_i) \neq \eta_i = \infty$,
then $\w(x_i)$ exists and $\w(x_i) \neq \w_i$.  
In other words, $(b) \implies (2)$.

{\bf Case II}:
Suppose that $\hat{\varphi}(x_i) = \infty$.
By (\ref{e4.120}) and (\ref{tt1tzi}), 
 $\Theta_{11}^{(1)}(z) \equiv \Theta_{11}(z)$ and
$\Theta_{21}^{(1)}(z) \equiv \Theta_{21}(z)$, and so (\ref{r6.35}) and 
(\ref{r6.36}) imply
\begin{equation}
\Theta_{11}^{(1)}(x_i) = \w_i \cdot \frac{\tilde{p}_{ii}}{\tilde{c}_i}
\quad\text{and}\quad
\Theta_{21}^{(1)}(x_i) = \frac{\tilde{p}_{ii}}{\tilde{c}_i}.
\label{t11t21lim}
\end{equation}
Since $\varphi$ satisfies $\tilde{\bC}_5$, it follows that 
$\tilde{p}_{ii} \neq 0$.  Therefore,
due to (\ref{e4.125}), (\ref{t11t21lim}), 
and the analyticity of
$\Theta^{(1)}(z)$ at $x_i$, 
we have
\begin{equation*}
\w(x_i) = \lim_{z \to x_i}\w(z) = 
\frac{\Theta_{11}^{(1)}(x_i)}{\Theta_{21}^{(1)}(x_i)} = \w_i.
\end{equation*}

To show that $(3)$ occurs, 
it remains to show
that 
\begin{equation}
\left|\lim_{z \to x_i}K_\w(z,z)\right| = \infty.
\label{6700}
\end{equation}
Substituting $\varphi_{-1}(x_i) = \frac{
\tilde{c}_i^2}{\tilde{p}_{ii}}$ into (\ref{r6.42}) yields,
after simplifications,
\begin{equation}
\lim_{z \to x_i}|z - x_i|^2\Psi_\varphi(z)P\Psi_\varphi(z) =
\frac{\tilde{c}_i^2}{\tilde{p}_{ii}}.
\label{6701}
\end{equation}
We remark that the formula (\ref{r6.42}) is valid
in the present situation since the only assumptions used
in deriving (\ref{r6.42}) were that $\tilde{e}_i = 0$ and
$\varphi \not \equiv \infty$.
Next, due to Lemma \ref{CJN2}, we have
\begin{equation}
\lim_{z \to x_i}|z - x_i|^2K_\varphi(z,z) =
\varphi_{-1}(x_i) = \frac{\tilde{c}_i^2}{\tilde{p}_{ii}}.
\label{6702}
\end{equation}
Formula (\ref{vlim6}) is also valid in the present situation,
since, as may be seen from the proof of the previous Theorem, 
this result ultimately only requires that $\tilde{e}_i = 0$.
Substituting $\varphi_{-1}(x_i) = \frac{\tilde{c}_i^2}{\tilde{p}_{ii}}$
into (\ref{vlim6}) yields
\begin{equation}
\lim_{z \to x_i}|z - x_i|^2V(x_i) = 0.
\label{6703}
\end{equation}
Since, by (\ref{6701}) and (\ref{6702}),
\begin{equation*}
\lim_{z \to x_i} |z - x_i|^2 \left(K_\varphi(z,z) +
\Psi_\varphi(z)P\Psi_\varphi(z)\right) =
2\frac{\tilde{c}_i^2}{\tilde{p}_{ii}} \neq 0,
\end{equation*}
using formulas (\ref{6701}) -- (\ref{6703}) in
(\ref{r6.7}) proves (\ref{6700}).
Thus, when $\hat{\varphi}(x_i) = \infty$, we have
$\w(x_i)$ and $|K_\w(x_i,x_i)| = \infty$.

We have shown that if $\hat{\varphi}(x_i) = \eta_i$
and $\hat{\varphi}_{-1}(x_i) = 0$, then
the limit $\w(x_i)$ exists and $\w(x_i) = \w_i$.  In
other words $(c) \implies (3)$.
Since $\w = \mathbf{T}_{\Theta^{(1)}}[\hat{\varphi}]$ and since
$\hat{\varphi}$ satisfies exactly one of $(a)$, $(b)$, $(c)$,
we have $(1) \iff (a)$, $(2) \iff (b)$, and $(3) \iff (c)$,
which completes the proof.
\end{proof}

\begin{Tm}
Let $i \in \{1, \ldots, \ell\}$ and suppose that $\tilde{e}_i = 0$.
Suppose that $\varphi \in \tN_0$ satisfies
condition $\btC_6$ at $x_i$.
Then the nontangential limit $\w(x_i)$ of $\w := {\bf
  T}_\Theta[\varphi]$ exists,
but is not equal to $\w_i$.
\label{T6.6}
\end{Tm}
\begin{proof}
The only function that may satisfy $\btC_6$ is $\varphi
\equiv \infty$.  Therefore,
\begin{equation}
\w \equiv \frac{\Theta_{11}}{\Theta_{21}}.
\label{wequivfrac}
\end{equation}
We shall use decompositions (\ref{PP}) and (\ref{XEC}) with the
understanding that $\tilde{p}_{ii} = 0$.  Thus
\begin{equation*}
\tilde{P}_{21}P_{12} = 1
\quad\text{and}\quad
P^{-1}\mathbf{e}_i =
\begin{bmatrix}\tilde{P}_{12} \\ 0\end{bmatrix}.
\end{equation*}
We shall also use the formula
\begin{equation*}
P_{21}(x_iI - X_1)^{-1} = \w_iE_1 - C_1
\end{equation*}
which follows 
upon substituting partitionings (\ref{PP}) and (\ref{XEC}) into
the Lyaponov identity (\ref{LI}) and comparing the
$(2,1)$ block entries.

By (\ref{r6.35}), (\ref{r6.36}), and the fact that $\tilde{p}_{ii} = 0$,
\begin{equation}
\Theta_{11}(x_i) := \lim_{z \to x_i}\Theta_{11}(z) = 0
\quad\text{and}\quad
\Theta_{21}(x_i) := \lim_{z \to x_i}\Theta_{21}(z) = 0. 
\label{6.84}
\end{equation}
Substituting (\ref{r6.18}) into (\ref{r6.20}) shows that
\begin{align*}
\Theta_{11}^\prime(z) &= C(zI - X)^{-2} \cdot
\frac{1}{\tilde{c}_i}(x_i I - X)P^{-1}\mathbf{e}_i ,\\
\Theta_{21}^\prime(z) 
&= E(zI-X)^{-2}\cdot \frac{1}{\tilde{c}_i}(x_iI-X)P^{-1}
\mathbf{e}_i.
\end{align*}
Taking limits yields
\begin{align}
\Theta_{11}^\prime(x_i) &=
\frac{1}{\tilde{c}_i}C_1(x_iI-X_1)^{-1}\tilde{P}_{12}, 
\label{6.86}\\
\Theta_{21}^\prime(x_i) &= \frac{1}{\tilde{c}_i}E_1
(x_iI-X_1)^{-1}\tilde{P}_{12}. \label{6.87}
\end{align}
Due to (\ref{6.84})--(\ref{6.87}), $\Theta_{11}$ and
$\Theta_{21}$ admit the Taylor series expansions
\begin{align*}
\Theta_{11}(z) &=
\frac{z - x_i}{\tilde{c}_i}C_1(x_iI-X_1)^{-1}\tilde{P}_{12} + O(|z - x_i|^2), \\
\Theta_{21}(z) &=  \frac{z - x_i}{\tilde{c}_i}E_1
(x_iI-X_1)^{-1}\tilde{P}_{12} + O(|z - x_i|^2),
\end{align*}
which, by (\ref{wequivfrac}), permit the following representation of $\w$:
\begin{equation*}
\w(z) = \frac{C_1(x_iI-X_1)^{-1}
\tilde{P}_{12} + O(|z - x_i|)}{
E_1(x_iI-X_1)^{-1}\tilde{P}_{12} + O(|z - x_i|)}.
\end{equation*}
Since $\w$ is rational, $|w(z)|$ has a limit at $x_i$.
When this limit is finite, it holds that
\begin{align*}
\w_i - \w(x_i) &= \w_i - \frac{C_1(x_iI - X_1)^{-1}
\tilde{P}_{12}}{E_1(x_iI - X_1)^{-1}\tilde{P}_{12}}
= \frac{(\w_i E_1 - C_1)(x_iI - X_1)^{-1}\tilde{P}_{12}}{
E_1(x_iI - X_1)^{-1}\tilde{P}_{12}} \\
&= \frac{P_{21}\tilde{P}_{12}}{E_1(x_iI - X_1)^{-1}
\tilde{P}_{12}}
= \frac{1}{E_1(x_iI - X_1)^{-1}\tilde{P}_{12}}.
\end{align*}
Now, if $E_1(x_iI - X_1)^{-1}\tilde{P}_{12} = 0$, then
we must have
$C_1(x_iI - X_1)^{-1}\tilde{P}_{12} = 0$ in order for
the limit $\w(x_i)$ to be finite.  However, this forces
$(\w_i E_1 - C_1)(X_1 - x_iI)^{-1}\tilde{P}_{12} = 0$, which is 
a contradiction by the above calculation.  Hence
$\w(x_i)$ exists but does not equal $\w_i$.
\end{proof}

\begin{Rk}
Theorems \ref{T6.4}, \ref{T6.5}, and
\ref{T6.6} together prove Theorem \ref{Main1}.
\end{Rk}

\chapter{Parameters and Interpolation Conditions II}

In Chapter VI, we initiated the characterization of
the boundary behavior of $\w = \mathbf{T}_\Theta[\varphi]$
around interpolation nodes $x_i$ 
(in terms of the behavior of the parameter $\varphi$ around these nodes)
in order to determine
under what conditions $\w$ would satisfy some or 
all of the interpolation conditions 
stated in Problem \ref{Prob1} or in Problem \ref{Prob2}, i.e.,
the conditions
(\ref{1.26})--(\ref{1.27c}) or
the conditions (\ref{1.29}) and (\ref{1.30}), respectively.
In particular, in Chapter VI we completed the characterization
around interpolation nodes $x_i$ for which $\tilde{e}_i = 0$
(i.e., $\eta_i = \infty$).  In this chapter we consider
the remaining case by
characterizing
the boundary behavior of $\w$ around interpolation nodes
$x_i$ for which $\tilde{e}_i \neq 0$ (i.e., $\eta_i \neq \infty$).

\section{Interpolation Conditions at Singular Interpolation Nodes}
Throughout this section we let $i \in \{\ell + 1, \ldots, n\}$
be fixed.  The introductory material in 
and preceding Section 6.1 is relevant here,
and we shall frequently refer to it.

\begin{Tm}
Let $i \in \{\ell + 1, \ldots, n\}$ and suppose that
$\tilde{e}_i \neq 0$ (i.e., $\eta_i \neq \infty$).  
Let $\varphi \in \tN_0$ 
and let $\w := \mathbf{T}_\Theta[\varphi]$.  If $\varphi$
satisfies condition $\bC_1$ at $x_i$, then
\begin{equation*}
\w_{-1}(x_i) = \xi_i.
\end{equation*}
\end{Tm}
\begin{proof}
By Lemma \ref{CJN2}, a function $\varphi \in \N_0$ satisfies
exactly one of the following:
\begin{equation*}
(a)\;\; \varphi_{-1}(x_i) = 0
\qquad\text{or}\qquad
(b)\;\; -\infty < \varphi_{-1}(x_i) < 0.
\end{equation*}
Suppose that the parameter $\varphi$ satisfies $(a)$.
By (\ref{r6.3}), (\ref{r6.13}), and (\ref{r6.1}),
we have
\begin{align}
(z - x_i)\w(z) &= (z - x_i)\frac{U_\varphi(z)}{V_\varphi(z)} \nonumber
\\
&=
\xi_i \frac{\Theta_{21}(z)\varphi(z) + \Theta_{22}(z)}{
\Theta_{21}(z)\varphi(z) + \Theta_{22}(z)} + 
(z - x_i)\frac{\tilde{\Theta}_{11}^{(2)}(z)\varphi(z) +
\tilde{\Theta}_{12}^{(2)}(z)}{
\Theta_{21}(z)\varphi(z) + \Theta_{22}(z)}.
\label{s7.1}
\end{align}
By (\ref{r6.15}) and (\ref{r6.16}), 
\begin{equation*}
-\lim_{z \to x_i} \frac{\Theta_{22}(z)}{\Theta_{21}(z)} =
 \frac{\tilde{c}_i}{\tilde{e}_i} = \eta_i,
\end{equation*}
which, combined with the fact that $\varphi$ meets condition
$\bC_1$ at $x_i$, shows that 
there exists an $\epsilon > 0$ and a sequence $\{z_\alpha\}$
such that
\begin{equation}
\left|\Theta_{21}(z_\alpha)\varphi(z_\alpha) + \Theta_{22}(z_\alpha)
\right| > \epsilon
\quad \text{for all}\quad
\alpha.
\label{s7.1a}
\end{equation}
By Remark \ref{PoleRemark}, the functions
$\tilde{\Theta}_{11}^{(2)}$ and $\tilde{\Theta}_{12}^{(2)}(z)$
are analytic at $x_i$.
In particular, 
\begin{equation}
\lim_{z \to x_i}(z - x_i)\tilde{\Theta}_{12}^{(2)}(z) = 0.
\label{s7.1b}
\end{equation}
Now (\ref{s7.1}), (\ref{s7.1a}), and (\ref{s7.1b}) show
that
\begin{equation}
\lim_{\alpha \to \infty}(z_\alpha - x_i)\w(z_\alpha) = \xi_i.
\label{s7.1c}
\end{equation}
Theorem \ref{CJGN2} applied to $\w$ and (\ref{s7.1c}) shows
that in fact that limit
\begin{equation*}
\w_{-1}(x_i) := \lim_{z \to x_i}\w(z)
\end{equation*}
exists and equals $\xi_i$.
Therefore the Theorem is proved for any $\varphi \in \N_0$
that satisfies Case 1.

Let us now suppose that $\varphi \in \N_0$ satisfies
$(b)$, i.e.,
\begin{equation}
-\infty < \varphi_{-1}(x_i) := \lim_{z \to x_i}(z - x_i)\varphi(z) < 0,
\label{phicase2}
\end{equation}
Instead of (\ref{s7.1}), we write
\begin{align}
(z - x_i)\w(z) =& (z - x_i)\frac{U_\varphi(z)}{V_\varphi(z)} \nonumber
\\
=&
\xi_i\frac{\Theta_{21}(z)\varphi(z)(z - x_i) + (z - x_i)\Theta_{22}(z)}{\Theta_{21}(z)
\varphi(z)(z - x_i) + (z - x_i)\Theta_{22}(z)} \nonumber \\
&+
(z - x_i)\frac{\tilde{\Theta}_{11}^{(2)}(z)\varphi(z)(z - x_i) + 
(z - x_i)\tilde{\Theta}_{12}^{(2)}(z)}{\Theta_{21}(z)\varphi(z)(z - x_i) +
(z - x_i)\Theta_{22}(z)},
\label{s7.4}
\end{align}
Due to (\ref{r6.15}) and (\ref{r6.16}),
\begin{equation*}
\lim_{z \to x_i}\Theta_{21}(z) = \tilde{e}_i
\quad\text{and}\quad
\lim_{z \to x_i}(z - x_i)\Theta_{22}(z) = 0.
\end{equation*}
Therefore
\begin{equation}
\lim_{z \to x_i}\left((z - x_i)\Theta_{21}(z) \varphi(z) +
(z - x_i)\Theta_{22}(z)\right) 
= \tilde{e}_i \varphi_{-1}(x_i) \neq 0.
\label{s7.5}
\end{equation}
Since $\tilde{\Theta}_{11}^{(2)}$ and $\tilde{\Theta}_{12}^{(2)}$
are analytic at $x_i$ (see Remark \ref{PoleRemark}), 
it follows that
\begin{equation}
\lim_{z \to x_i}\left(
(z - x_i)\tilde{\Theta}_{11}^{(2)}(z)\varphi(z) +
(z - x_i)\tilde{\Theta}_{12}^{(2)}(z)\right) = \tilde{e}_i\varphi_{-1}(x_i).
\label{s7.6}
\end{equation}
Pass to limits in (\ref{s7.4}).  On account of (\ref{s7.5})
and (\ref{s7.6}), the result is $\w_{-1}(x_i) = \xi_i$.
Therefore the Theorem is proved for $\varphi \in \N_0$
that satisfy Case 2.

The only case left to consider the case where
$\varphi \equiv \infty$.  Together (\ref{r6.7a}),
(\ref{r6.12}), and (\ref{r6.12a}) show that
\begin{equation}
(z - x_i)\w(z) = (z - x_i)\frac{U_\infty(z)}{V_\infty(z)}
= \xi_i \frac{\Theta_{21}(z)}{\Theta_{21}(z)} +
(z - x_i)\frac{\tilde{\Theta}_{11}^{(2)}(z)}{\Theta_{21}(z)}.
\label{s7.7}
\end{equation}
By Lemma \ref{L42},
\begin{equation*}
\lim_{z \to x_i}(z - x_i)\Theta_{21}(z) = \tilde{e}_i \neq 0.
\end{equation*}
This combined with the analyticity of
$\tilde{\Theta}_{11}^{(2)}$ at $x_i$ (see Remark \ref{PoleRemark})
shows that taking limits in (\ref{s7.7}) yields
$\w_{-1}(x_i) = \xi_i$,
which completes the proof.
\end{proof}

\begin{La}
Let $x_0 \in \mathbb{R}$ and let $\varphi \in \N_0$ be such
that
\begin{equation}
\lim_{z \to x_0}\varphi(z) = \varphi_0 \in \mathbb{R}
\quad\text{and}\quad
\lim_{z \to x_0}K_\varphi(z,z) = \infty.
\label{EastGermany}
\end{equation}
Then
\begin{equation}
\lim_{z \to x_0} K_\varphi(z,z) \cdot
\left|\frac{z - x_0}{\varphi(z) - \varphi_0}\right|^2 = 0
\quad\text{and}\quad
\lim_{z \to x_0} \frac{z - x_0}{\varphi(z) - \varphi_0} = 0.
\label{WestGermany}
\end{equation}
\label{slemma}
\end{La}
\begin{proof}
First note that
\begin{equation}
|\varphi(z) - \varphi_0| >  \Im (\varphi(z) - \varphi_0) =
\frac{\varphi(z) - \varphi(z)^*}{2i} =
\frac{1}{2}(\varphi(z) - \varphi(z)^*)
\label{Poland}
\end{equation}
If $D \subset \C^+$ is a nontangential neighborhood of $x_i$,
then we may find an $a > 0$ such that the domain
\begin{equation}
\Gamma_a := \left\{z \in \C^+: a(z - \bar{z}) > |z - x_i| \right\}
\label{Belgium}
\end{equation}
contains $D$.  Note that if $z$ tends to $x_i$ staying
inside $\Gamma_a$, then $z$ tends to $x_i$ nontangentially.
Let us assume now that $\Gamma_a$ is fixed and that
$z \in \Gamma_a$.  Then due
to (\ref{Belgium})
\begin{equation}
\frac{|z - x_i|}{z - \bar{z}} < a.
\label{Sweden}
\end{equation}
Combining (\ref{Poland}) and (\ref{Sweden}), we have
\begin{align}
K_\varphi(z,z) \cdot \left|\frac{z - x_0}{\varphi(z) -
    \varphi_0}\right|
&= \frac{\varphi(z) - \varphi(z)^*}{z - \bar{z}} \cdot
 \left|\frac{z - x_0}{\varphi(z) -
    \varphi_0}\right| \nonumber \\
&=
\frac{\varphi(z) - \varphi(z)^*}{|\varphi(z) - \varphi_0|}
\cdot \frac{|z - x_i|}{z - \bar{z}} < 2a.
\label{Ukraine}
\end{align}
In view of the second equality in (\ref{EastGermany}),
(\ref{Ukraine}) implies the second equality in
(\ref{WestGermany}).  The first equality in
(\ref{WestGermany}) thus follows as well.
\end{proof}

\begin{Tm}
Let $i \in \{\ell +1, \ldots, n\}$ and suppose that
$\tilde{e}_i \neq 0$ (i.e., $\eta_i \neq \infty$).  
Let $\varphi \in \tN_0$ 
and let $\w := \mathbf{T}_\Theta[\varphi]$.  If $\varphi$
satisfies condition $\bC_2$ at $x_i$, then
\begin{equation*}
\w_{-1}(x_i) = \xi_i.
\end{equation*}
\end{Tm}
\begin{proof}

By the definition (\ref{r2.2}) of $X$, we have
\begin{equation}
(z - x_i)(zI - X)^{-1} = {\bf e}_i{\bf e}_i^* + O(|z - x_i|)
\label{Korea}
\end{equation}
for $z$ near $x_i$.
Since
$\varphi(z)$ is bounded in a
neighborhood of $x_i$, using
(\ref{Korea}) in (\ref{r6.4}) shows that
\begin{align*}
(z - x_i)\Psi_\varphi(z) &= (z - x_i)(zI - X)^{-1}\left(
\tC^* - \varphi(z)\tE^*\right) \\
&= {\bf e}_i{\bf e}_i^*\left(
\tC^* - \varphi(z)\tE^*\right) + O(|z - x_i|)
\end{align*}
holds asymptotically.

Using the latter relation in 
(\ref{r6.4a}) yields
and taking into account
${\bf e}_i^*\tE^* = \tilde{e}_i$
and
${\bf c}_i^*\tC^* = \tilde{c}_i$
yields
the asymptotic relation
\begin{equation}
(z - x_i)U_\varphi(z) = \xi_i(\tilde{c}_i - \varphi(z)\tilde{e}_i)
+ O(|z - x_i|).
\label{Vietnam}
\end{equation} 
By (\ref{r6.15}) and (\ref{r6.16}),
\begin{equation*}
\Theta_{21}(x_i) = \tilde{e}_i
\quad\text{and}\quad
\Theta_{22}(x_i) = -\tilde{c}_i,
\end{equation*}
and so $\Theta_{21}$ and $\Theta_{22}$ admit 
at $x_i$ the Taylor expansions
\begin{equation*}
\Theta_{21}(z) = \tilde{e}_i + O(|z - x_i|)
\quad\text{and}\quad
\Theta_{22}(z) = -\tilde{c}_i + O(|z - x_i|).
\end{equation*}
Using the latter expansions in (\ref{r6.1})
yields
\begin{align}
V_\varphi(z) &= \left(\tilde{e}_i + O(|z - x_i|)\right)\varphi(z)
- \tilde{c}_i + O(|z - x_i|) \nonumber\\
&= \tilde{e}_i\varphi(z) - \tilde{c}_i + (\tilde{e}_i\varphi(z)
- \tilde{c}_i)O(|z - x_i|) + O(|z - x_i|).
\label{Congo}
\end{align}
Let us set
\begin{equation*}
\Delta_i(z) := \frac{\tilde{c}_i - \tilde{e}_i\varphi(z)}{
z - x_i}
\end{equation*}
for short and note that
\begin{equation*}
\Delta_i(z) \neq 0
\quad (z \in \C^+).
\end{equation*}
For, if $\varphi(z_0) = \eta_i \in \mathbb{R}$ at some
point $z_0 \in \C^+$, then $\varphi \equiv \eta_i$, which
would contradict the fact that $\varphi$ satisfies
$\bC_2$.  Indeed, since $\varphi$ does satisfy $\bC_2$ at $x_i$, 
Lemma \ref{slemma} is applicable with $\varphi_0 = \eta_i$,
and shows that
\begin{equation}
\lim_{z \to x_0}K_\varphi(z,z) \cdot \frac{1}{|\Delta_i(z)|^2}
= 0
\label{d1}
\end{equation}
and
\begin{equation}
\lim_{z \to x_0}\Delta_i(z)^{-1} = 0.
\label{d2}
\end{equation}
Note that
\begin{equation*}
(z - x_i)\w(z) = (z - x_i)\frac{U(z)}{V_\varphi(z)} =
\frac{\Delta_i(z)^{-1}U(z)}{V_\varphi(z)/(\tilde{c}_i - \varphi(z)\tilde{e}_i)}.
\end{equation*}
Substituting (\ref{Vietnam}) and (\ref{Congo}) into the latter
identity shows that
\begin{equation*}
(z - x_i)\w(z) = \frac{\xi_i + \Delta_i(z)^{-1} \cdot O(1)}{1
+ O(|z - x_i|) + \Delta_i(z)^{-1} \cdot O(1)}
\end{equation*}
and taking limits yields
\begin{equation*}
\w_{-1}(x_i) := \lim_{z \to x_i}(z - x_i)\w(z) = \xi_i.
\end{equation*}
\end{proof}

\begin{Tm}
Let $i \in \{\ell +1, \ldots, n\}$ and suppose that
$\tilde{e}_i \neq 0$ (i.e., $\eta_i \neq \infty$).  
Let $\varphi \in \tN_0$ 
and let $\w := \mathbf{T}_\Theta[\varphi]$.
\begin{enumerate}
\item If $\varphi$ satisfies $\bC_3$, then
$-\infty < -\w_{-1}^{-1}(x_i) < -\xi_i^{-1}.$
\item If $\varphi$ satisfies $\bC_4$, then
$-\xi_i^{-1} < -\w_{-1}^{-1}(x_i) < \infty.$
\item If $\varphi$ satisfies $\bC_{5-6}$, then
$\w_{-1}(x_i) = 0$.
\end{enumerate}
\end{Tm}
\begin{proof}
Let us assume that $\varphi \in \N_0$ satisfies
$\bC_{3-6}$.  By Theorem \ref{CJGN1}, the following limits
exist:
\begin{equation*}
\varphi^\prime(x_i) :=
\lim_{z \to x_i}\varphi^\prime(z) =
\lim_{z \to x_i}\frac{\varphi(z) - \eta_i}{z - x_i} =
\lim_{z \to x_i}K_\varphi(z,z).
\end{equation*}
By Corollary $7.9$ in \cite{DB},
the following asymptotic equality
holds:
\begin{equation}
\varphi(z) = \eta_i + (z - x_i)\varphi^\prime(x_i) +
o(|z - x_i|)
\label{r72.1}
\end{equation}
for $z \in \C^+$ close enough to $x_i$.
We shall use this asymptotic relation to show that the functions
$\Psi_\varphi$, $U_\varphi$, and $V_\varphi$,
defined (\ref{r6.4}) and (\ref{r6.4a}),
admit the
following nontangential boundary limits at the interpolation
node $x_i$:
\begin{equation}
\Psi_\varphi(x_i) = \frac{1}{\tilde{e}_i}
\left(P^{-1}\mathbf{e}_i - \mathbf{e}_i
(\tilde{p}_{ii} + \tilde{e}_i^2\varphi^\prime(x_i))\right),
\label{r72.2}
\end{equation}
\begin{equation}
U_\varphi(x_i) = \frac{\xi_i}{\tilde{e}_i}
(\tilde{p}_{ii} + \tilde{e}_i^2\varphi^\prime(x_i)),
\quad\text{and}\quad
V_\varphi(x_i) = 0
\label{r72.3}
\end{equation}
Multiply both sides of the
Lyaponov identity (\ref{e3.4}) by $\mathbf{e}_i$ on the right
to obtain
\begin{equation*}
XP^{-1}\mathbf{e}_i - P^{-1}X\mathbf{e}_i
= \tE^*\tilde{c}_i - \tC^*\tilde{e}_i,
\end{equation*}
which may be rewritten as
\begin{equation}
\tC^* - \eta_i\tE^* = \frac{1}{\tilde{e}_i}
(x_i I - X)P^{-1}\mathbf{e}_i.
\label{r72.4}
\end{equation}
Substitute (\ref{r72.1}) into (\ref{r6.4}):
\begin{align*}
\Psi_\varphi(z) &=
-(zI - X)^{-1}\left\{
\left[ \eta_i + (z - x_i)\varphi^\prime(x_i)
+ o(|z - x_i|)\right] \tE^* - \tC^* \right\}\\
&= -(zI - X)^{-1}\left(\eta_i \tE^* - \tC^*\right)
- (zI - X)^{-1}(z - x_i)\varphi^\prime(x_i)\tE^*
+ o(1).
\end{align*}
Substituting (\ref{r72.4}) into the latter equality yields
\begin{equation}
\Psi_\varphi(z) =
\frac{1}{\tilde{e}_i} (zI - X)^{-1}(x_iI - X)P^{-1}\mathbf{e}_i
- (zI - X)^{-1}(z - x_i)\varphi^\prime(x_i)\tE^* + o(1).
\label{r72.5}
\end{equation}
Using (\ref{r72lim1}) and (\ref{r72lim2}), we pass
to limits in (\ref{r72.5}) to obtain
\begin{align}
\Psi_\varphi(x_i) &= \frac{1}{\tilde{e}_i}(I - {\bf e}_i{\bf e}_i^*)
P^{-1}{\bf e}_i - {\bf e}_i{\bf e}_i^*\varphi^\prime(x_i)\tilde{E}^* 
\nonumber\\
&= \frac{1}{\tilde{e}_i}(P^{-1}{\bf e}_i - {\bf e}_i\tilde{p}_{ii})
- {\bf e}_i\tilde{e}_i\varphi^\prime(x_i) \nonumber \\
&= \frac{1}{\tilde{e}_i}\left(P^{-1}{\bf e}_i - {\bf e}_i(
\tilde{p}_{ii} + \tilde{e}_i^2\varphi^\prime(x_i))\right),
\label{PsiLimit}
\end{align}
where we have used the fact that
${\bf e}_i^*P^{-1}{\bf e}_i = \tilde{p}_{ii}$ and
${\bf e}_i^*\tilde{E}^* = \tilde{e}_i$.
This proves (\ref{r72.2}).

Using (\ref{r72.2}), we may pass to limits in (\ref{r6.4a}).
We get
\begin{align*}
U_\varphi(x_i) &= \varphi(x_i) - C\Psi_\varphi(x_i) \\
&= \eta_i - C \cdot \frac{1}{\tilde{e}_i}\left(P^{-1}{\bf e}_i
-{\bf e}_i(\tilde{p}_{ii} + \tilde{e}_i^2\varphi^\prime(x_i))\right) \\
&= \eta_i - \eta_i + \frac{\xi_i}{\tilde{e}_i}(\tilde{p}_{ii} +
\tilde{e}_i^2\varphi^\prime(x_i)) \\
&= \frac{\xi_i}{\tilde{e}_i}(\tilde{p}_{ii} +
\tilde{e}_i^2\varphi^\prime(x_i)),
\end{align*}
and, similarly,
\begin{align*}
V_\varphi(x_i) = 1 - E\Psi_\varphi(x_i)
&= 1 - E  \cdot \frac{1}{\tilde{e}_i}\left(P^{-1}{\bf e}_i
-{\bf e}_i(\tilde{p}_{ii} + \tilde{e}_i^2\varphi^\prime(x_i))\right) \\
&= 1 - 1 + 0 = 0,
\end{align*}
thus proving (\ref{r72.3}).

We now claim that the derivative of $V_\varphi$ satisfies
\begin{equation}
V_\varphi^\prime(x_i) := \lim_{z \to x_i}V_\varphi(z) =
\frac{1}{\xi_i\tilde{e}_i} + \frac{\tilde{p}_{ii}}{
\tilde{e}_i} + \tilde{e}_i\varphi^\prime(x_i).
\label{r72.6}
\end{equation}
Towards this end, we first differentiate $V$ with
respect to $z$:
\begin{equation}
V_\varphi^\prime= \left(\Theta_{21}\varphi +
  \Theta_{22}\right)^\prime
= \Theta_{21}^\prime\varphi +
\Theta_{22}^\prime + \Theta_{21}\varphi^\prime,
\label{r72.6a}
\end{equation}
where, on account of (\ref{ThetaFormula2}), we
have
\begin{equation}
\Theta_{21}^\prime(z) = -E(zI-X)^{-2}P^{-1}E^*
\quad\text{and}\quad
\Theta_{22}^\prime(z) = E(zI-X)^{-2}P^{-1}C^*.
\label{r72.7}
\end{equation}
Together (\ref{r72.1}) and (\ref{r72.7}) yield
\begin{align}
\Theta_{21}^\prime(z)\varphi(z) + \Theta_{22}^\prime(z)
=& E(zI - X)^{-2}(\tC^* - \eta_i \tE^*) \nonumber \\
&-E(zI - X)^{-2}(z - x_i)\varphi^\prime(x_i)\tE^* \nonumber \\
&-E(zI - X)^{-2}\tE^* \cdot o(|z - x_i|).
\label{r72.8}
\end{align}
Combining (\ref{r6.14}) and (\ref{r72.4}) gives
\begin{equation*}
E(zI - X)^{-2}(\tC^* - \eta_i \tE^*) =
\frac{1}{\xi_i \tilde{e}_i} \mathbf{e}_i^* P
(x_i I - X)(zI - X)^{-2}(x_i I - X)P^{-1}\mathbf{e}_i.
\end{equation*}
Using (\ref{r72lim1}), we may pass to limits in the
latter equality to obtain
\begin{equation}
\frac{1}{\xi_i\tilde{e}_i}\mathbf{e}_i^*P(I -
\mathbf{e}_i\mathbf{e}_i^*)^2P^{-1}\mathbf{e}_i
= \frac{1 + \xi_i\tilde{p}_{ii}}{\xi_i\tilde{e}_i}.
\label{r72.9}
\end{equation}
On account of (\ref{r6.14}), we have
\begin{align*}
E(zI - X)^{-2}(z - x_i)\varphi^\prime(x_i)\tE^* &=
\frac{1}{\xi_i}\mathbf{e}_i^*P(x_iI - X)(zI - X)^{-2}
(z - x_i)\varphi^\prime(x_i)\tE^*\\
E(zI - X)^{-2}\tE^* \cdot o(|z - x_i|) &=
\frac{1}{\xi_i}\mathbf{e}_i^*P(x_iI - X)(zI - X)^{-2}
\tE^* \cdot o(|z - x_i|) \\
&=\frac{1}{\xi_i}\mathbf{e}_i^*P(x_iI - X)(zI-X^{-2}
(z - x_i)\tE^* \cdot o(1).
\end{align*}
Taking limits in the latter equalities, we get, 
due to (\ref{r72lim1}) and (\ref{r72lim2}),
\begin{align}
\lim_{z \to x_i} E(zI - X)^{-2}(z - x_i)\varphi^\prime(x_i)\tE^*
&= \frac{1}{\xi_i}\mathbf{e}_i^*P(I - \mathbf{e}_i\mathbf{e}_i^*)
\mathbf{e}_i\mathbf{e}_i^*\varphi^\prime(x_i)\tE^* \nonumber \\
&=\frac{1}{\xi_i}\mathbf{e}_i^*P(\mathbf{e}_i\mathbf{e}_i^*-
\mathbf{e}_i\mathbf{e}_i^*)\varphi^\prime(x_i)\tE^* = 0,
\label{r72.10}
\end{align}
and, similarly
\begin{equation}
\lim_{z \to x_i} E(zI - X)^{-2}\tE^* \cdot o(|z - x_i|) = 0.
\label{r72.11}
\end{equation}
Using (\ref{r72.9})--(\ref{r72.11}), we may pass to limits
in (\ref{r72.8}) to get:
\begin{equation}
\lim_{z \to x_i} \left(\Theta_{21}^\prime(z)\varphi(z) + 
\Theta_{22}^\prime(z)\right) =
\frac{1}{\xi_i\tilde{e}_i} + \frac{\tilde{p}_{ii}}{\tilde{e}_i}.
\label{Denmark}
\end{equation}
Due to (\ref{r6.15}),
\begin{equation}
\lim_{z \to x_i}\Theta_{21}(z)\varphi^\prime(z) =
\tilde{e}_i\varphi^\prime(x_i).
\label{r72.12}
\end{equation}
Substituting (\ref{Denmark}) and (\ref{r72.12})
into (\ref{r72.6a}) and taking limits proves
(\ref{r72.6}).

By (\ref{r6.1}), (\ref{r72.1}), (\ref{r72.3}), and (\ref{r72.6}),
the function $V_\varphi$ admits the asymptotic formula
\begin{equation}
V_\varphi(z) = (z - x_i)V_\varphi^\prime(x_i) + o(|z - x_i|).
\label{r72.13}
\end{equation}
Combining (\ref{r72.13}) and (\ref{r72.6}) shows that
\begin{equation}
\lim_{z \to x_i}\frac{1}{z - x_i}V_\varphi(z) = V_\varphi^\prime(x_i) =
\frac{1}{\xi_i \tilde{e}_i} + \frac{\tilde{p}_{ii}}{\tilde{e}_i}
+ \tilde{e}_i \varphi_\varphi^\prime(x_i).
\label{r72.14}
\end{equation}
Thus, by (\ref{r6.3}), (\ref{r72.3}), and (\ref{r72.14}),
\begin{align}
\w_{-1}(x_i) := \lim_{z \to x_i}(z - x_i)\w(z) &=
\lim_{z \to x_i}\frac{U_\varphi(z)}{(z - x_i)^{-1}V_\varphi(z)} \nonumber \\
&= \frac{\xi_i^2(\tilde{p}_{ii} + \tilde{e}_i^2\varphi^\prime(x_i))}{
1 + \xi_i(\tilde{p}_{ii} + \tilde{e}_i^2\varphi^\prime(x_i))}.
\label{r72.15}
\end{align}

If $\varphi$ meets condition $\bC_{5-6}$ at $x_i$, then
$\tilde{p}_{ii} + \tilde{e}_i^2 \varphi^\prime(x_i) = 0$,
and so (\ref{r72.15}) reveals that $\w_{-1}(x_i) = 0$.

Now suppose that $\varphi$ meets condition $\bC_3$ or
$\bC_4$ at $x_i$.  In either case, it follows that
$\tilde{p}_{ii} + \tilde{e}_i^2\varphi^\prime(x_i) \neq 0$.
Thus, from (\ref{r72.15}), we get
\begin{equation*}
\frac{1}{\w_{-1}(x_i)} - \frac{1}{\xi_i} =
\frac{1}{\xi_i^2(\tilde{p}_{ii} + \tilde{e}_i^2\varphi^\prime(x_i))}.
\end{equation*}
It therefore follows that if $\varphi$ meets $\bC_3$ at $x_i$, then
$\w_{-1}^{-1}(x_i) < \xi_i^{-1}$, and if $\varphi$ meets $\bC_4$
at $x_i$, then $\w_{-1}^{-1}(x_i) > \xi_i^{-1}$.
\end{proof}

\section{Interpolation Conditions at Regular Interpolation Nodes}

Throughout this section we let $i \in \{1, \ldots, \ell\}$ be 
fixed.  The introductory material preceding Section 6.1 is
relevant here, and we shall frequently refer to it.

\begin{La}
Let $\varphi \in \N_0$, and let $i \in \{1, \ldots, \ell\}$.
Suppose that $\varphi_{-1}(x_i) = 0$.  
If the nontangential limit 
$\limsup_{z \to x_i}|\varphi(x_i)| < \infty$,  
then the following asymptotic relations
hold as $z$ tends to $x_i$ nontangentially:
\begin{align}
(z - x_i)\Psi_\varphi(z) &= {\bf e}_i(\tilde{c}_i - \varphi(z)\tilde{e}_i) +
O(|z - x_i|),
\label{e4.25} \\
(z - x_i)U_\varphi(z) &= \w_i(\tilde{c}_i - \varphi(z)\tilde{e}_i) +
O(|z - x_i|),
\label{e4.26} \\
(z - x_i)V_\varphi(z) &= (\tilde{c}_i - \varphi(z)\tilde{e}_i) +
O(|z - x_i|).
\label{e4.27}
\end{align}
If ${\displaystyle \limsup_{z \to x_i}|\varphi(x_i)| = \infty}$, then the following
asymptotic relations hold as $z$ tends to $x_i$ nontangentially:
\begin{align}
(z - x_i)\Psi_\varphi(z) &= {\bf e}_i(\tilde{c}_i - \varphi(z)\tilde{e}_i) +
o(1),
\label{e4.250}\\
(z - x_i)U_\varphi(z) &= \w_i(\tilde{c}_i - \varphi(z)\tilde{e}_i) +
o(1),
\label{e4.260}\\
(z - x_i)V_\varphi(z) &= (\tilde{c}_i - \varphi(z)\tilde{e}_i) +
o(1).
\label{e4.270}
\end{align}
\label{l4.5}
\end{La}
\begin{proof}
By the definition (\ref{r2.2}) of $X$, we have
\begin{equation}
(z - x_i)(zI - X)^{-1} = {\bf e}_i{\bf e}_i^* + O(|z - x_i|)
\label{e4.28}
\end{equation}
for $z$ near $x_i$.
Thus,
if $\varphi(z)$ is bounded in a
neighborhood of $x_i$, then
\begin{align*}
(z - x_i)\Psi_\varphi(z) &= (z - x_i)(zI - X)^{-1}\left(
\tilde{C}^* - \varphi(z)\tilde{E}^*\right) \\
&= {\bf e}_i{\bf e}_i^*\left(
\tilde{C}^* - \varphi(z)\tilde{E}^*\right) + O(|z - x_i|)
\end{align*}
holds asymptotically, 
which proves (\ref{e4.25}) since ${\bf e}_i^*\tilde{C}^* = \tilde{c}_i$ and
${\bf e}_i^*\tilde{E}^* = \tilde{e}_i$.
Using (\ref{e4.25}) in (\ref{r6.4a}) in turn
yields (\ref{e4.26}). 
If $\varphi(x)$ is not bounded in some neighborhood of $x_i$ and yet
\begin{equation*}
\varphi_{-1}(x_i) := \lim_{z \to x_i} (z - x_i)\varphi(x_i) = 0, 
\end{equation*}
then $O(|z - x_i|) \cdot \varphi(z)$ is
$o(1)$ asymptotically, which proves
the asymptotic relation
\begin{equation}
(z - x_i)\Psi_\varphi(z) =
{\bf e}_i{\bf e}_i^*\left(
\tilde{C}^* - \varphi(z)\tilde{E}^*\right) + o(1).
\label{e6.23a}
\end{equation}
Substituting (\ref{e6.23a}) into (\ref{r6.4a}) then yields
(\ref{e4.260}) and (\ref{e4.270}).
\end{proof}

\begin{La}
Let $\varphi \in \N_0$
and let $i \in \{1, \ldots, \ell\}$.  Suppose that
$\varphi(x_i) < 0$.  Then the following limits exist:
\begin{align}
\lim_{z \to x_i}(z - x_i)^2\Psi_\varphi(z) &= -{\bf
  e}_i\tilde{e}_i\varphi_{-1}(x_i), 
\label{bull1}\\
\lim_{z \to x_i}(z - x_i)^2 U_\varphi(z) &= \w_i\varphi_{-1}(x_i)\tilde{e}_i,
\label{e4.68} \\
\lim_{z \to x_i}(z - x_i)^2 V_\varphi(z) &= \varphi_{-1}(x_i)\tilde{e}_i.
\label{e4.69}
\end{align}
\label{l4.8}
\end{La}
\begin{proof}
By formula (\ref{r6.4}) for $\Psi_\varphi$, we may write
\begin{equation*}
(z - x_i)^2\Psi_\varphi(z) = (z - x_i)(zI - X)^{-1}\left((z - x_i)\tilde{C}^*  -
(z - x_i)\varphi(z)\tilde{E}^*\right).
\end{equation*}
Using (\ref{e4.28}) and the assumption
that $-\infty < \varphi_{-1}(x_i) < 0$, we obtain the asymptotic relation
\begin{align}
(z - x_i)^2\Psi_\varphi(z) &= \left({\bf e}_i{\bf e}_i^* +O(|z - x_i|)\right)
\left( -\varphi_{-1}(x_i)\tilde{e}_i + O(|z - x_i|)\right) \nonumber \\
&= -{\bf e}_i\tilde{e}_i\varphi_{-1}(x_i) + O(|z - x_i|).
\label{e4.72}
\end{align}
Taking limits as $z \to x_i$ in (\ref{e4.72}) yields (\ref{bull1}).
Substituting (\ref{bull1}) into (\ref{r6.4a})
yields (\ref{e4.68}) and (\ref{e4.69}).
\end{proof}

\begin{Tm}
Let $i \in \{1, \ldots, \ell\}$ and suppose that
$\tilde{e}_i \neq 0$ (i.e., $\eta_i \neq \infty$).
Let $\varphi \in \tN_0$ 
and let $\w := \mathbf{T}_\Theta[\varphi]$.  If $\varphi$
satisfies condition $\bC_1$ at $x_i$, then
\begin{equation*}
\w(x_i) = \w_i
\quad\text{and}\quad
\w^\prime(x_i) = \gamma_i.
\end{equation*}
\end{Tm}
\begin{proof}
Let us first suppose that $\varphi \in \N_0$ and
that $\varphi_{-1}(x_i) = 0$.
Then $\varphi$ is subject to either
\begin{equation*}
(a)\;\; \limsup_{z \to x_i}|\varphi(z)| < \infty
\qquad\text{or}\qquad
(b)\;\; \limsup_{z \to x_i}|\varphi(z)| = \infty.
\end{equation*}
Since $\varphi$ satisfies $\bC_1$ at $x_i$,
it does not hold that $\varphi(x_i) = \eta_i$,
and so 
there exists $\epsilon > 0$ and
a sequence of points $\{z_\alpha\}_{\alpha = 1}^\infty$ tending
to $x_i$ from $\C^+$ such that
\begin{equation}
|\tilde{c}_i - \tilde{e}_i\varphi(z_\alpha)| \geq \epsilon
\quad\text{for every}\quad
\alpha.
\label{e4.32}
\end{equation}
Note that ${\bf e}_i^*P{\bf e}_i =
\gamma_i$.
When $\varphi$ is subject to $(a)$, we
obtain from (\ref{e4.25}) the following asymptotic
equality:
\begin{equation}
|z - x_i|^2\Psi(z)^*P\Psi(z) = 
|\tilde{c}_i - \tilde{e}_i\varphi(z)|^2\gamma_i
+ O(|z-x_i|).
\label{asy1psi}
\end{equation}
Due to (\ref{e4.27}),
\begin{equation}
|z - x_i|^2|V_\varphi(z)|^2 = |\tilde{c}_i - \tilde{e}_i\varphi(z)|^2
+ O(|z - x_i|)
\label{asy1v}
\end{equation}
holds.  By assumption,
$\varphi(z)$ is bounded near $x_i$, and therefore
\begin{equation}
\lim_{z_\alpha \to x_i}|z - x_i|^2 K_\varphi(z,z) = 0
\label{asy1k}
\end{equation}
by Theorem \ref{CJN2}.
When $\varphi$ is subject to $(b)$, we obtain
from (\ref{e4.250})
the following asymptotic
estimate:
\begin{equation}
|z - x_i|^2\Psi_\varphi(z)^*P\Psi_\varphi(z) = 
|\tilde{c}_i - \tilde{e}_i\varphi(z)|^2\gamma_i
+ o(1).
\label{asy2psi}
\end{equation}
From (\ref{e4.270}) we get
\begin{equation}
|z - x_i|^2|V_\varphi(z)|^2 = |\tilde{c}_i - \tilde{e}_i\varphi(z)|^2
+ o(1).
\label{asy2v}
\end{equation}
Since $\varphi_{-1}(x_i) = 0$, the limit
(\ref{asy1k}) still follows from Theorem
\ref{CJN2}.

Upon substituting either (\ref{asy1psi}),
(\ref{asy1v}), and (\ref{asy1k})
or (\ref{asy2psi}), (\ref{asy2v}), and (\ref{asy1k})
into (\ref{r6.7}) and letting $z_\alpha \to
x_i$, we obtain, due to (\ref{e4.32}) (i.e., $V_\varphi$ does not vanish on our
 sequence), that
\begin{align*}
\lim_{z = z_\alpha \to x_i}K_\w(z,z)
 &= \lim_{z=z_\alpha \to x_i}\frac{|z - x_i|^2
K_\varphi(z,z) + |z - x_i|^2\Psi_\varphi(z)^*
P\Psi_\varphi(z)}{|z - x_i|^2|V_\varphi(z)|^2}\\
&= \frac{0 + \gamma_i}{1} = \gamma_i.
\end{align*}
In view of Theorem \ref{CJGN1}, we conclude
that the nontangential limit ${\displaystyle \lim_{z \to x}K_\w(z,z)}$ 
exists and equals $\gamma_i$, whence $\w^\prime(x_i)$
exists and equals $\gamma_i$ also.  Theorem \ref{CJGN1}
also implies that the limit $\w(x_i)$ exists.
Since, by (\ref{e4.26}) and (\ref{e4.27}),
\begin{equation*}
\lim_{z=z_\alpha \to x_i}\w(z) = 
\lim_{z=z_\alpha \to x_i}\frac{(z-x_i)U_\varphi(z)}{
(z-x_i)V_\varphi(z)} = \w_i,
\end{equation*}
the nontangential limit $\w(x_i)$ equals $\w_i$.

Now let us assume that $\varphi \in \N_0$ is such that
$\varphi_{-1}(x_i) < 0$.  Due to (\ref{e4.68}) and
(\ref{e4.69}),
\begin{equation*}
\w(x_i) := \lim_{z \to x_i}\w(z) =
\lim_{z \to x_i}\frac{(z - x_i)^2U_\varphi(z)}{
(z - x_i)^2V_\varphi(z)} = \frac{\w_i\varphi_{-1}(x_i)
\tilde{e}_i}{\varphi_{-1}(x_i)\tilde{e}_i}.
\end{equation*}
Note that, by Theorem \ref{CJN2}, we have
\begin{equation}
\lim_{z \to x_i}|z - x_i|^4K_\varphi(z,z) = 0.
\label{babel}
\end{equation}
Using (\ref{babel}), (\ref{bull1}), and (\ref{e4.69}),
we may pass to limits in (\ref{r6.7}):
\begin{align*}
\lim_{z \to x_i}K_\w(z,z) &=
\lim_{z \to x_i}\frac{|z - x_i|^4K_\varphi(z,z) +
|z - x_i|^4\Psi_\varphi(z)^*P\Psi_\varphi(z)}{
|z - x_i|^4|V_\varphi(z)|^2} \nonumber \\
& =
\frac{\gamma_i\tilde{e}_i^2\varphi_{-1}^2(x_i)}{
\tilde{e}_i^2\varphi_{-1}^2(x_i)} = \gamma_i.
\end{align*}
By Theorem \ref{CJGN1}, we have
$\w^\prime(x_i) = \gamma_i$ as well.

Finally, we consider the case where $\varphi \equiv \infty$.
We refer to (\ref{r6.7c}) and (\ref{r6.7b}) for the expressions for
$\Psi$, $U_\infty$, and $V_\infty$.
These readily yield the following asymptotic formulas as
$z$ tends to $x_i$ nontagentially:
\begin{equation*}
(z - x_i)\Psi(z) = -{\bf e}_i\tilde{e}_i^* + O(|z - x_i|),
\end{equation*}
\begin{equation*}
(z - x_i)U_\infty(z) = \w_i\tilde{e}_i^* + O(|z - x_i|),
\end{equation*}
\begin{equation*}
(z - x_i)V_\infty(z) = \tilde{e}_i^* + O(|z - x_i|).
\end{equation*}
Observe that $\w(z)$ is rational, that $V_\infty(z) \neq 0$
near $x_i$, and that
\begin{equation*}
\w(z) = \frac{\Theta_{11}(z)}{\Theta_{21}(z)} = \frac{U_\infty(z)}{
V_\infty(z)}.
\end{equation*}
We therefore may pass to limits, obtaining
\begin{equation*}
\lim_{z \to x_i}\w(z) = \lim_{z \to x_i} 
\frac{(z - x_i)U_\infty(z)}{(z - x_i)V_\infty(z)} =
\w_i.
\end{equation*}
Due to (\ref{r6.7}),
\begin{equation*}
\lim_{z \to x_i}K_\w(z,z) =
\lim_{z \to x_i}\frac{|z - x_i|^2\Psi_\varphi(z)^*P\Psi_\varphi(z)}{|z - x_i|^2
|V_\infty(z)|^2} = \gamma_i.
\end{equation*}
By Theorem \ref{CJGN1}, $\w^\prime(x_i) = \gamma_i$ as well.
\end{proof}

\begin{Tm}
Let $i \in \{1, \ldots, \ell\}$ and suppose
that $\tilde{e}_i \neq 0$ (i.e., $\eta_i \neq \infty$).
Let $\varphi \in \tN_0$ 
and let $\w := \mathbf{T}_\Theta[\varphi]$.  If $\varphi$
satisfies condition $\bC_2$ at $x_i$, then
\begin{equation*}
\w(x_i) = \w_i
\quad\text{and}\quad
\w^\prime(x_i) = \gamma_i.
\end{equation*}
\end{Tm}
\begin{proof}
Let us set
\begin{equation*}
\Delta_i(z) := \frac{\tilde{c}_i - \tilde{e}_i\varphi(z)}{
z - x_i}
\end{equation*}
for short and note that
\begin{equation*}
\Delta_i(z) \neq 0
\quad (z \in \C^+).
\end{equation*}
For, if $\varphi(z_0) = \eta_i \in \mathbb{R}$ at some
point $z_0 \in \C^+$, then $\varphi \equiv \eta_i$, which
would contradict the fact that $\varphi$ satisfies
$\bC_2$.   Indeed, since $\varphi$ does satisfy $\bC_2$ at $x_i$, 
Lemma \ref{slemma} is applicable with $\varphi_0 = \eta_i$,
and shows that
\begin{equation}
\lim_{z \to x_0}K_\varphi(z,z) \cdot \frac{1}{|\Delta_i(z)|^2}
= 0
\label{d1again}
\end{equation}
and
\begin{equation}
\lim_{z \to x_0}\Delta_i(z)^{-1} = 0.
\label{d2again}
\end{equation}
Since we assume that $\varphi$ satisfies $\bC_2$,
asymptotic relations (\ref{e4.25}) -- (\ref{e4.27})
apply.  Dividing these by $(\tilde{c}_i -
\tilde{e}_i\varphi(z))$ and writing the resulting
equalities in terms of $\Delta_i$ gives
\begin{align*}
\Delta_i(z)^{-1}\Psi_\varphi(z) &=
-\mathbf{e}_i + \Delta_i(z)^{-1} \cdot O(1),\\
\Delta_i(z)^{-1}U_\varphi(z) &=
-\w_i + \Delta_i(z)^{-1}\cdot O(1) \\
\Delta_i(z)^{-1}V_\varphi(z) &=
-1 + \Delta_i(z)^{-1} \cdot O(1).
\end{align*}
By (\ref{d2}), the following limits exist:
\begin{equation*}
\lim_{z \to x_i}\Delta_i(z)^{-1}\Psi_\varphi(z)
= -\mathbf{e}_i,\quad
\lim_{z \to x_i}\Delta_i(z)^{-1}U_\varphi(z)
= -\w_i,\quad
\lim_{z \to x-i}\Delta_i(z)^{-1}V_\varphi(z)
= -1.
\end{equation*}
Using the preceding limits along with (\ref{d1}),
we may pass to limits in (\ref{r6.7}):
\begin{align*}
\lim_{z \to x_i}K_\w(z,z) &=
\frac{|\Delta_i(z)|^{-2}K_\varphi(z,z)+|\Delta_i(z)|^{-2}
\Psi_\varphi(z)^*P\Psi_\varphi(z)}{
|\Delta_i(z)|^{-2}|V_\varphi(z)|^2} \\
&= \frac{0 + \mathbf{e}_i^*P\mathbf{e}_i}{1} = \gamma_i.
\end{align*}
Finally,
\begin{equation*}
\lim_{z \to x_i}\w(z) = \lim_{z \to x_i}
\frac{\Delta_i(z)^{-1}U_\varphi(z)}{
\Delta_i(z)^{-1}V_\varphi(z)} =
\frac{\w_i}{1} = \w_i.
\end{equation*}
which completes the proof.
\end{proof}

\begin{Tm}
Let $i \in \{1, \ldots, \ell\}$ and suppose that
$\tilde{e}_i \neq 0$ (i.e., $\eta_i \neq \infty$).
Let $\varphi \in \tN_0$ 
and let $\w := \mathbf{T}_\Theta[\varphi]$.
\begin{enumerate}
\item If $\varphi$ satisfies $\bC_3$, then
$\w(x_i) = \w_i$ and $-\infty < \w^\prime(x_i) < \gamma_i$.
\item If $\varphi$ satisfies $\bC_4$, then
$\w(x_i) = \w_i$ and $\gamma_i < \w^\prime(x_i) < \infty$.
\end{enumerate}
\end{Tm}
\begin{proof}
Let us assume that $\varphi \in \N_0$ satisfies
$\bC_{3-4}$.  By Theorem \ref{CJGN1}, the following limits
exist:
\begin{equation*}
\varphi^\prime(x_i) :=
\lim_{z \to x_i}\varphi^\prime(z) =
\lim_{z \to x_i}\frac{\varphi(z) - \eta_i}{z - x_i} =
\lim_{z \to x_i}K_\varphi(z,z).
\end{equation*}
By Corollary $7.9$ in \cite{DB},
the following asymptotic equality
holds:
\begin{equation}
\varphi(z) = \eta_i + (z - x_i)\varphi^\prime(x_i) +
o(|z - x_i|)
\label{r73.1}
\end{equation}
for $z \in \C^+$ close enough to $x_i$.
We shall use this asymptotic relation to show that the functions
$\Psi_\varphi$, $U_\varphi$, and $V_\varphi$,
defined (\ref{r6.4}) and (\ref{r6.4a}),
admit the
following nontangential boundary limits at the interpolation
node $x_i$:
\begin{equation}
\Psi_\varphi(x_i) = \frac{1}{\tilde{e}_i}
\left(P^{-1}\mathbf{e}_i - \mathbf{e}_i
(\tilde{p}_{ii} + \tilde{e}_i^2\varphi^\prime(x_i))\right),
\label{r73.2}
\end{equation}
\begin{equation}
U_\varphi(x_i) = \frac{\w_i}{\tilde{e}_i}
(\tilde{p}_{ii} + \tilde{e}_i^2\varphi^\prime(x_i)),
\quad\text{and}\quad
V_\varphi(x_i) = \frac{1}{\tilde{e}_i}
(\tilde{p}_{ii} + \tilde{e}_i^2\varphi^\prime(x_i)).
\label{r73.3}
\end{equation}
Substitute (\ref{r73.1}) into (\ref{r6.4}):
\begin{align*}
\Psi_\varphi(z) &=
-(zI - X)^{-1}\left\{
\left( \eta_i + (z - x_i)\varphi^\prime(x_i)
+ o(|z - x_i|)\right) \tE^* - \tC^* \right\}\\
&= -(zI - X)^{-1}\left(\eta_i \tE^* - \tC^*\right)
- (zI - X)^{-1}(z - x_i)\varphi^\prime(x_i)\tE^*
+ o(1).
\end{align*}
Substituting (\ref{r72.4}) into the latter equality yields
\begin{equation}
\Psi_\varphi(z) =
\frac{1}{\tilde{e}_i} (zI - X)^{-1}(x_iI - X)P^{-1}\mathbf{e}_i
- (zI - X)^{-1}(z - x_i)\varphi^\prime(x_i)\tE^* + o(1).
\label{r73.5}
\end{equation}
Using (\ref{r72lim1}) and (\ref{r72lim2}), we pass
to limits in (\ref{r73.5}) to obtain
\begin{align}
\Psi_\varphi(x_i) &= \frac{1}{\tilde{e}_i}(I - {\bf e}_i{\bf e}_i^*)
P^{-1}{\bf e}_i - {\bf e}_i{\bf e}_i^*\varphi^\prime(x_i)\tilde{E}^*
\nonumber \\ 
&= \frac{1}{\tilde{e}_i}(P^{-1}{\bf e}_i - {\bf e}_i\tilde{p}_{ii})
- {\bf e}_i\tilde{e}_i\varphi^\prime(x_i) \nonumber \\
&= \frac{1}{\tilde{e}_i}\left(P^{-1}{\bf e}_i - {\bf e}_i(
\tilde{p}_{ii} + \tilde{e}_i^2\varphi^\prime(x_i))\right),
\label{PsiLimit2}
\end{align}
where we have used the fact that
${\bf e}_i^*P^{-1}{\bf e}_i = \tilde{p}_{ii}$ and
${\bf e}_i^*\tilde{E}^* = \tilde{e}_i$.
This proves (\ref{r73.2}).
Using (\ref{r73.2}), we may pass to limits in (\ref{r6.4a}).
We get
\begin{align*}
U_\varphi(x_i) &= \varphi(x_i) - C\Psi_\varphi(x_i) \nonumber \\
&= \eta_i - C \cdot \frac{1}{\tilde{e}_i}\left(P^{-1}{\bf e}_i
-{\bf e}_i(\tilde{p}_{ii} + \tilde{e}_i^2\varphi^\prime(x_i))\right) \\
&= \eta_i - \eta_i + \frac{\w_i}{\tilde{e}_i}(\tilde{p}_{ii} +
\tilde{e}_i^2\varphi^\prime(x_i)) \nonumber \\
&= \frac{\w_i}{\tilde{e}_i}(\tilde{p}_{ii} +
\tilde{e}_i^2\varphi^\prime(x_i))
\end{align*}
and, similarly,
\begin{align*}
V_\varphi(x_i) &= 1 - E\Psi_\varphi(x_i) \nonumber \\
&= 1 - E  \cdot \frac{1}{\tilde{e}_i}\left(P^{-1}{\bf e}_i
-{\bf e}_i(\tilde{p}_{ii} + \tilde{e}_i^2\varphi^\prime(x_i))\right) \\
&= 1 - 1 + \frac{1}{\tilde{e}_i}(\tilde{p}_{ii} + \tilde{e}_i^2
\varphi^\prime(x_i)) = \frac{1}{\tilde{e}_i}(\tilde{p}_{ii} + \tilde{e}_i^2
\varphi^\prime(x_i)),
\end{align*}
thus proving (\ref{r73.3}).

Because $\varphi$ satisfies condition $\bC_{3-4}$ at $x_i$, it follows
that $\tilde{p}_{ii} + \tilde{e}_i^2 \varphi^\prime(x_i) \neq 0$.
Hence (\ref{r6.3}) and (\ref{r73.3}) imply that
\begin{equation*}
\w(x_i) := \lim_{z \to x_i} = \lim_{z \to x_i}
\frac{U_\varphi(z)}{V_\varphi(z)} =
\frac{U_\varphi(x_i)}{V_\varphi(x_i)} = \w_i.
\end{equation*}
Since the limits $\varphi(x_i)$ and $\varphi^\prime(x_i)$ both exist
and are real, it follows from Theorem \ref{CJGN1} that
$K_\varphi(x_i,x_i) = \varphi^\prime(x_i)$.  This combined with 
(\ref{r73.2}) and (\ref{r73.3}) shows that
we may pass to the limit in (\ref{r6.7}):
\begin{equation*}
K_\w(x_i,x_i) := \lim_{z \to x_i}K_\w(z,z) =
\frac{\varphi^\prime(x_i) + \Psi_\varphi(x_i)^*P
\Psi_\varphi(x_i)}{|V_\varphi(x_i)|^2}.
\end{equation*}
Substituting in the explicit formulas yields
\begin{equation*}
K_\w(x_i,x_i)
=
\frac{\tilde{e}_i^2 \varphi^\prime(x_i) +
({\bf e}_i^*P^{-1} - (\tilde{p}_{ii} + \tilde{e}_i^2
\varphi^\prime(x_i)){\bf e}_i^*)
P
(P^{-1}{\bf e}_i - {\bf e}_i(\tilde{p}_{ii} + \tilde{e}_i^2
\varphi^\prime(x_i))) }{
(\tilde{p}_{ii}+\tilde{e}_i^2 \varphi^\prime(x_i))^2},
\end{equation*}
which, after simplifications, becomes
\begin{equation}
K_\w(x_i,x_i) =  \gamma_i - \frac{1}{\tilde{p}_{ii} + \tilde{e}_i^2
\varphi^\prime(x_i)}.
\label{7000}
\end{equation}
By Theorem \ref{CJGN1} and (\ref{7000}), we have that
\begin{equation*}
\w^\prime(x_i) = \lim_{z \to x_i}\w^\prime(z) =
\gamma_i - \frac{1}{\tilde{p}_{ii} + \tilde{e}_i^2
\varphi^\prime(x_i)}
\end{equation*}
Thus if $\varphi$ meets condition $\bC_3$ at $x_i$, then
$\w^\prime(x_i) < \gamma_i$, and if $\varphi$ meets condition
$\bC_4$ at $x_i$, then $\w^\prime(x_i) > \gamma_i$.
\end{proof}

\begin{Tm}
Let $\Omega$ be such that $\ell = n = 1$:
\begin{equation*}
\Omega = \left\{ \{ x_1 \}, \{ \w_1 \}, \{\gamma_1\} \right\}.
\end{equation*}
Let $P$ be the Pick matrix associated to $\Omega$
and suppose that $P$ is positive definite.  Then
all solutions of Problem \ref{Prob1} are parameterized
by $\mathbf{T}_\Theta[\varphi]$, where $\varphi \in \tN_0$
and satisfies $\bC_{1-2}$ at each interpolation node.
\label{ellis1}
\end{Tm}
\begin{proof}
Since $P > 0$ and $E = 1$, it follows that $\tilde{e}_1 \neq 0$.
Thus we need only examine under what conditions parameters
$\varphi$ satisfy $\bC_1$--$\bC_6$ (since conditions
$\btC_1$--$\btC_6$ only apply when $\tilde{e}_i = 0$).  Note
that positive definiteness of $P$ (and hence that of $P^{-1}$)
precludes $\varphi$ from satisfying $\bC_{4-6}$.
Observe that Theorem \ref{Main4} has been proven,
the relevant parts of Theorem \ref{Main3} dealing with
conditions $\bC_1$--$\bC_3$ has already been proven, and the
necessity part of Theorem \ref{lessprecise} has been proven.
Note that, by Remark \ref{Remark6.1}, a function of
the form $\mathbf{T}_\Theta[\varphi]$ for $P > 0$ belongs
to the class $\N_0$.  Therefore, due to Theorems
\ref{lessprecise}, \ref{Main3}, and \ref{Main4},
a function $\w = \mathbf{T}_\Theta[\varphi]$ is a solution
of Problem \ref{Prob1} (with $\ell = n = 1$ and $P>0$)
if and only if $\varphi \in \tN_0$ satisfies
$\bC_{1-2}$ at each interpolation node.
\end{proof}

\begin{Tm}
Let $i \in \{1, \ldots, \ell\}$ and suppose that
$\tilde{e}_i \neq 0$.  Suppose that $\varphi \in \tN_0$
satisfies condition $\bC_5$ at $x_i$.
The function $\w := \mathbf{T}_\Theta[\varphi]$ is subject to
one of the following:
\begin{enumerate}
\item \text{The nontangential limit $\w(x_i)$ does not exist.}
\item \text{The nontangential limit $\w(x_i)$ exists and $\w(x_i) \neq
  \w_i$.}
\item \text{The nontangential limit $\w(x_i)$ exists and equals
  $\w_i$, and ${\displaystyle \left|\lim_{z \to x_i}K_\w(z,z)\right| = \infty.}$}
\end{enumerate}
\label{bC5}
\end{Tm}
\begin{proof}
The conditions imposed by
condition $\bC_5$
form a well-posed one-point interpolation problem; in particular, they
form a problem of type Problem \ref{Prob1} with one
interpolation node and with $\ell = 1$.  By Theorem
\ref{ellis1},
any solution $\varphi$ of this problem
admits a representation
\begin{equation*}
\varphi = {\bf T}_{\hat\Theta} [ \hat{\varphi} ]
\end{equation*}
with the coefficient matrix $\hat{\Theta}$ defined via formula
(\ref{ThetaFormula}), but with $P, X, E$ and $C$ replaced by
$-\frac{\tilde{p}_{ii}}{\tilde{e}_i^2}, x_i, 1$ and $\eta_i$,
respectively:
\begin{equation}
\hat{\Theta}(z) = I_2 + \frac{i}{z - x_i}\begin{bmatrix}\eta_i \\ 1
\end{bmatrix}\frac{\tilde{e}_i^2}{\tilde{p}_{ii}}\begin{bmatrix}
\eta_i & 1\end{bmatrix}J,
\label{e4.114}
\end{equation}
and where $\hat{\varphi} \in \tN_0$ and is subject to
$\bC_{1-2}$ at $x_i$.
Since $\eta_i = \frac{\tilde{c}_i}{\tilde{e}_i}$ and $\tilde{e}_i \neq
0$, (\ref{e4.114}) can be written as
\begin{equation*}
\hat{\Theta}(z) = I_2 + \frac{i}{z - x_i}\begin{bmatrix}\tilde{c}_i
\\ \tilde{e}_i\end{bmatrix}\frac{1}{\tilde{p}_{ii}}\begin{bmatrix}
\tilde{c}_i & \tilde{e}_i\end{bmatrix}J.
\end{equation*}
By the reflection principle, $\Theta(z)^{-1} = J\Theta(\bar{z})^*J$,
and so,
\begin{equation}
\hat{\Theta}(z)^{-1} = I_2 - \frac{i}{z - x_i}\begin{bmatrix}
\tilde{c}_i \\ \tilde{e}_i\end{bmatrix}\frac{1}{\tilde{p}_{ii}}
\begin{bmatrix}\tilde{c}_i & \tilde{e}_i\end{bmatrix}J.
\label{4.120}
\end{equation}
Note that $\hat{\Theta}(z)^{-1}$ coincides with the function $\hat{\Theta}^{(2)}$ introduced
in (\ref{e3.24}) (if we rearrange
interpolation nodes so that $x_i \mapsto x_n$, which may be done
without loss of generality) where we take $\tC_2 = \eta_i$, 
$\tE_2 = 1$, and $\widetilde{P}_{22}^{-1} = -\tilde{e}_i^2/\tilde{p}_{ii}$.  
Therefore, by (\ref{4.120}) and (\ref{e3.25}),
\begin{equation*}
\Theta(z) = \Theta^{(1)}(z)\hat{\Theta}(z)^{-1}
\end{equation*}
where $\Theta^{(1)}$ is given in (\ref{e3.23}).  
Thus the following chain of equalities holds:
\begin{equation*}
w:= {\bf T}_\Theta[\varphi] = {\bf T}_\Theta[{\bf T}_{\hat\Theta}
[ \hat{\varphi}]]= {\bf T}_{\Theta \hat{\Theta}}[ \hat{\varphi}] =
{\bf T}_{\Theta^{(1)}}[\hat{\varphi}].
\end{equation*}
Upon setting
\begin{equation*}
U_{\hat{\varphi}}(z) = \Theta_{11}^{(1)}(z)\hat{\varphi}(z) +
\Theta_{12}^{(1)}(z)
\quad\text{and}\quad
V_{\hat{\varphi}}(z) = \Theta_{21}^{(1)}(z)\hat{\varphi}(z) +
\Theta_{22}^{(1)}(z)
\end{equation*}
so that
\begin{equation*}
\begin{bmatrix}U_{\hat{\varphi}} \\ V_{\hat{\varphi}}\end{bmatrix} =
\Theta^{(1)}\begin{bmatrix}\hat{\varphi} \\ 1\end{bmatrix}
\end{equation*}
we arrive at
\begin{equation}
\w(z) = {\bf T}_{\Theta^{(1)}}[\hat{\varphi}](z) =
\frac{\Theta_{11}^{(1)}(z)\hat{\varphi}(z) + \Theta_{12}^{(1)}(z)}{
\Theta_{21}^{(1)}(z)\hat{\varphi}(z) + \Theta_{22}^{(1)}(z)} =
\frac{U_{\hat{\varphi}}(z)}{V_{\hat{\varphi}}(z)}.
\label{4.125}
\end{equation}
Note that $\hat{\varphi}$ satisfies precisely one of the following:
\begin{align*}
&(a)\;\;\text{The limit $\hat{\varphi}(x_i)$ does not exist.}\\
&(b)\;\;\text{The limit $\hat{\varphi}(x_i)$ exists, is finite,
and not equal to $\eta_i$.}\\
&(c)\;\;\text{It holds that $\hat{\varphi}(x_i) = \infty$.}\\
&(d)\;\;\text{It holds that $\hat{\varphi}(x_i) = \eta_i$ and
$K_{\hat{\varphi}}(x_i,x_i) = \infty$.}
\end{align*}
The function $\Theta^{(1)}$ is a rational function invertible
at $x_i$, and so
the limit $\w(x_i)$ exists if and only if the limit
$\hat{\varphi}(x_i)$ exists, i.e., $(1) \iff (a)$.

Let us henceforth assume that the limit $\hat{\varphi}(x_i)$
exists.  Now
suppose that $V_{\hat{\varphi}}(x_i) = 0$.
If $U_{\hat{\varphi}}(x_i) \neq 0$, then $\w(x_i) = \infty \neq \w_i$.
On the other hand, if $U_{\hat{\varphi}}(x_i) = 0$, then
\begin{equation*}
\Theta^{(1)}(x_i)\begin{bmatrix}\hat{\varphi}(x_i) \\ 1\end{bmatrix} =
\begin{bmatrix}U_{\hat{\varphi}}(x_i) \\
 V_{\hat{\varphi}}(x_i)\end{bmatrix} = 0, 
\end{equation*}
which implies that $\Theta^{(1)}(x_i)$ is singular, 
which is a contradiction.
We henceforth assume that $V_{\hat{\varphi}}(x_i) \neq 0$,
and split the rest of the proof into two complementary
cases:

\nline
{\bf Case I:}
Suppose that $\hat{\varphi}(x_i) \neq \infty$, i.e.,
that $\hat{\varphi}$ is subject to either $(b)$ or $(d)$.
Observe that
\begin{align}
\w(x_i) - \w_i &= \frac{U_{\hat{\varphi}}(x_i) - \w_i
V_{\hat{\varphi}}(x_i)}{V_{\hat{\varphi}}(x_i)} \nonumber \\
&= \frac{1}{V_{\hat{\varphi}}(x_i)}\begin{bmatrix}1 & -\w_i
\end{bmatrix}\Theta^{(1)}(x_i)\begin{bmatrix}
\hat{\varphi}(z) \\ 1\end{bmatrix}.
\label{sec6abo}
\end{align}
We shall prove that
\begin{equation}
\begin{bmatrix}1 & -\w_i\end{bmatrix} = \frac{1}{\tilde{p}_{ii}}
\begin{bmatrix}\tilde{e}_i & -\tilde{c}_i\end{bmatrix}
\Theta^{(1)}(x_i)^{-1},
\label{e4.129}
\end{equation}
which, combined with (\ref{sec6abo}), yields
\begin{equation}
\w(x_i) - \w_i = \frac{\tilde{e}_i}{\tilde{p}_{ii}
V_{\hat{\varphi}}(x_i)}\left(\hat{\varphi}(x_i) - \eta_i\right),
\label{yield1}
\end{equation}
which in turn shows that  $\w(x_i) \neq \w_i$ 
provided the $\hat{\varphi}(x_i) \neq \eta_i$.
To prove (\ref{e4.129}), note that, by Lemma
\ref{L42},
\begin{equation*}
{\rm Res}_{z = x_i} \Theta(z) = 
-\begin{bmatrix}\w_i \\ 1\end{bmatrix}
\begin{bmatrix}-\tilde{e}_i & \tilde{c}_i\end{bmatrix}
\end{equation*}
and
\begin{equation*}
{\rm Res}_{z = x_i} \hat{\Theta}(z)^{-1} = 
\frac{-1}{\tilde{p}_{ii}}
\begin{bmatrix}
\tilde{c}_i \\ \tilde{e}_i\end{bmatrix}
\begin{bmatrix}
-\tilde{e}_i & \tilde{c}_i\end{bmatrix}.
\end{equation*}
Taking into account the fact that $\Theta^{(1)}$ is analytic
at $x_i$, we compare residues on boths sides of
$\Theta(z) = \Theta^{(1)}(z)\hat{\Theta}(z)$ at $x_i$:
\begin{equation*}
-\begin{bmatrix}\w_i \\ 1\end{bmatrix}
\begin{bmatrix}
-\tilde{e}_i & \tilde{c}_i\end{bmatrix} = 
\frac{-1}{\tilde{p}_{ii}}\Theta^{(1)}(x_i)
\begin{bmatrix}
\tilde{c}_i \\ \tilde{e}_i\end{bmatrix}
\begin{bmatrix}
-\tilde{e}_i & \tilde{c}_i\end{bmatrix}.
\end{equation*}
Since $\tilde{e}_i \neq 0$, it follows that
\begin{equation*}
\begin{bmatrix}\w_i \\ 1\end{bmatrix} =
\frac{1}{\tilde{p}_{ii}}\Theta^{(1)}(x_i)
\begin{bmatrix}\tilde{c}_i \\ \tilde{e}_i\end{bmatrix}.
\end{equation*}
Therefore
\begin{equation*}
\begin{bmatrix}1 \\ -\w_i\end{bmatrix} =
iJ\begin{bmatrix}\w_i \\ 1\end{bmatrix}
= \frac{1}{\tilde{p}_{ii}}J\Theta^{(1)}(x_i)J \cdot iJ
\begin{bmatrix}\tilde{c}_i \\ \tilde{e}_i\end{bmatrix} =
\frac{1}{\tilde{p}_{ii}}\Theta^{(1)}(x_i)^{-1}
\begin{bmatrix}\tilde{e}_i \\ -\tilde{c}_i\end{bmatrix}.
\end{equation*}
Taking adjoints in the above
equalities gives (\ref{e4.129}).

Due to (\ref{yield1}), it follows that
$(b) \implies (2)$, and that $(d)$ implies
that $\w(x_i) = \w_i$.  Suppose that
$\hat{\varphi}$ is subject to $(d)$.
Since $\w = \mathbf{T}_{\Theta^{(1)}}[\hat{\varphi}]$,
we may write (\ref{r6.7}) in terms of $\hat{\varphi}$:
\begin{equation}
K_\w(z,z) = \frac{1}{|V_{\hat{\varphi}}(z)|^2}
(K_{\hat{\varphi}}(z,z) +
\hat{\Psi}_{\hat{\varphi}}(z)^*P\hat{\Psi}_{\hat{\varphi}}(z)),
\label{Portugal}
\end{equation}
where 
\begin{equation*}
\hat{\Psi}(z) = (zI - X_1)^{-1}(C_1^* - \hat{\varphi}(z)E_1^*).
\end{equation*}
The latter equality shows that the limit $\hat{\Psi}(x_i)$
exists and is finite.
By assumption, $0 < |V_{\hat{\varphi}}(z)| < \infty$
and $K_{\hat{\varphi}}(x_i,x_i) = \infty$.
Hence, taking limits in (\ref{Portugal}) shows that
\begin{equation*}
\lim_{z \to x_i}K_w(z,z) = \infty.
\end{equation*}
Therefore $(d) \implies (3)$.

\nline
{\bf Case II:}
Suppose that $\hat{\varphi}(x_i) = \infty$, i.e.,
that $\hat{\varphi}$ is subject to $(c)$.
Due to (\ref{4.125}), we have
\begin{equation}
\w(x_i) = \frac{\Theta_{11}^{(1)}(x_i)}{\Theta_{21}^{(1)}(x_i)}.
\label{SouthAfrica}
\end{equation}
Note that
\begin{equation*}
\Theta^{(1)}(x_i)\begin{bmatrix}1 \\ 0\end{bmatrix} =
\begin{bmatrix}\Theta_{11}^{(1)}(x_i) \\
\Theta_{21}^{(1)}(x_i)\end{bmatrix}.
\end{equation*}
The preceding equality and the fact that $\Theta^{(1)}(x_i)$ is
invertible imply that
we cannot have $\Theta_{11}^{(1)} = \Theta_{21}^{(1)}(x_i) = 0.$
If $\Theta_{21}^{(1)}(x_i) = 0$, then $|\w(x_i)| = \infty$ and
so $\w(x_i) \neq \w_i$.  Let us henceforth assume that
$\Theta_{21}^{(1)}(x_i) \neq 0$.
By (\ref{SouthAfrica})
\begin{align*}
\w(x_i) - \w_i &=
\frac{\Theta_{11}^{(1)}(x_i) - \w_i\Theta_{21}^{(1)}(x_i)}{
\Theta_{21}^{(1)}(x_i)} \\
&= \frac{1}{\Theta_{21}^{(1)}(x_i)}\begin{bmatrix}1 &
  -\w_i\end{bmatrix}
\Theta^{(1)}(x_i)\begin{bmatrix}1 \\ 0\end{bmatrix}.
\end{align*}
We may then use equality (\ref{e4.129}) to show that
\begin{equation*}
\w(x_i) - \w_i =
\frac{\tilde{e}_i}{\tilde{p}_{ii}\Theta_{21}^{(1)}(x_i)} \neq 0.
\end{equation*}
Hence $(c) \implies (2)$.

We have shown that $(1) \iff (a)$, $(b) \vee (c) \implies
(2),$ and $(d) \implies (3)$.
Since $\w = \mathbf{T}_{\Theta^{(1)}}[\hat{\varphi}]$ and
since
$\hat{\varphi}$ satisfies exactly one of $(a)$, $(b)$, $(c)$,
or $(d)$,
we have $(1) \iff (a)$, $(2) \iff (b) \vee (c)$, and $(3) \iff (d)$,
which completes the proof.
\end{proof}

\begin{Tm}
Let $i \in \{1, \ldots, \ell\}$ and suppose that
$\tilde{e}_i \neq 0$ (i.e., $\eta_i \neq \infty$).
Let $\varphi \in \tN_0$ 
and let $\w := \mathbf{T}_\Theta[\varphi]$.  If $\varphi$
satisfies condition $\bC_6$ at $x_i$, then the limit
$\w(x_i)$ exists but does not equal $\w_i$.
\end{Tm}
\begin{proof}
Since $\varphi(x_i) = \eta_i$ and $\varphi^\prime(x_i) = 0$ 
by $\mathbf{C}_6$, we in fact have
\begin{equation}
\varphi \equiv \eta_i.
\label{China}
\end{equation}
Without loss of generality, we may
rearrange interpolation nodes so that $x_i \mapsto x_n$.  
Then, using decompositions
(\ref{PP}) and (\ref{XEC}) with $\tilde{p}_{ii} = 0$, we obtain
\begin{equation*}
\tilde{P}_{21}P_{12} = 1 \text{ and }
P^{-1}{\bf e}_i =
\begin{bmatrix}\tilde{P}_{12} \\ 0\end{bmatrix}
\end{equation*}
The identity
\begin{equation*}
P_{12} = (x_i- X_1)^{-1}(\w_i E _1^* - C_1^*),
\end{equation*}
follows from the Lyapunov identity (\ref{LI}) upon substituting
the partitionings (\ref{PP}) and (\ref{XEC})
into equation (\ref{LI}) and comparing the
$(1,2)$ block entries.
Multiplying both sides of the identity (\ref{e3.4}) by
$\mathbf{e}_i$ on the right and rearranging terms yields
\begin{equation*}
\tilde{e}_i\tC^* - \tilde{c}_i\tE^* = -(x_iI - X)P^{-1}\mathbf{e}_i.
\end{equation*}
We may divide the latter equality by $\tilde{e}_i \neq 0$ to get
\begin{equation}
\tC^* - \eta_i\tE^* = -\frac{1}{\tilde{e}_i}(x_iI - X)P^{-1}\mathbf{e}_i.
\label{Japan}
\end{equation}
Using (\ref{China})
and (\ref{Japan}) in (\ref{r6.4}) shows that
\begin{equation}
\Psi_\varphi(z)
= -\frac{1}{\tilde{e}_i}(zI - X)^{-1}(x_iI - X)P^{-1}
{\bf e}_i.
\label{e4.140}
\end{equation}
Substituting (\ref{e4.140}) into (\ref{r6.4a}), we obtain
\begin{align*}
U_\varphi(z) &= \eta_i - \frac{1}{\tilde{e}_i}C (zI - X)^{-1}
(x_iI - X)P^{-1}{\bf e}_i \\
&= \frac{1}{\tilde{e}_i}CP^{-1}{\bf e}_i - \frac{1}{\tilde{e}_i}
C(zI - X)^{-1}(x_iI - X)P^{-1}{\bf e}_i \\
&= \frac{1}{\tilde{e}_i}C\left( I - (zI - X)^{-1}(x_iI - X)
\right)P^{-1}{\bf e}_i \\
&= \frac{1}{\tilde{e}_i} C (zI - X)^{-1}\left((zI - X) -
(x_iI - X)\right)P^{-1}{\bf e}_i\\
&= \frac{z - x_i}{\tilde{e}_i}C(zI - X)^{-1}P^{-1}{\bf e}_i,
\end{align*}
and
\begin{align*}
V_\varphi(z) &= 1 - \frac{1}{\tilde{e}_i}E (zI - X)^{-1}
(x_iI - X)P^{-1}{\bf e}_i \\
&= \frac{1}{\tilde{e}_i}E\left(I - (zI - X)^{-1}(x_iI - X)
\right)P^{-1}{\bf e}_i \\
&=\frac{z - x_i}{\tilde{e}_i}E(zI - X)^{-1}P^{-1}{\bf e}_i.
\end{align*}
In terms of the partitions (\ref{PP}) and (\ref{XEC}), the 
preceding equalities become
\begin{equation*}
U_\varphi(z) = \frac{z - x_i}{\tilde{e}_i}C_1 (zI - X_1)^{-1}
\tilde{P}_{12}
\quad\text{and}\quad
V_\varphi(z) = \frac{z - x_i}{\tilde{e}_i}E_1 (zI - X_1)^{-1}
\tilde{P}_{12}.
\end{equation*}
The latter equalities thus show that
\begin{equation*}
\w(z) = \frac{U_\varphi(z)}{V_\varphi(z)} = \frac{C_1(
zI - X_1)^{-1}\tilde{P}_{12}}{E_1(zI - X_1)^{-1}\tilde{P}_{12}}.
\end{equation*}
Since $\w$ is rational, the limit ${\displaystyle \lim_{z \to
    x_i}|\w(x_i)|}$
exists.  
We shall show
that when this limit is finite, it holds that
\begin{align*}
\w_i - \w(x_i) &= \w_i - \frac{C_1(x_iI - X_1)^{-1}
\tilde{P}_{12}}{E_1(x_iI - X_1)^{-1}\tilde{P}_{12}} \\
&= \frac{(\w_i E_1 - C_1)(x_iI - X_1)^{-1}\tilde{P}_{12}}{
E_1(x_iI - X_1)^{-1}\tilde{P}_{12}} \\
&= \frac{P_{21}\tilde{P}_{12}}{E_1(x_iI - X_1)^{-1}
\tilde{P}_{12}}\\
&= \frac{1}{E_1(x_iI - X_1)^{-1}\tilde{P}_{12}}.
\end{align*}
If $E_1(x_iI - X_1)^{-1}\tilde{P}_{12} = 0$, then
we must have
$C_1(x_iI - X_1)^{-1}\tilde{P}_{12} = 0$ in order for
the limit $\w(x_i)$ to be finite.  However, this forces
$(\w_i E_1 - C_1)(X_1 - x_iI)^{-1}\tilde{P}_{12} = 0$, which is 
a contradiction by the above calculation.  Hence
the limit of $\w$ at $x_i$ does not equal $\w_i$.
\end{proof}

\section{A Special Case}
Let us consider the Problem \ref{Prob1} under the assumption that
$\ka := {sq}_-P = 0$, i.e., the case where $P$ is positive semidefinite.

\begin{Tm}
Let the matrices $P$, $X$, $E$, and $C$ be associated to the data
set $\Omega$ as in (\ref{1.19}) -- (\ref{120}), (\ref{r2.2}), and
(\ref{r2.3}).
If $P$ is positive definite, then a function $\w$ meromorphic on
$\C^+$
is a solution of Problem \ref{Prob1} if and only if it
is of the form
\begin{equation*}
\w(z) = \mathbf{T}_\Theta[\varphi](z) :=
\frac{\Theta_{11}(z)\varphi(z) + \Theta_{12}(z)}{
\Theta_{21}(z)\varphi(z) + \Theta_{22}(z)}
\end{equation*}
for some function $\varphi \in \tN_0$ satisfying
$\bC_{1-3} \vee \btC_{1-3}$ at each interpolation node
$x_i$.
\label{Canada}
\end{Tm}
\begin{proof}
This theorem is a special case of Theorem \ref{lessprecise},
the necessity part of which was proven in Chapter V.
Thus the necessity part of the theorem is proved.
Sufficiency for the
case $\ell = 0$ was established in Theorem \ref{ellis0}.
Sufficiency for the case $n = \ell = 1$ was established in
Theorem \ref{ellis1}.  Since Theorems \ref{Main1}--\ref{Main4}
are now completely proved, Theorems \ref{ellis0} and
\ref{ellis1} are applicable to the more general situation
where we only assume $P > 0$.
\end{proof}

\begin{Tm}
Let the matrices $P$, $X$, $E$, and $C$ be associated to the data
set $\Omega$ as in (\ref{1.19}) -- (\ref{120}), (\ref{r2.2}), and
(\ref{r2.3}).
If $P$ is positive semidefinite and singular, 
then there is a rational solution
of Problem \ref{Prob1} given by
\begin{equation}
\w(z) = \frac{x^*(zI - X)^{-1}C^*}{x^*(zI - X)^{-1}E^*},
\label{degsol}
\end{equation}
where $x$ is any nonzero vector such that $Px = 0$.  Moreover,
this solution is unique.
\label{T7.8}
\end{Tm}
\begin{proof}
By the necessity part in Theorem \ref{lessprecise}, 
proven in Chapter V, for $\w$ to be solution, it is necessary
that the kernel
\begin{equation*}
\mathbf{K}_\w(z,\zeta) :=
\begin{bmatrix}
P & (zI_n - X)^{-1}(\w(z)E^* - C^*) \\
(\w(\zeta)^* E - C)(\bar{\zeta}I_n - X)^{-1} & K_\w(z, \zeta)
\end{bmatrix}
\end{equation*}
be positive $\C^+ \times \C^+$:
\begin{equation}
\mathbf{K}_\w(z, \zeta) \succeq 0.
\label{Kwpos}
\end{equation}
Let $y$ be any nonzero vector such that $Py = 0$, and set
\begin{equation*}
Y = \begin{bmatrix}y & 0 \\ 0 & 1\end{bmatrix}.
\end{equation*}
Then
\begin{equation*}
Y^*\mathbf{K}_\w(z, \zeta)Y
= \begin{bmatrix} 0 & y^*(zI_n - X)^{-1}(\w(z)E^* - C^*) \\
(\w(\zeta)^* E - C)(\bar{\zeta}I_n - X)^{-1}y & K_\w(z, \zeta)
\end{bmatrix},
\end{equation*}
since $y^*Py = 0$.
Positivity condition (\ref{Kwpos}) therefore forces
\begin{equation}
y^*(zI_n - X)^{-1}(\w(z) E^* - C^*) \equiv 0.
\label{pre}
\end{equation}
We may rearrange (\ref{pre})
so that
\begin{equation}
\w(z) = \frac{y^* (zI_n - X)^{-1}C^*}{y^* (zI_n - X)^{-1}E^*}.
\label{wasin}
\end{equation}
We claim that the function $\w$ defined in (\ref{wasin}) is
independent of the choice of $y$ so long as $y$ is nonzero
and satisfies $Py = 0$.  Let $\tilde{\w}(z)$ be defined
by the formula (\ref{wasin}), but with $\tilde{y}$
in place of $y$, where we assume $\tilde{y} \neq 0$ and
$P\tilde{y} = 0$.
By taking adjoints and replacing $\bar{z}$ by $z$
in (\ref{wasin}), we get
a representation of $\w$ equivalent to (\ref{wasin}):
\begin{equation*}
\w(z) = \frac{C(zI_n - X)^{-1}y}{E(zI_n - X)^{-1}y}.
\end{equation*}
The numerator of the difference $\w(z) - \tilde{w}(z)$
is then given by
\begin{equation}
\tilde{y}^*(zI_n - X)^{-1}(E^*C - C^*E)(zI_n - X)^{-1}y
\label{numerator}
\end{equation}
Due to the Lyaponov identity (\ref{LI}) 
we may rewrite (\ref{numerator}), which becomes
\begin{align*}
&\tilde{y}^*(zI_n - X)^{-1}(PX - XP)(zI_n - X)^{-1}y \\
=&\; \tilde{y}^*(zI_n - X)^{-1}\left[
(zI_n - X)P - P(zI_n - X)\right](zI_n - X)^{-1}y \\
=&\; \tilde{y}^*P(zI_n - X)^{-1}y - 
\tilde{y}^*(zI_n - X)^{-1}Py \\
=&\; 0.
\end{align*}
The last equality follows since $Py = 0$ and
$(\tilde{y}^*P)^* = P\tilde{y} = 0$.
Therefore $\w(z) \equiv \tilde{w}(z)$.

We shall show in Chapter IX
that the unique $\w$, given in (\ref{wasin}), is in fact
a solution to Problem \ref{Prob1}.
\end{proof}

\chapter{Negative Squares of the Kernel $K_{\mathbf{T}_\Theta[\varphi]}$}
Throughout this section, we assume that the Pick matrix $P$ associated
to the data set $\Omega$, given in (\ref{1.18}), is invertible.
\begin{Tm}
There exists a function $\varphi \in \tN_0$ satisfying
conditions $\bC_{4-6} \vee \tilde{\bC}_{4-6}$
at $x_{i_1}, \ldots, x_{i_k}$ if and only if the $k \times k$
matrix
\begin{equation*}
\mathcal{P} := \left[\tilde{p}_{i_\alpha,i_\beta}\right]_{\alpha,\beta=1}^k
\end{equation*}
is negative semidefinite.  If $\mathcal{P}$ is nonsingular, then there
are infinitely many such functions, and if $\mathcal{P}$ is singular,
then there is a unique such function and moreover this function is rational.
\label{Thm1}
\end{Tm}
\begin{proof}
Rearrange the interpolation nodes $x_i$ so that
$i_\alpha = \alpha$ for $\alpha \in \{1, \ldots, k\}$.
Assume that
$\varphi$ satisfies conditions $\bC_{4-6}$
at interpolation nodes $x_1, \ldots, x_{m}$,
$\tilde{\bC}_{4-6}$ at interpolation nodes
$x_{m + 1}, \ldots, x_{k}$, and
$\bC_{1-3} \vee \tilde{\bC}_{1-3}$ at 
$x_{k + 1}, \ldots,
x_n$.  Thus, we assume that
\begin{align}
\varphi(x_i) &= \eta_i \neq \infty \quad\text{and}\quad
\varphi^\prime(x_i) \leq -\frac{\tilde{p}_{ii}}{\tilde{e}_i^2}
&(i = 1, \ldots, m), \nonumber \\
\varphi(x_j) &= \eta_j = \infty \quad\text{and}\quad
\frac{1}{\varphi_{-1}(x_j)} \geq \frac{\tilde{p}_{jj}}{\tilde{c}_j^2}
&(j = m + 1 , \ldots, k).
\label{conditions}
\end{align}
Let us partition the inverse of $P$ as follows:
\begin{equation*}
P^{-1} = \begin{bmatrix}\tilde{P}_{11} & \tilde{P}_{12} \\
\tilde{P}_{21} & \tilde{P}_{22}\end{bmatrix} \quad\text{with}\quad
\tilde{P}_{11} \in \C^{k \times k}
\end{equation*}
and observe that $\tilde{P}_{11}$ coincides with $\mathcal{P}$.

Conditions (\ref{conditions}) form a well-posed boundary
Nevanlinna-Pick problem in the Nevanlinna class $\N_0$.  
Let us start the proof by first making the following
technical assumption:
\begin{equation}
\tilde{p}_{ii} \neq 0 \quad\quad (i = 1, \ldots, k).
\label{technical}
\end{equation}
Then, this Nevanlinna-Pick problem has
a solution $\varphi$ if and only if the corresponding
Pick matrix $\mathbb{P} = [\mathbb{P}_{ij}]_{i,j = 1}^k$ 
is positive semidefinite:
\begin{equation*}
\mathbb{P} = \begin{bmatrix}
\mathbb{P}^{(11)} & \mathbb{P}^{(12)} \\
\mathbb{P}^{(21)} & \mathbb{P}^{(22)}\end{bmatrix},
\end{equation*}
where
\begin{equation*}
\mathbb{P}^{(11)} = [\mathbb{P}_{ij}]_{i,j=1}^{m}
\qquad
\mathbb{P}_{ij} = 
\left\{ \begin{matrix}
\frac{\eta_j - \eta_i}{x_j - x_i} & (i \neq j) \\
-\frac{\tilde{p}_{ii}}{\tilde{e}_i^2} & (i = j)
\end{matrix}\right.  ,
\end{equation*}

\begin{equation*}
\mathbb{P}^{(22)} = [\mathbb{P}_{ij}]_{i,j= m + 1}^n \qquad
\mathbb{P}_{ij} = 
\left\{ \begin{matrix}
0 & (i \neq j) \\
-\frac{\tilde{c}_i^2}{\tilde{p}_{ii}} & (i = j)
\end{matrix} \right. ,
\end{equation*}
\begin{equation*}
\mathbb{P}^{(12)} = [\mathbb{P}_{ij}]_{\stackrel{i = 1, \ldots, m}{
\mathop{}^{j = m + 1, \ldots, k}}} \qquad
\mathbb{P}_{ij} = \frac{1}{x_j - x_i} \cdot \frac{\tilde{c}_j^2}{\tilde{p}_{jj}},
\end{equation*}
and 
\begin{equation*}
\mathbb{P}^{(21)} = \left(\mathbb{P}^{(12)}\right)^*.
\end{equation*}
There exist infinitely many $\varphi \in \N_0$ satisfying (\ref{conditions}) if 
$\mathbb{P} > 0$ and there exists a unique $\varphi \in \N_0$ satisfying
(\ref{conditions}) if $\mathbb{P} \geq 0$ is singular, by
Theorem \ref{T7.8}.  To complete the
proof under assumption (\ref{technical}), it suffices to show that
\begin{equation}
\mathbb{P} > 0 \iff \tilde{P}_{11} < 0, \quad
\mathbb{P} \geq 0 \iff \tilde{P}_{11} \leq 0, \quad\text{and}\quad
{\rm rank}\;\mathbb{P} = {\rm rank}\;\tilde{P}_{11}.
\label{claims}
\end{equation}
Note that
\begin{align*}
\tilde{e}_i \cdot \mathbb{P}_{ij} \cdot \tilde{e}_j &= -\tilde{p}_{ij} 
& (i,j \in \mathbb{N}_m), \\
\tilde{e}_i \cdot \mathbb{P}_{ij} \cdot
-\frac{\tilde{p}_{jj}}{\tilde{c}_j} &= -\tilde{p}_{ij}
& (i \in \mathbb{N}_m; j \in \mathbb{N}_k \setminus \mathbb{N}_m), \\
-\frac{\tilde{p}_{ii}}{\tilde{c}_i} \cdot \mathbb{P}_{ij} \cdot
-\frac{\tilde{p}_{jj}}{\tilde{c}_j} &= -\tilde{p}_{ij}
& (i,j \in \mathbb{N}_k \setminus \mathbb{N}_m).
\end{align*}
These equalities may be rewritten as
\begin{equation}
\bC^*\mathbb{P}\bC = -\tilde{P}_{11}
\label{CP}
\end{equation}
where $\bC$ is the $k \times k$ diagonal matrix
given by
\begin{equation}
\bC = {\rm diag}\;\left(\tilde{e}_1, \ldots, \tilde{e}_m,
-\frac{\tilde{p}_{m+1,m+1}}{\tilde{e}_{m+1}}, \ldots,
  -\frac{\tilde{p}_{k k}}{\tilde{e}_k}\right).
\label{CP2}
\end{equation}
Since $\bC$ is invertible, all statements in (\ref{claims})
follow from (\ref{CP}).

It now remains to be shown that the technical assumptions 
(\ref{technical}) may be
lifted.  Assume that $\tilde{p}_{ii} = 0$ for some $i \in \mathbb{N}_k
\setminus \mathbb{N}_m$.  
Then, in order for $\varphi_{-1}(x_i)^{-1} \geq
0$ to hold, we must have $\varphi \equiv \infty$.  In this case,
$\varphi$ cannot satisfy $\bC_{4-6}$ at any node.  On the other
hand, if $x_j$ is a regular interpolation node,
then the Hermitian matrix
$\begin{bmatrix}0 & \tilde{p}_{ij} \\ \tilde{p}_{ji} &
\tilde{p}_{jj}\end{bmatrix}$
has real nonzero $\tilde{p}_{ij} = \tilde{p}_{ji}$, and therefore has
eigenvalues of both signs.  Hence, under the assumption that
$\tilde{p}_{ii} = 0$ for some $i \in \mathbb{N}_k \setminus
\mathbb{N}_m$,
we may assume that $m = 0$, i.e.,
$\Omega$ only has singular interpolation nodes.  
In this case, the matrix $\tilde{P}_{11}$
is diagonal and singular, and it is easy to see that each condition in
(\ref{conditions}) is satisfied if and only if each $\tilde{p}_{ii} \leq
0$, i.e., if and only if $\tilde{P}_{11}$ is negative semidefinite.
\begin{Rk}
In particular, if $\tilde{p}_{ii} = 0$ for some $i \in \mathbb{N}_k
\setminus \mathbb{N}_m$, then $\varphi \equiv \infty$. 
\label{Yes}
\end{Rk}
\end{proof}

\begin{La}
Let $P \in \C^{n \times n}$ be an invertible Hermitian matrix and let
\begin{equation*}
P = \begin{bmatrix}P_{11} & P_{12} \\ P_{21} & P_{22}\end{bmatrix}
\quad \text{and} \quad
P^{-1} = \begin{bmatrix}\tilde{P}_{11} & \tilde{P}_{12} \\
\tilde{P}_{21} & \tilde{P}_{22}\end{bmatrix}
\end{equation*}
be conformal block decompositions of $P$ and of $P^{-1}$, respectively, with
$P_{22}, \tilde{P}_{22} \in \C^{k \times k}$.  If
$\tilde{P}_{22}$ is negative semidefinite, then
\begin{equation*}
{\rm sq}_-P_{11} = {\rm sq}_-P - k.
\end{equation*}
\label{Lemma}
\end{La}
\begin{proof}
For the proof, see \cite{BK}.
\end{proof}

\begin{Tm}
If the Pick matrix $P$ is invertible and has $\ka$ negative
eigenvalues, then a Nevanlinna function $\varphi \in \tN_0$
may satisfy conditions $\bC_{4-6} \vee \tilde{
\bC}_{4-6}$ at at most $\ka$ interpolation nodes.
Furthermore, if $\varphi$ meets conditions
$\bC_{4-6} \vee \tilde{\bC}_{4-6}$ at
exactly $k \leq \ka$ interpolation nodes, then the function
$\w = \mathbf{T}_\Theta[\varphi]$ belongs to the class
$\tN_{\ka - k}$.
\label{Thm2}
\end{Tm}

\begin{proof}
Set $p = n - k$.
Without loss of generality, we assume that
$\varphi \in \tN_0$ satisfies conditions $\bC_{1-3} \vee
\tilde{\bC}_{1-3}$ at
interpolation nodes $x_1, \ldots, x_{p}$, conditions
$\bC_{4-6}$ at $x_{p + 1}, \ldots, x_{q}$, and
conditions $\tilde{\bC}_{4-6}$ at $x_{q + 1}, \ldots,
x_n$. 

By Theorem \ref{Thm1}, the block $\tilde{P}_{22}$ is negative
semidefinite and we have ${\rm sq}_-P_{11} = \ka - k$.
Since $\varphi$ meets conditions $\bC_{1-3} \vee
\tilde{\bC}_{1-3}$ at $x_1, \ldots, x_{n - k}$,
it follows from Remark \ref{Russia} that
the function $\w = \mathbf{T}_\Theta[\varphi]$ is such that 
the kernel $K_\w$ has
at least $\ka - k$ negative squares.

Let us start by making the following technical assumption
(similar to (\ref{technical})):
\begin{equation}
\tilde{p}_{ii} \neq 0 \quad\quad (i = \mathbb{N}_n \setminus
\mathbb{N}_q).
\label{technical2}
\end{equation}
Then we have, as in Theorem \ref{Thm1},
a truncated boundary Nevanlinna-Pick
interpolation problem with interpolation nodes $x_{p+1}, \ldots, x_n$.
The Pick matrix $\mathbb{P}$ constructed as in Theorem \ref{Thm1}
is positive semidefinite, and so the interpolation problem
has at least one solution $\varphi$ (see Section 7.3).
We split the rest of the proof into two cases: $(1)\; \tilde{P}_{22} < 0$ is 
invertible and $(2)\; \tilde{P}_{22} \leq 0$ is singular.  Note
that, by the same reasoning used in Theorem \ref{Thm1}, case $(1)$
cannot occur if $\tilde{p}_{ii} = 0$ for some $i \in \mathbb{N}_n
\setminus \mathbb{N}_q$.

{\bf Case I: }Here $\mathbb{P} > 0$ and, by Theorem \ref{Canada},
$\varphi$ admits a
representation
\begin{equation*}
\varphi = \mathbf{T}_{\hat{\Theta}}[\hat{\varphi}]
\end{equation*}
for some $\hat{\varphi} \in \tN_0$, where $\hat{\Theta}$
is given by
\begin{equation}
\hat{\Theta}(z) = I_2 - i\begin{bmatrix}M \\ E_2\end{bmatrix}
(zI - X_2)^{-1}\mathbb{P}^{-1}
\begin{bmatrix}M^* & E_2^*\end{bmatrix}J
\label{bolstar}
\end{equation}
with
\begin{align*}
M &= \begin{bmatrix}\eta_{p+1} & \ldots & \eta_{q} &
-\frac{\tilde{p}_{q+1,q+1}}{\tilde{c}_{q+1}} & \ldots &
-\frac{\tilde{p}_{nn}}{\tilde{c}_n}\end{bmatrix},\\
E_2 &= \begin{bmatrix}1 & \ldots & 1 & 0 & \ldots & 0\end{bmatrix},
\end{align*}
which play the role of $C$ and $E$ in Theorem \ref{Canada}.
With $\bC$ as in (\ref{CP2}), it holds that
\begin{equation*}
\begin{bmatrix}M \\ E_2\end{bmatrix} (zI-X_2)^{-1}
\bC =
\begin{bmatrix}\tC_2 \\ \tE_2\end{bmatrix}
(zI - X_2)^{-1}
\end{equation*}
since
\begin{align*}
\begin{bmatrix}\eta_i \\ 1\end{bmatrix} \cdot
\frac{1}{z - x_i} \cdot \tilde{e}_i &=
\begin{bmatrix}\tilde{c}_i \\ \tilde{e}_i\end{bmatrix}
\cdot \frac{1}{z - x_i} 
\quad& (i \in \mathbb{N}_q \setminus \mathbb{N}_p),\\
\begin{bmatrix} -\frac{\tilde{p}_{jj}}{\tilde{c}_j} \\
0\end{bmatrix} \cdot \frac{1}{z - x_j} \cdot
-\frac{\tilde{c}_j^2}{\tilde{p}_{jj}} &=
\begin{bmatrix}\tilde{c}_j \\ 0\end{bmatrix}
\cdot \frac{1}{z - x_j}
\quad &(j \in \mathbb{N}_n \setminus \mathbb{N}_q).
\end{align*}
Then (\ref{bolstar}) may be written as
\begin{align*}
\hat{\Theta}(z) &= I_2 - i\begin{bmatrix}M \\ E_2\end{bmatrix}
(zI - X_2)^{-1}\bC \bC^{-1} \mathbb{P}^{-1}
\bC^{-*}\bC^*
\begin{bmatrix}M^* & E_2^*\end{bmatrix}J\\
&= I_2 + i\begin{bmatrix}\tC_2 \\ \tE_2\end{bmatrix}
(zI - X_2)^{-1}\tilde{P}_{22} 
\begin{bmatrix}\tC_2^* & \tE_2^*\end{bmatrix}J,
\end{align*}
and by the reflection principle, we have
\begin{equation*}
\hat{\Theta}(z)^{-1} = I_2 - i
\begin{bmatrix}\tC_2 \\ \tE_2\end{bmatrix}
\tilde{P}_{22}(zI-X_2)^{-1}
\begin{bmatrix}\tC_2^* & \tE_2^*\end{bmatrix}J,
\end{equation*}
which coincides with $\tilde{\Theta}^{(2)}$ introduced in
(\ref{e3.24}).
It follows from Lemma \ref{InvLemma} that
\begin{equation}
\Theta(z) = \Theta^{(1)}(z)\hat{\Theta}(z)^{-1},
\label{T8Fact}
\end{equation}
where $\Theta^{(1)}$ is given by
\begin{equation*}
\Theta^{(1)}(z) = I_2 -i
\begin{bmatrix}\tC_1 \\ \tE_1\end{bmatrix}
(zI - X_1)^{-1}P_{11}^{-1}
\begin{bmatrix}\tC_1^* & \tE_1^*\end{bmatrix}J.
\end{equation*}
By Lemma \ref{InvLemma},
\begin{equation}
\Theta^{(1)} \in \mathcal{W}_{\ka_1} \quad\text{where}\quad
\ka_1 = {\rm sq}_-P_{11} = \ka - k.
\label{T1Fact}
\end{equation}
Due to (\ref{T8Fact}), the following chain of equalities holds:
\begin{equation*}
\w := \mathbf{T}_\Theta[\varphi] = \mathbf{T}\left[
\mathbf{T}_{\hat{\Theta}}[\hat{\varphi}]\right] =
\mathbf{T}_{\Theta \hat{\Theta}}[\hat{\varphi}] =
\mathbf{T}_{\Theta^{(1)}}[\hat{\varphi}].
\end{equation*}
This implies, in conjuction with (\ref{T1Fact})
and the fact that $\varphi \in \tN_0$, that
the kernel
$K_\w$ has at most $\ka - k$ negative squares.
Since we have already shown that the kernel
$K_\w$ has at least $\ka - k$ negative squares,
the proof is complete for the first case.

{\bf Case II: }Here $\mathbb{P}$ is positive definite
and singular, and so, by Theorem
\ref{T7.8}, $\varphi$ admits a representation
\begin{equation*}
\varphi(z) = \frac{x^*(zI - X_2)^{-1}M^*}{x^*(zI - X_2)^{-1}E_2^*}
\end{equation*}
where $x$ is any nonzero vector such that $\mathbb{P}x = 0$.
Set $y = \bC^{-1}x$ and note that $\tilde{P}_{22}y = 0$.
We may express $\varphi$ alternately:
\begin{align*}
\varphi(z) &= \frac{y^*\bC^*(zI - X_2)^{-1}M^*}{
y^*\bC^*(zI - X_2)^{-1}E_2^*} \\
&= \frac{y^*(zI - X_2)^{-1}\tC_2^*}{
y^*(zI - X_2)^{-1}\tE_2^*}.
\end{align*}
Since $\varphi(\bar{z}) = \varphi(z)^*$, we also have
\begin{equation*}
\varphi(z) = \frac{\tC_2(zI - X_2)^{-1}y}{
\tE_2(zI-X_2)^{-1}y}.
\end{equation*}
Set
\begin{equation*}
u(z) = \tC_2(zI - X_2)^{-1}y
\quad\text{and}\quad
v(z) = \tE_2(zI - X_2)^{-1}y,
\end{equation*}
and note that
\begin{align*}
u(z)v(\zeta)^* - v(z)u(\zeta)^* =&
y^*(\bar{\zeta}I - X_2)^{-1}\tE_2^*\tC_2
(zI-X_2)^{-1}y \\
& - y^*(\bar{\zeta}I - X_2)^{-1}
\tC_2^*\tE_2(zI - X_2)^{-1}y\\
=& y^*(\bar{\zeta}I-X_2)^{-1}\left(
\tE_2^*\tC_2 - \tC_2^*\tE_2\right)
(zI-X_2)^{-1}y\\
=& y^*(\bar{\zeta}I-X_2)^{-1}(X_2\tilde{P}_{22} -
\tilde{P}_{22}X_2)(zI-X_2)^{-1}y\\
=& y^*(\bar{\zeta}I - X_2)^{-1}\left(-\bar{\zeta}\tilde{P}_{22}
+ X_2\tilde{P}_{22} - \tilde{P}_{22}X_2 + z\tilde{P}_{22} \right.\\
&+ \left. \bar{\zeta}\tilde{P}_{22} - z\tilde{P}_{22}\right)
(zI - X_2)^{-1}y\\
=& -(z - \bar{\zeta})y^*(\bar{\zeta}I - X_2)^{-1}
\tilde{P}_{22}(zI-X_2)^{-1}y,
\end{align*}
so that
\begin{align*}
\frac{\varphi(z) - \varphi(\zeta)^*}{z - \bar{\zeta}} &=
\frac{1}{z - \bar{\zeta}} \cdot
\frac{u(z)v(\zeta)^* - v(z)u(\zeta)^*}{v(\zeta)^*v(z)}\\
&= - \frac{y^*}{v(\zeta)^*} \cdot (\bar{\zeta}I - X_2)^{-1}
\tilde{P}_{22}(zI - X_2)^{-1} \cdot \frac{y}{v(z)}.
\end{align*}
We also have for $\Psi_\varphi$ defined in (\ref{r6.4})
\begin{align}
\Psi_\varphi(z) &= (zI - X)^{-1}\left(\tC^* -
  \varphi(z)\tE^*\right)\nonumber \\
&= (zI - X)^{-1}\left(\tC^*v(z) - u(z)\tE^*\right) \cdot
\frac{1}{v(z)} \nonumber \\
&= (zI - X)^{-1}\left(\tC^*\tE_2 -
  \tE^*\tC_2\right)(zI - X_2)^{-1} \cdot
\frac{y}{v(z)}.
\label{Germany}
\end{align}
Substituting the partitionings for $P^{-1}$, $X$,
$E$, and $C$ into the 
Lyaponov identity (\ref{e3.4}) and comparing the
right block entries yields
\begin{equation*}
X\begin{bmatrix}\tilde{P}_{12} \\ \tilde{P}_{22}\end{bmatrix}
- \tilde{P}\begin{bmatrix}0 \\ X_2\end{bmatrix} =
\tE^*\tC_2 - \tC^*\tE_2.
\end{equation*}
This allows us to rewrite (\ref{Germany}) as
\begin{equation*}
\Psi_\varphi(z) = \begin{bmatrix}\tilde{P}_{12} \\ \tilde{P}_{22}\end{bmatrix}
(zI - X_2)^{-1}\cdot \frac{y}{v(z)} - (zI - X)^{-1}
\begin{bmatrix}\tilde{P}_{12} \\ \tilde{P}_{22}\end{bmatrix}
\cdot \frac{y}{v(z)}.
\end{equation*}
Since $\tilde{P}_{22}y = 0$, this becomes
\begin{equation*}
\Psi_\varphi(z) = \begin{bmatrix}\tilde{P}_{12} \\ \tilde{P}_{22}\end{bmatrix}
(zI - X_2)^{-1}\cdot \frac{y}{v(z)} - (zI - X)^{-1}
\begin{bmatrix}\tilde{P}_{12} \\ 0\end{bmatrix}
\cdot \frac{y}{v(z)}.
\end{equation*}
Using the latter equality, we have,
\begin{align}
\Psi_\varphi(\zeta)^*P\Psi_\varphi(z) =& \frac{y^*}{v(\zeta)^*}\left(
(\bar{\zeta}I - X_2)^{-1}\tilde{P}_{22}(zI -
X_2)^{-1}\right. \nonumber \\
& \; \left. + \begin{bmatrix}\tilde{P}_{12}^* & 0 \end{bmatrix}
(\bar{\zeta}I - X_1)^{-1}P(zI - X)^{-1}\begin{bmatrix}\tilde{P}_{12}
  \\ 0\end{bmatrix}\right) \cdot \frac{y}{v(z)} \nonumber \\
=& \frac{y^*}{v(\zeta)^*}\left((\bar{\zeta}I - X_2)^{-1}\tilde{P}_{22}
(zI - X_2)^{-1} \right. \nonumber \\
&+ \left. \tilde{P}_{21}(\bar{\zeta}I - X_1)^{-1}P_{11}
(zI - X_1)^{-1}\tilde{P}_{12}\right)\cdot \frac{y}{v(z)}.
\label{England}
\end{align}
Substituting (\ref{England})
into (\ref{r6.7}),
we have
\begin{equation*}
K_\w(z,z) = \frac{y^*}{(V_\varphi(\zeta)v(\zeta))^*} \cdot
\tilde{P}_{12}^* (\bar{\zeta}I - X_1)^{-1}P_{11}
(zI - X_1)^{-1}\tilde{P}_{12} \cdot
\frac{y}{V_\varphi(z)v(z)}.
\end{equation*}
This shows that ${\rm sq}_-K_\w \leq {\rm sq}_-P_{11} = \ka - k$.
Since we have already shown that $K_\w$ has at least $\ka - k$
negative squares, the proof for the second case is complete.

It now only remains to show that the technical assumption 
(\ref{technical2}) may be lifted.  As shown in 
Remark \ref{Yes}, when this condition is lifted,
$\w = \mathbf{T}_\Theta[\infty]$.  By (\ref{r6.8}) and
(\ref{r6.7c}),
\begin{equation*}
K_\w(z,z) = \frac{\Psi_\infty(z)^*P\Psi_\infty(z)}{|V_\infty(z)|^2},
\end{equation*}
with
\begin{equation*}
\Psi_\infty(z) = -(zI - X)^{-1}\tE^*.
\end{equation*}
Since $E_2 = 0$ in this situation, we in fact have,
\begin{equation*}
K_\w(z,z) = \frac{1}{|V_\infty(z)|^2}\cdot
(zI - X_1)^{-1}\tE_1^*P_{11}\tE_1
(zI - X_1)^{-1},
\end{equation*}
which shows that ${\rm sq}_-K_\w \leq {\rm sq}_-P_{11} = 
\ka - k$.
We have already shown that $K_\w$ has at least $\ka - k$
negative squares, and so the proof is complete.
\end{proof}

\chapter{The Degenerate Case}
Here we study Problem \ref{Prob3} in the case where
the Pick matrix $P$ is singular.

\begin{La}
Let $A \in \C^{n \times n}$ be a Hermitian matrix with
$\ka$ negative eigenvalues.  Suppose $B$ is a 
$k \times k$ invertible principal submatrix of $A$.
Then there exists a principal submatrix $C$
of the matrix $A$ such that
\begin{enumerate}
\item \text{$B$ is a principal submatrix of $C$}.
\item ${\rm rank}\;C = {\rm rank}\;A.$
\item ${\rm sq}_-C = {\rm sq}_-A.$
\end{enumerate}
\label{La9}
\end{La}
\begin{proof}
Without loss of generality (by permutation similarity), we may assume that $B$
is the upper $k \times k$ principal submatrix of $A$.
Partition $A$ so that
\begin{equation*}
A = \begin{bmatrix}
B & Y \\ Y^* & Z
\end{bmatrix}
\end{equation*}
The matrix $\tilde{A}$ given by
\begin{equation*}
\tilde{A} =
\begin{bmatrix}I_k & 0 \\ -Y^*B^{-1} & I_{n-k}\end{bmatrix}
\begin{bmatrix}B & Y \\ Y^* & Z\end{bmatrix}
\begin{bmatrix}I_k & -B^{-1}Y \\ 0 & I_{n-k}\end{bmatrix}
=
\begin{bmatrix}B & 0 \\ 0 & Z - Y^*B^{-1}Y\end{bmatrix}
\end{equation*}
has the same inertia as does $A$.
Note that
\begin{align*}
s &:= {\rm sq}_-(Z - Y^*B^{-1}Y) =
{\rm sq}_-A - {\rm sq}_-B,\\
r &:= {\rm rank}\;(Z - Y^*B^{-1}Y) =
{\rm rank}\;A - {\rm rank}\;B.
\end{align*}
It is possible to choose a principal
submatrix
$W$ of $Z - Y^*B^{-1}Y$ so that
\begin{equation*}
{\rm sq}_-W = s
\quad\text{and}\quad
{\rm rank}\;W = r.
\end{equation*}
Let $m$ denote the size of the square
matrix $W$
Without loss of generality, we may
assume that $W$ is the upper $m \times m$
principal submatrix of $Z - Y^*B^{-1}Y$.
Then the matrix
\begin{equation*}
C = \begin{bmatrix}I_k & 0 \\ Y^*B^{-1} & I_{n-k}\end{bmatrix}
\begin{bmatrix}
B & 0 & 0\\
0 & W & 0\\
0 & 0 & 0_{n - k - m}\end{bmatrix}
\begin{bmatrix}I_k & B^{-1}Y \\ 0 & I_{n-k}\end{bmatrix}
\end{equation*}
is a principal submatrix of $A$ and has all of the
desired properties.
\end{proof}

\begin{Tm}
Suppose that the Pick matrix $P$ is singular with
${\rm rank}\;P = k < n$ and
${\rm sq}_-P = \ka$.  Then there is a unique
generalized Nevanlinna function
$\w$ such that
\begin{equation}
{\rm sq}_-\mathbf{K}_\w(z, \zeta) = \ka.
\label{9.1}
\end{equation}
This $\w$ is the unique solution of
Problem \ref{Prob3}, and satisfies the interpolation
conditions stated in Problem \ref{Prob1} at at least
$n - {\rm rank}\;P$ interpolation nodes.
\end{Tm}
\begin{proof}
Without loss of generality we may rearrange the interpolation
nodes $x_i$ so that the upper
$k \times k$ principal submatrix $P_{11}$ of $P$
is invertible and has $\ka$ negative eigenvalues.
Due to Lemma \ref{La9} and the fact
that each $\xi_i \neq 0$ (the matrix with all $-\xi_i$'s along
the diagonal corresponds to $B$ in Lemma \ref{La9}), 
we may choose $P_{11}$
so that $\w_i$ and $\gamma_i$ are specified at each node
$x_i$ where $i \in \mathbb{N}_n \setminus \mathbb{N}_k$;
in other words, we ensure that the data corresponding
to all singular interpolation nodes is included $P_{11}$.

We use the conformal partitionings
\begin{equation}
X = \begin{bmatrix}X_1 & 0 \\ 0 & X_2\end{bmatrix}, \quad
E = \begin{bmatrix}E_1 & E_2\end{bmatrix}, \quad
C = \begin{bmatrix}C_1 & C_2\end{bmatrix}
\label{Switzerland}
\end{equation}
and
\begin{equation}
P = \begin{bmatrix}P_{11} & P_{12} \\ P_{21} & P_{22}\end{bmatrix},
\quad {\rm det}\;P_{11} \neq 0, \quad
{\rm sq}_-P_{11} = \ka = {\rm sq}_-P.
\label{Austria}
\end{equation}
Since ${\rm rank}\;P_{11} = {\rm rank}\;P$, the Schur complement
of $P_{11}$ in $P$ is zero; that is,
\begin{equation}
P_{22} = P_{21}P_{11}^{-1}P_{12}.
\label{iszero}
\end{equation}
By substituting (\ref{Switzerland}) and (\ref{Austria})
into the Lyaponov identity (\ref{LI}), multiplying
both sides by $\mathbf{e}_i^*$ on the left, and rearranging,
we obtain the following representation of
the $i$-th row of the block $P_{21}$:
\begin{equation}
\mathbf{e}_i^*P_{21} = (\w_{k + i}E_1 - C_1)(x_{k+i}I - X_1)^{-1}.
\label{9.5}
\end{equation}
We split the rest of the proof into six steps.

\nline
{\bf Step 1:} If $\w$ is a meromorphic function such that (\ref{9.1}) holds,
then it is necessarily of the form
\begin{equation}
\w = \mathbf{T}_{\Theta^{(1)}}[\varphi] :=
\frac{\Theta_{11}^{(1)}\varphi + \Theta_{12}^{(1)}}{
\Theta_{21}^{(1)}\varphi + \Theta_{22}^{(1)}}
\label{9.6}
\end{equation}
for some function $\varphi \in \tN_0$, where
$\Theta^{(1)}$ is given in (\ref{e3.23}).

\nline
{\bf Proof of Step 1:} Write the kernel $K_\w(z, \zeta)$ in the block
form as
\begin{equation}
\mathbf{K}_\w(z, \zeta) = 
\begin{bmatrix}
P_{11} & P_{12} & F_1(z)\\
P_{21} & P_{22} & F_2(z)\\
F_1{\zeta}^* & F_2(\zeta)^* & K_\w(z, \zeta)
\end{bmatrix}
\label{9.7}
\end{equation}
where $F_1$ and $F_2$ are given in (\ref{r5.7}).  The kernel
\begin{align*}
\mathbf{K}_\w^1 &:= 
\begin{bmatrix}
P_{11} & F_1(z)\\
F_1(\zeta)^* & K_\w(z, \zeta)
\end{bmatrix}\\
&= \begin{bmatrix}
P_{11} & (zI - X)^{-1}(\w(z)E^* - C^*) \\
(\w(\zeta)^* E - C)(\bar{\zeta}I - X)^{-1} & K_\w(z, \zeta)
\end{bmatrix}
\end{align*}
is contained in $\mathbf{K}_\w(z, \zeta)$ as a principal submatrix
and therefore ${\rm sq}_-\mathbf{K}_\w^1 \leq \ka$.  On the other
hand, $\mathbf{K}_\w^1$ contains $P_{11}$ as a principal submatrix
and therefore ${\rm sq}_-\mathbf{K}_\w^1 \geq {\rm sq}_-P_{11} = \ka$.
Thus,
\begin{equation}
{\rm sq}_-\mathbf{K}_\w^1 = \ka.
\label{9.8}
\end{equation}
The block $P_{11}$ is an invertible Hermitian matrix with $\ka$ negative
eigenvalues and satisfies the Lyaponov identity in (\ref{LI}).  We may
therefore apply Theorem \ref{lessprecise} to the FMI (\ref{9.8}) to conclude
that $\w$ is of the form (\ref{9.6}) for some $\varphi \in \tN_0$.

\nline
{\bf Step 2:} Every function of the form (\ref{9.6}) solves the
truncated Problem \ref{Prob3}, which is restricted to the
interpolation nodes $x_1, \ldots, x_k$.

\nline
{\bf Proof of Step 2:} The Pick matrix for this truncated
interpolation problem is $P_{11}$, which is invertible and
has $\ka$ negative eigenvalues.  We may then apply
Theorem \ref{lessprecise} (already proven for the nondegenerate case)
to get the desired statement.

\nline
The rational function $\Theta^{(1)}$ is analytic and $J$-unitary
at $x_i$ for every $i \in \mathbb{N}_n \setminus \mathbb{N}_k$.
Let us consider the numbers $a_i$ and $b_i$ defined by
\begin{equation}
\begin{bmatrix}a_i \\ b_i\end{bmatrix} =
\Theta^{(1)}(x_i)^{-1}
\begin{bmatrix}\w_i \\ 1\end{bmatrix}
\quad (i \in \mathbb{N}_n \setminus \mathbb{N}_k) \\
\label{9.9}
\end{equation}
Note that $a_i$ and $b_i$ are real and that
$|a_i| + |b_i| > 0$.

\nline
{\bf Step 3:} It holds that
$
a_ib_j = a_jb_i \quad\text{for} \quad i,j \in \mathbb{N}_n
  \setminus \mathbb{N}_k.
$

\nline
{\bf Proof of Step 3:} Let $i \in \mathbb{N}_n \setminus \mathbb{N}_k$.  
Because of (\ref{9.9}), it holds that
\begin{equation}
\begin{bmatrix}a_i & b_i\end{bmatrix}J
\begin{bmatrix}a_j \\ b_j\end{bmatrix} =
\begin{bmatrix}\w_i & 1\end{bmatrix}
\Theta^{(1)}(x_i)^{-1}
J
\Theta^{(1)}(x_j)^{-1}
\begin{bmatrix}\w_j \\ 1\end{bmatrix}.
\label{9.11}
\end{equation}
Due to Remark \ref{ThetaRemark}, we have
\begin{equation}
\frac{\Theta^{(1)}(\zeta)^{-*}J\Theta^{(1)}(z)^{-1} - J}{
-i(z - \bar{\zeta})} =
J\begin{bmatrix}C_1 \\ E_1\end{bmatrix}
(\bar{\zeta}I - X_1)^{-1}P_{11}^{-1}(zI - X_1)^{-1}
\begin{bmatrix}C_1^* & E_1^*\end{bmatrix}J.
\label{9.12}
\end{equation}
Substituting this into (\ref{9.11}) and taking into account
$\begin{bmatrix}\w_i & 1\end{bmatrix}J\begin{bmatrix}
\w_j \\ 1\end{bmatrix} = i(\w_j - \w_i)$, we obtain
\begin{align}
-a_ib_j + a_jb_i = -\w_i + \w_j - &(x_j - x_i)
(\w_i E_1 - C_1)(x_i I - X_1)^{-1}P_{11}^{-1} \\
&\times (x_j I - X_1)^{-1}(\w_j E_1^* - C_1^*).
\label{NewZealand}
\end{align}
Substituting (\ref{9.5}) into the right hand side of
(\ref{NewZealand},
we get
\begin{align*}  
a_jb_i - a_ib_j = -\w_i + \w_j - (x_j - x_i)
\mathbf{e}_{i-k}^*P_{21}P_{11}^{-1}P_{12}
\mathbf{e}_{j-k}.
\end{align*}
Due to (\ref{iszero}) and
\begin{equation*}
\mathbf{e}_{i-k}P_{22}\mathbf{e}_{j-k} =
\frac{\w_j - \w_i}{x_j - x_i},
\end{equation*}
the right hand side of (\ref{NewZealand}) is zero.
Hence
$a_ib_j = a_jb_i$.  Since $|a_i| + |b_i| > 0$,
$a_i = 0$ implies $a_j = 0$ and $b_i = 0$ implies $b_j = 0$.
So if $b_i \neq 0$, then
\begin{equation*}
\frac{a_i}{b_i} = \frac{a_j}{b_j}.
\end{equation*}

\nline
{\bf Step 4:} The row vectors
\begin{equation*}
A := \begin{bmatrix}a_{k + 1} & \ldots & a_n\end{bmatrix}, \quad
B := \begin{bmatrix}b_{k + 1} & \ldots & b_n\end{bmatrix}
\end{equation*}
can be represented as follows:
\begin{equation}
\begin{bmatrix}A \\ B\end{bmatrix} =
\begin{bmatrix}C \\ E\end{bmatrix}
\begin{bmatrix}-P_{11}^{-1}P_{12} \\ I\end{bmatrix}.
\label{9.15}
\end{equation}

\nline
{\bf Proof of Step 4:} We first substitute formula (\ref{E423}) for
the inverse of $\Theta^{(1)}$ into (\ref{9.9}) to get
\begin{equation}
\begin{bmatrix}a_i \\ b_i\end{bmatrix} =
\begin{bmatrix}\w_i \\ 1\end{bmatrix} -
\begin{bmatrix}C_1 \\ E_1\end{bmatrix}
P_{11}^{-1}(x_iI - X_1)^{-1}(w_iE_1^* - C_1^*)
\quad (i \in \mathbb{N}_n \setminus \mathbb{N}_k).
\label{Australia}
\end{equation}
The equality
\begin{equation}
P_{12}\mathbf{e}_{i} = (x_{k+i}I - X_1)^{-1}(w_{k+i}E_1^* - C_1^*)
\label{Italy}
\end{equation}
follows from (\ref{9.5}).
Substituting (\ref{Italy}) into (\ref{Australia}) yields
\begin{align*}
\begin{bmatrix}A \\ B\end{bmatrix}\mathbf{e}_i =
\begin{bmatrix}\w_{k+i} \\ 1\end{bmatrix} -
\begin{bmatrix}C_1 \\ E_1\end{bmatrix}P_{11}^{-1}P_{12}\mathbf{e}_i
&= \begin{bmatrix}C_2 \\ E_2\end{bmatrix}\mathbf{e}_i -
\begin{bmatrix}C_1 \\ E_1\end{bmatrix}P_{11}^{-1}P_{12}\mathbf{e}_i\\
&= \begin{bmatrix}C \\ E\end{bmatrix}
\begin{bmatrix} -P_{11}^{-1}P_{12} \\ I\end{bmatrix}\mathbf{e}_i.
\end{align*}
Since the latter equality holds for $i = 1, \ldots, n-k$,
(\ref{9.15}) follows.

\begin{Rk}
Comparing (\ref{9.15}) and (\ref{e3.36}), we conclude that
\begin{equation*}
\begin{bmatrix}A \\ B\end{bmatrix} = \Theta^{(1)}(z)^{-1}
\begin{bmatrix}C \\ E\end{bmatrix}(zI - X)^{-1}
\begin{bmatrix}-P_{11}^{-1}P_{12} \\ I\end{bmatrix}
(zI - X_2).
\end{equation*}
\end{Rk}
Using the symmetry relation
$\Theta^{(1)}(z)^{-1} = J\Theta^{(1)}(\bar{z})J$
in the above equality shows that
\begin{equation*}
J\begin{bmatrix}A \\ B\end{bmatrix}(zI - X_2)^{-1} =
\Theta^{(1)}(\bar{z})J
\begin{bmatrix}C \\ E\end{bmatrix}(zI - X)^{-1}
\begin{bmatrix}-P_{11}^{-1}P_{12}\\ I\end{bmatrix}.
\end{equation*}
Taking adjoints and replacing $\bar{z}$ by $z$, we
obtain
\begin{equation}
(zI - X_2)^{-1} \begin{bmatrix}A^* & B^*\end{bmatrix} J =
\begin{bmatrix}-P_{11}^{-1}P_{12} & I\end{bmatrix}
(zI - X)^{-1}\begin{bmatrix}C^* & E^*\end{bmatrix}J
\Theta^{(1)}(z).
\label{9.23}
\end{equation}

\nline
{\bf Step 5:} A function $\w$ of the form (\ref{9.6}) satisfying
FMI (\ref{9.1}) is necessarily such that the corresponding
parameter $\varphi$ is a real constant (or
$\infty$):
\begin{equation}
\varphi(z) \equiv \varphi_0 := \frac{a_{k + 1}}{b_{k + 1}} =
\dots = \frac{a_n}{b_n}.
\label{9.24}
\end{equation}

\nline
{\bf Proof of Step 5:} Let us consider the Schur complement
$\mathbf{S}$ of the block $P_{11}$ in (\ref{9.7}):
\begin{equation*}
\mathbf{S}(z, \zeta) =
\begin{bmatrix}P_{22} & F_2(z) \\ F_2(\zeta)^* & K_\w(z, \zeta)
\end{bmatrix} -
\begin{bmatrix}P_{21} \\ F_1(\zeta)^*\end{bmatrix}
P_{11}^{-1}
\begin{bmatrix}P_{12} & F_1(z)\end{bmatrix}.
\end{equation*}
Since
\begin{equation*}
{\rm sq}_-\mathbf{K}_\w = {\rm sq}_-P_{11} + 
{\rm sq}_-\mathbf{S} = \ka + {\rm sq}_-\mathbf{S},
\end{equation*}
it follows that FMI (\ref{9.1}) is equivalent to positivity
of $\mathbf{S}$ on $\rho(\w)$:
\begin{equation*}
\mathbf{S}(z, \zeta) \succeq 0.
\end{equation*}
Since the ``11'' block in $\mathbf{S}(z, \zeta)$ equals
$P_{22} - P_{21}P_{11}^{-1}P_{12} = 0$, positivity
guarantees that the offdiagonal entries in $\mathbf{S}$
vanish everywhere on $\C^+$:
\begin{equation*}
F_2(z) - P_{21}P_{11}^{-1}F_1(z) \equiv 0.
\end{equation*}
This may be rewritten as
\begin{equation}
\begin{bmatrix}-P_{21}P_{11}^{-1} & I\end{bmatrix}
(zI - X)^{-1}(\w(z)E^* - C^*) \equiv 0.
\label{9.29}
\end{equation}
We shall now show that (\ref{9.29}) holds for $\w$ of the form
(\ref{9.6}) if and only if $\varphi$ is subject to
\begin{equation}
\varphi(z)B \equiv A.
\label{9.30}
\end{equation}
We first assume that 
\begin{equation}
B \neq 0.
\label{Bnot0}
\end{equation}
For $\w$ of the form (\ref{9.6}), it holds that
\begin{equation*}
\begin{bmatrix}C^* & E^*\end{bmatrix}J
\begin{bmatrix}\w(z) \\ 1\end{bmatrix} =
\left(\Theta_{21}^{(1)}\varphi + \Theta_{22}^{(1)}\right)^{-1}
\begin{bmatrix}C^* & E^*\end{bmatrix} J
\begin{bmatrix}\Theta_{11}^{(1)} & \Theta_{12}^{(1)} \\
\Theta_{21}^{(1)} & \Theta_{22}^{(1)}\end{bmatrix}
\begin{bmatrix}\varphi \\ 1\end{bmatrix}
\end{equation*}
and therefore identity (\ref{9.29}) may be written equivalently
in terms of the parameter $\varphi$ as
\begin{equation*}
\begin{bmatrix}-P_{21}P_{11}^{-1} & I\end{bmatrix}
(zI - X)^{-1}\begin{bmatrix}C^* & E^*\end{bmatrix}J
\Theta^{(1)}(z)
\begin{bmatrix}\varphi(z) \\ 1\end{bmatrix} \equiv 0,
\end{equation*}
which, due to (\ref{9.23}), is the same as
\begin{equation*}
(zI - X_1)^{-1}\begin{bmatrix}A^* & B^*\end{bmatrix}J
\begin{bmatrix}\varphi(z) \\ 1\end{bmatrix} \equiv 0.
\end{equation*}
This identity is easily seen to be equivalent to
(\ref{9.30}).  Writing (\ref{9.30}) entrywise, we obtain
the system
\begin{equation}
b_i\varphi(z) \equiv a_i \quad\quad (i = k + 1, \ldots, n).
\label{system}
\end{equation}
By (\ref{Bnot0}), Step 3, and the fact that
$|a_i| + |b_i| > 0$, it follows that
each $b_i \neq 0$, which shows that the system (\ref{system})
admits at most one solution.  Invoking Step 3 again proves
that (\ref{system}) is consistent, and hence admits a solution.
Combining Steps 1 and 5 thus shows
that, for $B \neq 0$, FMI (\ref{9.1}) has at most one solution, the
only candidate being
\begin{equation}
\w = \mathbf{T}_{\Theta^{(1)}}[\varphi_0]
\label{9.35}
\end{equation}
where $\varphi_0$ is the constant defined in (\ref{9.24}).
If $B = 0$, then by interpreting the above expressions projectively
(i.e., take $\begin{bmatrix}1 & 0\end{bmatrix}$ instead of
$\begin{bmatrix}\varphi(z)^* & 1\end{bmatrix}$), we obtain that
(\ref{9.29}) holds if and only if $\varphi(z) \equiv \varphi_0 \equiv
\infty.$  Therefore, without restriction, it holds that FMI
(\ref{9.1}) has at most one solution and that the only candidate is
(\ref{9.35}).

\nline
{\bf Step 6:} The function (\ref{9.35}) satisfies FMI (\ref{9.1}) and
interpolation conditions
\begin{equation}
\w(x_i) = \w_i \quad\text{and}\quad \w^\prime(x_i) = \gamma_i
\quad (i = k+1, \ldots, n).
\label{Step6Con}
\end{equation}

\nline
{\bf Proof of Step 6:} In the proof of Step 5 the equivalence of
equations (\ref{9.30}) and (\ref{9.29}) was established, which implies that
for the function (\ref{9.35}), it holds that
\begin{equation*}
F_2(z) \equiv P_{21}P_{11}^{-1}F_1(z),
\end{equation*}
which, in view of (\ref{9.29}), may be rewritten as
\begin{equation}
(zI - X_2)^{-1}(\w(z)E_2^* - C_2^*) \equiv
P_{21}P_{11}^{-1}(zI - X_1)^{-1}(\w(z)E_1^* - C_1^*).
\label{9.38}
\end{equation}
Our next object is to show that
\begin{equation}
K_\w(z, \zeta) \equiv F_1(\zeta)^* P_{11}^{-1} F_1(z),
\label{9.39}
\end{equation}
or equivalently,
\begin{equation}
\frac{\w(z) - \w(\zeta)^*}{z - \bar{\zeta}} \equiv
(\w(\zeta)^*E_1 - C_1)(\bar{\zeta}I - X_1)^{-1}
P_{11}^{-1}
(zI - X_1)^{-1}(\w(z)E_1^* - C_1^*).
\label{9.40}
\end{equation}
On account of (\ref{9.12}),
\begin{align}
&(\w(z)E^* - C^*)(\bar{\zeta}I - X_1)^{-1}
P_{11}^{-1}(zI - X_1)^{-1}(\w(z)E_1^* - C_1^*) \nonumber\\
& \; = 
\begin{bmatrix}\w(\zeta)^* & 1\end{bmatrix}
\frac{\Theta^{(1)}(\zeta)^{-*}J\Theta^{(1)}(z)^{-1} - J}{-i(z - \bar{\zeta})}
\begin{bmatrix}\w(z) \\ 1\end{bmatrix} \nonumber\\
& \; =\frac{\w(z) - \w(\zeta)^*}{z - \bar{\zeta}} +
\begin{bmatrix}\w(\zeta)^* & 1\end{bmatrix}
\frac{\Theta^{(1)}(\zeta)^{-*}J\Theta^{(1)}(z)^{-1}}{
-i(z - \bar{\zeta})}
\begin{bmatrix}\w(z) \\ 1\end{bmatrix}.
\label{9.43}
\end{align}

If $\varphi_0 \neq \infty$, then representation (\ref{9.35}) is
equivalent to
\begin{equation*}
\begin{bmatrix}\w(z) \\ 1\end{bmatrix} = \Theta^{(1)}(z)
\begin{bmatrix}\varphi_0 \\ 1\end{bmatrix}
\frac{1}{v(z)}, \quad \text{where}\;
v(z) = \Theta_{21}^{(1)}(z)\varphi_0 + \Theta_{22}^{(1)}(z),
\end{equation*}
and therefore
\begin{equation}
\begin{bmatrix}\w(\zeta)^* & 1\end{bmatrix}
\Theta^{(1)}(\zeta)^{-*}J\Theta^{(1)}(z)^{-1}
\begin{bmatrix}\w(z) \\ 1\end{bmatrix} =
\frac{0}{v(\zeta)^*v(z)} \equiv 0.
\label{9.45}
\end{equation}
Using this in (\ref{9.43}) proves (\ref{9.39}).
If $\varphi_0 = \infty$, then representation (\ref{9.35})
is equivalent to
\begin{equation*}
\begin{bmatrix}\w(z) \\ 1\end{bmatrix} =
\Theta^{(1)}(z) \begin{bmatrix}1 \\ 0\end{bmatrix}
\frac{1}{\Theta_{21}^{(1)}(z)},
\end{equation*}
which implies (\ref{9.45}) with $v(z) = \Theta_{21}^{(1)}(z)$, thus
showing that (\ref{9.39}) also holds for $\varphi_0 = \infty$.

Note that the kernel $\mathbf{K}_\w(z, \zeta)$ now admits
the representation
\begin{equation*}
\mathbf{K}_\w(z, \zeta) =
\begin{bmatrix}P_{11} \\ P_{21} \\ F_1(\zeta)^*\end{bmatrix}
P_{11}^{-1}
\begin{bmatrix}P_{11} & P_{12} & F_1(z)\end{bmatrix}.
\end{equation*}

Let $i \in \mathbb{N}_n \setminus \mathbb{N}_k$ and 
compare residues of both sides of identity (\ref{9.38}) at
$z - x_i$:
\begin{equation*}
\mathbf{e}_i\mathbf{e}_i^*(\w(x_i)E_2^* - C_2^*) = 0,
\end{equation*}
which is equivalent to
\begin{equation*}
\w(x_i) - \w_i = 0.
\end{equation*}
Now let $z, \zeta \to x_i$ in (\ref{9.40}) bearing in mind
that $\w(x_i) = \w_i$:
\begin{equation*}
\lim_{z \to x_i}K_\w(z,z) = (\w_iE_1 - C_1)(x_iI - X_1)^{-1}
P_{11}^{-1}(x_iI - X_1)^{-1}(\w_iE_1^* - C_1^*),
\end{equation*}
Comparing $i$-th diagonal entries in (\ref{iszero})
with the last equation implies the second
statement in (\ref{Step6Con}).
\end{proof}

\chapter{Examples}
We present some numerical examples to illustrate the methods
employed in this thesis.

\section{Interpolation With Two Regular Interpolation Nodes}

Let
\begin{equation*}
x_1 = 0, \quad x_2 = 1, \quad
\w_1 = 0, \quad \w_2 = 1, \quad
\gamma_1 = -1, \quad \gamma_2 = 1.
\end{equation*}
Then, according to (\ref{r2.2}), (\ref{r2.3}), and (\ref{1.19}) -- (\ref{120}),
\begin{equation*}
X = \begin{bmatrix}0 & 0 \\ 0 & 1\end{bmatrix}, \quad
E = \begin{bmatrix}1 & 1\end{bmatrix}, \quad
C = \begin{bmatrix}0 & 1\end{bmatrix},
\end{equation*}
\begin{equation*}
P = \begin{bmatrix}-1 & 1 \\ 1 & 1\end{bmatrix},
\quad\text{and}\quad
P^{-1} = \frac{1}{2}\begin{bmatrix}-1 & 1 \\ 1 & 1\end{bmatrix}.
\end{equation*}
Note that $\ka := {\rm sq}_-P = 1$.
Problems \ref{Prob1}, \ref{Prob2}, and \ref{Prob3} thus take
the following forms.

\nline
{\bf Problem \ref{Prob1}:} Find all functions $\w \in \N_1$
such that
\begin{equation*}
\w(0) = 0, \quad \w^\prime(0) = -1, \quad
\w(1) = 1, \quad \w^\prime(1) = 1.
\end{equation*}

\nline
{\bf Problem \ref{Prob2}:} Find all functions $\w \in \N_1$
that satisfy
\begin{equation}
\w(0) = 0, \quad -\infty < \w^\prime(0) \leq -1, \quad
\w(1) = 1, \quad -\infty < \w^\prime(1) \leq 1.
\label{10P2con}
\end{equation}

\nline
{\bf Problem \ref{Prob3}:} Find all functions $\w$ such that
either
\begin{enumerate}
\item $\w \in \N_1$ and satisfies all of the conditions in
(\ref{10P2con}) or
\item $\w \in \N_0$ and satisfies the first two conditions
in (\ref{10P2con}) or
\item $\w \in \N_0$ and satisfies the last two conditions
in (\ref{10P2con}).
\end{enumerate}

By formula (\ref{ThetaFormula2}) for $\Theta$, we have
\begin{align*}
\Theta(z) &=
I_2 - i \begin{bmatrix}0 & 1\\1 & 1\end{bmatrix}
\begin{bmatrix} \frac{1}{z} & 0 \\ 0 & \frac{1}{z - 1}\end{bmatrix}
\begin{bmatrix}-\frac{1}{2} & \frac{1}{2} \\ \frac{1}{2} &
  \frac{1}{2}\end{bmatrix}
\begin{bmatrix}0 & 1 \\ 1 & 1\end{bmatrix}
\begin{bmatrix}0 & -i \\ i & 0\end{bmatrix} \\
&= \frac{1}{2z(z-1)}
\begin{bmatrix}2z^2 & -z\\ 2z & 2z^2 - 4z + 1\end{bmatrix},
\end{align*}
and so, by Theorem \ref{lessprecise}, 
the set of all solutions $\w$ of
Problem \ref{Prob3} is parameterized by the linear
fractional formula
\begin{equation}
\w(z) = \frac{2z^2 \cdot \varphi(z) - z}{
2z \cdot \varphi(z) + 2z^2 - 4z + 1},
\label{Ex1form}
\end{equation}
where the parameter $\varphi$ runs through the extended
Nevanlinna class $\tN_0$.
By (\ref{e3.2}), we have
\begin{equation*}
\tE = \begin{bmatrix}0 & 1\end{bmatrix}, \quad
\tC = \begin{bmatrix}\frac{1}{2} & \frac{1}{2} \end{bmatrix}.
\end{equation*}
The diagonal entries of $P^{-1}$ are 
$\tilde{p}_{11} = -1/2$ and $\tilde{p}_{22} = 1/2$, and so
\begin{equation*}
\eta_1 := \frac{\tilde{c}_1}{\tilde{e}_1} = \infty, \quad
\eta_2 := \frac{\tilde{c}_2}{\tilde{e}_2} = \frac{1}{2}, \quad
\frac{\tilde{p}_{11}}{\tilde{c}_1^2} = -2, \quad
\frac{\tilde{p}_{22}}{\tilde{e}_2^2} = \frac{1}{2}.
\end{equation*}
It follows from Theorem \ref{Tm27} that every function $\w$ of the form
(\ref{Ex1form}) also solves Problem \ref{Prob2}, unless the
parameter $\varphi$ is subject to
\begin{equation}
\varphi(0) = \infty \quad\text{and}\quad
-\frac{1}{\varphi_{-1}(0)} \leq 2
\label{subject1}
\end{equation}
or to
\begin{equation}
\varphi(1) = \frac{1}{2} \quad\text{and}\quad
\varphi^\prime(1) \leq -\frac{1}{2}. 
\label{subject2}
\end{equation}
Note that a parameter $\varphi$ can never satisfy condition
(\ref{subject2}), due to Theorems \ref{CJN1} and \ref{CJGN1}.
It follows from Theorem \ref{Tm26} that
every function $\w$ of the form (\ref{Ex1form}) solves Problem
\ref{Prob1}, unless the parameter $\varphi$ is subject to
\begin{equation*}
\varphi(0) = \infty \quad\text{and}\quad
-\frac{1}{\varphi_{-1}(0)} < \infty
\end{equation*}
or to
\begin{equation*}
\varphi(1) = \frac{1}{2} \quad\text{and}\quad
\varphi^\prime(1) < \infty.
\end{equation*}
We note in particular that every parameter $\varphi \in \tN_0$
that satisfies conditions (\ref{subject1}) or (\ref{subject2}) leads
to a solution $\w$ of Problem \ref{Prob3} that is not a solution
of Problem \ref{Prob2}.

\section{Interpolation With One Regular and One Singular Interpolation
Node}
Let
\begin{equation*}
x_1 = 1, \quad x_2 = 0, \quad
\w_1 = 0, \quad \gamma_1 = -1, \quad
\xi_2 = -1.
\end{equation*}
Then, according to (\ref{r2.2}), (\ref{r2.3}), and (\ref{1.19}) -- (\ref{120}),
\begin{equation*}
X = \begin{bmatrix}1 & 0 \\ 0 & 0\end{bmatrix}, \quad
E = \begin{bmatrix}1 & 0\end{bmatrix}, \quad
C = \begin{bmatrix}0 & -1\end{bmatrix},
\end{equation*}
\begin{equation*}
P = \begin{bmatrix}-1 & 1 \\ 1 & 1\end{bmatrix},
\quad\text{and}\quad
P^{-1} = \frac{1}{2}\begin{bmatrix}-1 & 1 \\ 1 & 1\end{bmatrix}.
\end{equation*}
Note that the Pick matrix $P$ is the same for this data
as it was for the data used in the previous example, even
though the two data sets are quite different.
As before,
$\ka := {\rm sq}_-P = 1$.
Problems \ref{Prob1}, \ref{Prob2}, and \ref{Prob3} thus take
the following forms.

\nline
{\bf Problem \ref{Prob1}:} Find all functions $\w \in \N_1$
such that
\begin{equation*}
\w(1) = 0, \quad \w^\prime(1) = -1, \quad
\w_{-1}(0) = -1.
\end{equation*}

\nline
{\bf Problem \ref{Prob2}:} Find all functions $\w \in \N_1$
that satisfy
\begin{equation}
\w(1) = 0, \quad -\infty < \w^\prime(1) \leq -1, \quad
-\infty < -\frac{1}{\w_{-1}(0)} \leq 1.
\label{10P2con2}
\end{equation}

\nline
{\bf Problem \ref{Prob3}:} Find all functions $\w$ such that
either
\begin{enumerate}
\item $\w \in \N_1$ and satisfies all of the conditions in
(\ref{10P2con2}) or
\item $\w \in \N_0$ and satisfies the first two conditions
in (\ref{10P2con2}) or
\item $\w \in \N_0$ and satisfies the last condition
in (\ref{10P2con2}).
\end{enumerate}

By formula (\ref{ThetaFormula2}) for $\Theta$, we have
\begin{align*}
\Theta(z) &=
I_2 - i \begin{bmatrix}0 & -1\\1 & 0\end{bmatrix}
\begin{bmatrix} \frac{1}{z-1} & 0 \\ 0 & \frac{1}{z}\end{bmatrix}
\begin{bmatrix}-\frac{1}{2} & \frac{1}{2} \\ \frac{1}{2} &
  \frac{1}{2}\end{bmatrix}
\begin{bmatrix}0 & 1 \\ -1 & 0\end{bmatrix}
\begin{bmatrix}0 & -i \\ i & 0\end{bmatrix} \\
&= \frac{1}{2z(z-1)} \begin{bmatrix} 2z^2 - 3z + 1 & 1-z \\
-z & 2z^2 - z\end{bmatrix},
\end{align*}
and so, by Theorem \ref{lessprecise}, the set of all solutions $\w$ of
Problem \ref{Prob3} is parameterized by the linear
fractional formula
\begin{equation}
\w(z) = \frac{ (-2z^2 + 3z -1)\varphi(z) + z - 1}{ z \varphi(z) - 2z^2 +
z},
\label{Ex1form2}
\end{equation}
where the parameter $\varphi$ runs through the extended
Nevanlinna class $\tN_0$.
By (\ref{e3.2}), we have
\begin{equation*}
\tE = \begin{bmatrix}-\frac{1}{2} & \frac{1}{2} \end{bmatrix}, \quad
\tC = \begin{bmatrix}-\frac{1}{2} & -\frac{1}{2} \end{bmatrix}.
\end{equation*}
The diagonal entries of $P^{-1}$ are 
$\tilde{p}_{11} = -1/2$ and $\tilde{p}_{22} = 1/2$, and so
\begin{equation*}
\eta_1 := \frac{\tilde{c}_1}{\tilde{e}_1} = 1, \quad
\eta_2 := \frac{\tilde{c}_2}{\tilde{e}_2} = -1, \quad
\frac{\tilde{p}_{11}}{\tilde{e}_1^2} = 1, \quad
\frac{\tilde{p}_{22}}{\tilde{e}_2^2} = 1.
\end{equation*}
From  Theorem \ref{Tm27} it follows that every function $\w$ of the form
(\ref{Ex1form}) also solves Problem \ref{Prob2}, unless the
parameter $\varphi$ is subject to
\begin{equation}
\varphi(1) = 1 \quad\text{and}\quad
\varphi^\prime(1) \leq 1
\label{2subject1}
\end{equation}
or to
\begin{equation}
\varphi(0) = -1 \quad\text{and}\quad
\varphi^\prime(0) \leq 1. 
\label{2subject2}
\end{equation}
On the other hand, it follows from Theorem \ref{Tm26} that
every function $\w$ of the form (\ref{Ex1form}) solves Problem
\ref{Prob1}, unless the parameter $\varphi$ is subject to
\begin{equation*}
\varphi(1) = 1 \quad\text{and}\quad
\varphi^\prime(1) < \infty
\end{equation*}
or to
\begin{equation*}
\varphi(0) = -1 \quad\text{and}\quad
\varphi^\prime(0) < \infty. 
\end{equation*}
We note in particular that every parameter $\varphi \in \tN_0$
that satisfies conditions (\ref{2subject1}) or (\ref{2subject2}) leads
to a solution $\w$ of Problem \ref{Prob3} that is not a solution
of Problem \ref{Prob2}.

\section{Two Point Degenerate Interpolation Problem}

Let
\begin{equation*}
x_1 = -\frac{1}{2}, \quad x_2 = \frac{1}{2}, \quad
\w_1 = 0, \quad \gamma_1 = -1, \quad
\xi_2 = 1.
\end{equation*}
Then, according to (\ref{r2.2}), (\ref{r2.3}), and (\ref{1.19}) -- (\ref{120}),
\begin{equation*}
X = \begin{bmatrix}-\frac{1}{2} & 0 \\ 0 & \frac{1}{2}\end{bmatrix}, \quad
E = \begin{bmatrix}1 & 0\end{bmatrix}, \quad
C = \begin{bmatrix}0 & 1\end{bmatrix}, \quad \text{and}
\quad
P = 
\begin{bmatrix}-1 & 1 \\ 1 & -1\end{bmatrix}.
\end{equation*}
Note that $P$ is singular, with spectrum 
$\sigma(P) = \{-2, 0\}$.
With the conditions
\begin{equation}
\w(-1/2) = 0, \quad -\infty < \w^\prime(-1/2) \leq -1, \quad
-\infty < -\frac{1}{\w_{-1}(1/2)} \leq -1,
\label{10P2con3}
\end{equation}
Problem \ref{Prob3} thus takes the form

\nline
{\bf Problem \ref{Prob3}:} Find all functions $\w$ such that
either
\begin{enumerate}
\item $\w \in \N_1$ and satisfies all of the conditions in
(\ref{10P2con3}) or
\item $\w \in \N_0$ and satisfies the first two conditions
in (\ref{10P2con3}) or
\item $\w \in \N_0$ and satisfies the last condition
in (\ref{10P2con3}).
\end{enumerate}
By Theorem \ref{TmDeg}, Problem \ref{Prob3} has a unique solution.
According to Chapter IX, the solution $\w$ is given by
\begin{equation*}
\w(z) := \mathbf{T}_{\Theta^{(1)}}[\varphi_0],
\end{equation*}
where $\Theta^{(1)}(z)$ and $\varphi_0$ are described below.
The matrix-valued function
$\Theta^{(1)}(z)$
is given (after rearranging the two interpolation nodes) by (\ref{e3.23}):
\begin{equation*}
\Theta^{(1)}(z) = 1 - i \begin{bmatrix}1 \\ 0\end{bmatrix}
\cdot \frac{2}{2z - 1} \cdot -1 \cdot
\begin{bmatrix}1 & 0\end{bmatrix}
\begin{bmatrix}0 & -i \\ i & 0\end{bmatrix}
= \begin{bmatrix}1 & \frac{2z+1}{2z-1} \\ 1 & 1\end{bmatrix}.
\end{equation*}
The constant $\varphi_0$ is equal to the the ratio $a/b$, which
we obtain from
\begin{equation*}
\begin{bmatrix}a \\ b\end{bmatrix} =
\Theta^{(1)}(-1/2)\begin{bmatrix}0 \\ 1\end{bmatrix}
= \begin{bmatrix}1 & 0 \\ -1 & 1\end{bmatrix}
\begin{bmatrix}0 \\ 1\end{bmatrix} =
\begin{bmatrix}0 \\ 1\end{bmatrix}.
\end{equation*}
Therefore the unique solution is
\begin{equation*}
\w(z) = \frac{2z + 1}{2z - 1}.
\end{equation*}
It is easy to check that $\w$ satisfies each condition
in (\ref{10P2con3}) and belongs to the class
$\N_1$.

\end{document}